\documentclass[a4paper,12pt,twoside]{article}
\usepackage[latin1]{inputenc}
\usepackage[english]{babel}
\usepackage{amsfonts,amssymb,amsmath,exscale}
\usepackage[dvips]{graphicx,color,psfrag}
\usepackage{fancyhdr}
\usepackage{rotating}
\usepackage{pifont}
\usepackage{hyperref}

\usepackage{natbib}
\usepackage[font=small,labelfont=bf]{caption}

\fancypagestyle{plain}{%
  \fancyhead{}%
  \fancyfoot[c]{\sffamily\thepage}%
}

\newcommand*{\inputfig}[3][htb]{{%
    \def\fps@figure{#1}%
    \def\DIR{#2}%
    \def\LABEL{#3}%
    \graphicspath{{\DIR/}}%
    \begin{figure}
\centering
{
  \footnotesize
  \psfrag{H}     [cc][cc] {$2L$}
  \psfrag{L}     [cc][cc] {$L$}
  \psfrag{u}     [cc][cc] {$\bar u$}
  \psfrag{X}     [cc][cc] {{\scriptsize $X$}}
  \psfrag{Y}     [cc][cc] {{\scriptsize $Y$}}
  \psfrag{B0}     [cc][cc] {$\calB$}
  \includegraphics*[scale=0.8]{fig.eps}
}
\caption {Cross Shear Localization. Geometry and mechanical loading. The process of heat conduction is neglected. Due to the symmetry of the boundary value problem, only the top right quarter of the domain is discretized by finite elements. To trigger plasticity in the center, the initial yield stress $y_0$ is reduced by 3\% in the dark grey element.
}
\label{\LABEL}
\end{figure}
}}%

\newcommand*{\inputtab}[3][htb]{{%
    \def\fps@table{#1}%
    \def\nametable{#2}%
    \def\LABEL{#3}%
    \input{\nametable}
}}%

\def\boldsymboli#1{{%
   \let\@nomath\@gobble\mathversion{bold}%
   \mathchoice%
      {\hbox{$\mathsurround0pt\displaystyle#1$}}%
      {\hbox{$\mathsurround0pt\textstyle#1$}}%
      {\hbox{$\mathsurround0pt\scriptstyle#1$}}%
      {\hbox{$\mathsurround0pt\scriptscriptstyle#1$}}%
}}

\DeclareRobustCommand{\Bbeta      }{{\boldsymbol{\beta}}}
\DeclareRobustCommand{\Bgamma     }{{\boldsymbol{\gamma}}}

\DeclareRobustCommand{\Bsigma     }{{\boldsymbol{\sigma}}}

\DeclareRobustCommand{\Bvarepsilon}{{\boldsymbol{\varepsilon}}}

\DeclareRobustCommand{\Bvarphi    }{{\boldsymbol{\varphi}}}

\DeclareRobustCommand{\BA}{{\boldsymbol{\mathnormal A}}}
\DeclareRobustCommand{\BB}{{\boldsymbol{\mathnormal B}}}
\DeclareRobustCommand{\BC}{{\boldsymbol{\mathnormal C}}}

\DeclareRobustCommand{\BE}{{\boldsymbol{\mathnormal E}}}
\DeclareRobustCommand{\BF}{{\boldsymbol{\mathnormal F}}}
\DeclareRobustCommand{\BG}{{\boldsymbol{\mathnormal G}}}

\DeclareRobustCommand{\BP}{{\boldsymbol{\mathnormal P}}}

\DeclareRobustCommand{\BX}{{\boldsymbol{\mathnormal X}}}

\DeclareRobustCommand{\Ba}{{\boldsymbol{\mathnormal a}}}
\DeclareRobustCommand{\Bb}{{\boldsymbol{\mathnormal b}}}

\DeclareRobustCommand{\Be}{{\boldsymbol{\mathnormal e}}}

\DeclareRobustCommand{\Bg}{{\boldsymbol{\mathnormal g}}}

\DeclareRobustCommand{\Bn}{{\boldsymbol{\mathnormal n}}}

\DeclareRobustCommand{\Bq}{{\boldsymbol{\mathnormal q}}}

\DeclareRobustCommand{\Bs}{{\boldsymbol{\mathnormal s}}}
\DeclareRobustCommand{\Bt}{{\boldsymbol{\mathnormal t}}}

\DeclareRobustCommand{\Bx}{{\boldsymbol{\mathnormal x}}}

\DeclareRobustCommand{\calB}{{\mathcal B}}

\DeclareRobustCommand{\calD}{{\mathcal D}}
\DeclareRobustCommand{\calE}{{\mathcal E}}

\DeclareRobustCommand{\calO}{{\mathcal O}}

\DeclareRobustCommand{\calT}{{\mathcal T}}

\DeclareMathAlphabet{\IKbb}{U}{bbm}{m}{sl}
\DeclareMathAlphabet{\Ikbb}{U}{bbmss}{m}{it}

\DeclareRobustCommand{\KIB}{{\IKbb B}}

\DeclareRobustCommand{\KIH}{{\IKbb H}}

\DeclareRobustCommand{\Bzero }{{\boldsymbol{\mathit 0}}}

\DeclareMathAlphabet{\Bgothic}{U}{euf}{b}{n}

\DeclareRobustCommand{\Bfrakc}{{\Bgothic c}}
\DeclareRobustCommand{\Bfrakd}{{\Bgothic d}}

\DeclareRobustCommand{\Bfrakf}{{\Bgothic f}}
\DeclareRobustCommand{\Bfrakg}{{\Bgothic g}}

\DeclareRobustCommand{\Bfrakm}{{\Bgothic m}}

\DeclareRobustCommand{\Bfrakq}{{\Bgothic q}}

\DeclareMathAlphabet{\gothic}{U}{euf}{m}{n}

\DeclareRobustCommand{\frakq}{{\gothic q}}

\DeclareMathAlphabet{\IKbb}{U}{bbm}{m}{sl}
\DeclareRobustCommand{\KIg}{{\IKbb g}}

\newcommand{\tr}{\mathop{\mbox{tr}}}                  
\newcommand{\Div}{\mathop{\mbox{Div}}}                

\newcommand{\AND}{\quad\mbox{and}\quad}                    
\newcommand{\WITH}{\quad\mbox{with}\quad}                  

\newcommand{\req }[1]{(\ref{#1})}

\newcommand{\fterm}[1]{\fbox{$ \, \displaystyle#1 \, $}}

\newcommand{\substackrel}[2]{%
   \mathop{ #2}\limits_{{\centering #1}}%
}

\begin{document}

\thispagestyle{empty}

\begin{center} 
\textbf{\large {A Variational Framework for the Thermomechanics of Gradient-Extended Dissipative Solids -- with Applications to Diffusion, Damage and Plasticity}}
\end{center}


\begin{center}
Stephan Teichtmeister$^1$, Marc-Andr\'e Keip$^2$
\end{center}


\begin{center}
$^{1}$ Graz University of Technology, Institute of Biomechanics, Austria \\
$^{2}$ University of Stuttgart, Institute of Applied Mechanics (Chair I), Germany
\end{center}

\bigskip

{\footnotesize \noindent \textbf{Abstract.} 
The paper presents a versatile framework for solids which undergo
nonisothermal processes with irreversibly changing microstructure at
large strains. It outlines rate-type and incremental variational
principles for the full thermomechanical coupling in gradient-extended
dissipative materials. It is shown that these principles yield as
Euler equations essentially the macro- and micro-balances as well as
the energy equation. Starting point is the incorporation of the
entropy and entropy rate as  canonical arguments into constitutive
energy and dissipation functions, which  additionally depend on the gradient-extended
mechanical state and its rate, respectively. By means of (generalized)
Legendre transformations, extended variational principles with thermal
as well as mechanical driving forces can be constructed. On the
thermal side, a rigorous distinction between the quantity conjugate to
the entropy and the quantity conjugate to the entropy rate is
essential here. Formulations with mechanical driving forces are
especially suitable when considering possibly temperature-dependent
threshold mechanisms. With regard to variationally consistent
incrementations, we suggest an update scheme which renders the exact
form of the intrinsic dissipation and is highly suitable when considering adiabatic
processes.  It is shown that this proposed numerical algorithm has the structure of an operator split. To underline the broad applicability of the proposed
framework, we  set up three model problems as applications: 
Cahn-Hilliard diffusion coupled with temperature evolution, where we propose a new variational principle in terms of the species flux vector, as well as thermomechanics of gradient damage and gradient plasticity. In a 
numerical example we study the formation of a cross shear
band. 

}

\bigskip

\thispagestyle{empty}

\section[Section1]{Introduction}
%
In order to model dissipative size effects in solid materials that are, for example, related to the width of shear bands or the grain size in polycrystals, nonstandard continuum theories have to be elaborated that are based on characteristic length scale parameters. The idea is to incorporate at least first-order spatial gradients of micro-structural variables that describe (possibly in a homogenized sense) the irreversibly evolving microstructure. It is of great importance here to account for thermomechanical coupling effects such as thermal softening and temperature increase due to intrinsic dissipation. Practical applications include technological processes like sheet-metal forming or extrusion of metals. Within this work, we present a {\it unified framework} for the fully coupled thermomechanics of gradient-extended dissipative solids that is applicable to a wide range of model problems such as diffusion, (crystal) plasticity and damage. The formulation is embedded into the concept of standard dissipative solids that is characterized by energy and dissipation functions, see \citet{biot65}, \citet{ziegler63} and \citet{halphen+nguyen75}. In particular, the strongly coupled multifield problem will exhibit an incremental variational structure which is an extension of the framework of gradient-extended dissipative solids presented in \citet{miehe11,miehe14} towards nonisothermal processes.

Historically, the consideration of long-range effects via length scales can be traced back to the work of the brothers \citet{cosserat+cosserat1896}. They investigated a microstructure with rigid particles by introducing an additional microrotation governed by three additional degrees of freedom. Since then, a lot of research has been done on so called micropolar, micromorphic and microstrain continua, i.e., by \citet{eringen12}, see also \citet{leismann+mahnken15} for a comparative study. The book of \citet{capriz13} provides a general setting for materials with microstructure where order parameters are considered as components of elements of an abstract manifold. Here, we also refer to the work of \citet{mariano02}. In \citet{maugin90} the standard thermodynamic theory of local internal variables is extended by an additional dependency of the energy function  on the first-order spatial gradient of a (not necessarily scalar) internal variable, see also \citet{maugin+muschik94} and \citet{fremond13}. An important ingredient of such theories is a proper energetic treatment of the microstructural processes yielding additional balance-type equations that are coupled with the standard macroscopic balance equations (mass, linear/angular momentum and energy). For an embedding into the theory of micro-forces we refer to \citet{gurtin03}.

Typical applications of above mentioned nonlocal theories including
additional balance equations are theories of phase transformation,
gradient plasticity and gradient damage. Theories of gradient plasticity are dealt with on a
phenomenological level in \citet{aifantis87},
\citet{fleck+hutchinson01}, \citet{gurtin03},
\citet{forest+sievert03}, \citet{gurtin+anand05} just to name a few. Gradient damage theories that are basically in the same spirit
are considered in \citet{peerlings+etal96}, \citet{dimitrijevic+hackl08}, \citet{pham+etal11}, \citet{fremond13} among many others.

An important role in the modeling of dissipative processes in solids is played by {\it incremental variational principles}. They offer an elegant way to couple different physical fields, possess inherent symmetry and require only Jacobian and Hessian matrices of the underlying potential density for an implementation into typical finite element codes. In addition, if a minimization structure is at hand, structural and material stability analysis is possible, i.e., the formation of complex microstructures can be described. For contributions on the variational theory of local plasticity we refer to \citet{hill50}, \citet{simo+honein90}, \citet{hackl97}, \citet{ortiz+stainier99}, \citet{ortiz+repetto99}, \citet{carstensen+hackl+mielke02} and \citet{mosler+bruhns10}. The works of \citet{miehe02}, \citet{petryk03} and \citet{hackl+fischer07} deal with general inelastic behavior. Extensions towards incorporation of (at least) the first gradient of the internal variable into constitutive functions can be found in \citet{muhlhaus+aifantis91}, \citet{fleck+willis09}, \citet{mainik+mielke09}, \citet{nguyen11}, \citet{miehe14} and \citet{lancioni+etal15} for gradient plasticity and in \citet{lorentz+andrieux99}, \citet{mielke+roubicek06}, \citet{dimitrijevic+hackl08} and \citet{pham+etal11} for gradient damage. In addition, we mention the general formulations of gradient-extended dissipative solids by \citet{svendsen04}, \citet{francfort+mielke06} and \citet{miehe11}. A gradient plasticity theory that incorporates a fractional derivative of the plastic strain is presented in \citet{dahlberg+ortiz19}.  Alternatively, in damage mechanics one can account for nonlocality  by considering the gradient of the total strain measure (instead of the damage variable) in the constitutive functions. Non-smooth (hemi)variational formulations of such an approach are presented in, e.g., \citet{placidi+barchiesi18}, \citet{placidi+misra+barchiesi19} and \citet{timofeev+etal2021}.  Studies comparing these formulations with the gradient damage approach  can be found in \citet{le+etal18} and \citet{placidi+barchiesi+misra18}. Note that all works cited so far in this paragraph deal with dissipative processes under isothermal conditions except \citet{hackl97} in the last section. In the latter work, an original idea of \citet{simo+miehe92} is taken up, namely  to introduce a plastic configurational entropy as an additional internal variable that, as part of the total entropy, leaves the internal energy unchanged. 
From the viewpoint of variational principles however, the main difficulty of fully coupled thermomechanics is to account for the internal dissipation in the energy equation. The construction of a {\it single unified potential} from which all coupled thermomechanical field equations follow is outlined in the seminal work of \citet{yang+stainier+ortiz06}. They introduce  a specific integrating factor that makes the linearized integral expressions, showing up in the fully coupled weak form, symmetric. This factor is based on distinguishing the equilibrium  temperature from a so-called external temperature. Both temperatures coincide only at equilibrium. With this variational principle at hand, the fully coupled thermomechanical boundary value problem  can be solved monolithically by a sequence of symmetric algebraic systems. Alternatively, staggered  numerical algorithms were mainly used before which lead to symmetric finite element stiffness matrices for each subproblem, see \citet{simo+miehe92} and \citet{armero+simo92}. In the recent years, the work of \citet{yang+stainier+ortiz06} has been the basis for further model developments in thermoplasticity, see \citet{stainier+ortiz10}, \citet{canadija+mosler11}, \citet{su+stainier+mercier14}, \citet{bartels+etal15}, \citet{mielke16} and \citet{fohrmeister+bartles+mosler2018}. In the latter, a specific phenomenological model of gradient plasticity together with a micromorphic extension in the sense of \citet{forest09} is considered. 
Besides incorporating the first gradient of a generic internal variable into the energy function, \citet{mielke16} also discusses in detail the gradient structure of thermoplasticity. In \citet{bartels+etal15} special focus is put on the thermodynamic (in)consistency of the Taylor-Quinney factor. This factor, usually introduced in an ad-hoc way, takes into account that only a certain amount of the plastic power is transformed into heat as observed in the experiments of \citet{taylor+quinney34}, see also \citet{rosakis+rosakis+ravichandran+hodowany00} and \citet{oliferuk+raniecki18}.  

According to the authors' knowledge, a {\it universal variational framework} for the thermomechanics of solids with energetic as well as dissipative gradient extensions is still missing in literature. The present work outlines a {\it generalization} of the versatile framework
of \citet{miehe11,miehe14} for gradient-extended dissipative solids
towards {\it nonisothermal processes}.  Point of departure is the
definition of two constitutive functions, namely (i) the internal
energy function that depends on the (gradient-extended) mechanical
state and the entropy and (ii) a (nonsmooth) dissipation potential function that
depends on the rates of the (gradient-extended) mechanical state and the 
entropy. Then, at current time the rates of the macro- and
micro-motion and the entropy rate are governed by a canonical saddle
point principle in case of an adiabatic process and a canonical
minimization principle in case of a process including heat
conduction. However, since specific forms of both constitutive
functions are difficult to access, (generalized) Legendre
transformations of those functions are considered. Here, the concept of 
{\it dual variables} is rigorously pursued, i.e., the quantity
conjugate to the entropy (the temperature) is distinguished from the
quantity conjugate to the entropy rate (the thermal driving
force). Then, the necessary condition of the (generalized) Legendre
transformation of the dissipation potential function must render the
structure of the energy equation in form of an entropy evolution
equation. The obtained rate-type variational principle is in line with
\citet{yang+stainier+ortiz06} but derived in an { \it alternative} way and
extended by an additional long-range micro-motion that accounts for
effects related to length scales. In addition, extended rate-type
variational principles are constructed that incorporate dissipative
mechanical driving forces and a temperature-dependent scalar threshold function via a
Lagrange multiplier. By construction of consistent  numerical time integration
algorithms, variational principles valid for the time increment under
consideration are obtained. Here, we especially point out a {\it
  new semi-explicit numerical algorithm} that, in contrast to a fully implicit scheme, evaluates the integration factor of \citet{yang+stainier+ortiz06} to the value one
and, hence, gives the {\it
algorithmically  consistent form} of the intrinsic dissipation.

In the subsequent sections, the above aspects are carried out in detail and the novel key contributions of the presented work for the thermomechanics of dissipative multifield problems are
\begin{itemize}
\parindent0mm
 \item a canonical variational principle based on a dissipation potential function that depends on the {\it entropy rate}, 

 \item the construction of rate-type and incremental mixed variational principles via {\it (generalized) Legendre transformations}, where we propose a {\it new semi-explicit} numerical update scheme that has the structure of an operator split and is highly suitable when considering adiabatic processes, 

\item a {\it general} and {\it versatile framework} for the full thermomechanical coupling in gradient-extended dissipative continua and

\item applications to gradient plasticity, gradient damage and Cahn-Hilliard diffusion, where for the latter we propose a {\it new minimization structure} with respect to the species flux vector.
\end{itemize}

\vspace*{-2mm}
\noindent
The paper is organized as
follows: In Section~\ref{Section2} we consider  as motivating example a very simple adiabatic
thermomechanical material element and derive rate-type as well as
incremental variational principles  based on the entropy and entropy rate as canonical variables in the energy and dissipation potential functions.  Two update schemes, namely the known implicit and a new semi-explicit numerical algorithm are presented that govern the time-discrete
evolution of the thermomechanical state. In Section~\ref{Section3} we
set up general forms of governing equations for the thermomechanics of
gradient-extended dissipative solids in a three-dimensional
large-strain continuum setting. Section~\ref{Section4} shows that this
fully coupled system is related to a variational statement. In
addition, we construct extended rate-type and incremental variational
principles that incorporate dissipative mechanical driving forces and
threshold mechanisms that depend on the thermal state. The formulations outlined in
Section~\ref{Section4} are general and can be applied to a wide
spectrum of problems  in thermomechanics. As examples, we specify in Section~\ref{Section5}
the developed setting to Cahn-Hilliard diffusion, gradient damage and
(additive) gradient plasticity with strong coupling to temperature
evolution.  The former is given a {\it new variational structure} in terms of the species flux vector. Finally,  to demonstrate the capability of the newly proposed semi-explicit variational update scheme, we show a numerical example
that is concerned with  adiabatic shear band localization in softening
plasticity.

\section[Section2]{Variational Principles for Nonisothermal Rheology}
\label{Section2}
%
\subsection[Section21]{An Adiabatic Thermomechanical Material Element}
\label{Section21}
%
\subsubsection[Section211]{Constitutive Functions}
We consider a very simple rheological material element depicted in Figure~\ref{rheo} describing a smooth linear thermo-visco-elastic behavior. The material element is understood to be embedded in a thermal environment characterized by the absolute temperature $\theta>0$. The spring  exhibits thermoelastic behavior by Hooke's law coupled with thermal expansion. The dashpot device characterizes  viscous response via Newton's law. In addition, the material element is able to store heat.
We have the fundamental constitutive relationships
\begin{equation*}
\begin{array}{llll}
\bullet & \mbox{Hooke's law}  & \sigma^e &= E \, [\, \varepsilon^e - \alpha_T(\theta-\theta_0)\, ]  \, , \\
\bullet & \mbox{Newton's law}  & \sigma^d& = H \dot \varepsilon^d \,  , \\
\bullet & \mbox{Heat storage}  & e^{\theta}& = C(\theta-\theta_0) \, .
\end{array}
\end{equation*}
Here, $\sigma$ denotes as usual a stress, $\varepsilon$ a strain  and $e$ an internal energy. The four involved material parameters are Young's modulus $E$, the thermal expansion coefficient $\alpha_T$, the viscosity $H$ and the heat capacity $C$ at constant internal and external deformation. Without loss of generality, these material parameters are assumed to be constant, i.e., they do not depend on the temperature. Finally, $\theta_0$ stands for a reference temperature. Such a simple model is able to describe classical thermomechanical coupling effects, such as the Gough-Joule effect or thermal expansion, as well as heating due to dissipation in the dashpot. The mechanical state of the rheological device is described by the total strain $\varepsilon$  and the internal variable $\frakq$ characterizing the strain of the inner dashpot. 
%
\inputfig[t]{floats/rheo}{rheo}%
%
We summarize these two quantities in the 
\begin{equation*}
\mbox{mechanical state array: } \quad \Bfrakc = ( \, \varepsilon, \frakq   \, ) \, .
\end{equation*}
Identifying the strain in the spring device by the simple kinematic relationship $\varepsilon^e=\varepsilon-\frakq$, we define the well known thermoelastic free energy function
\begin{equation}
\widehat \psi(\Bfrakc,\theta)=\underbrace{\vphantom{\ln\frac{\theta}{\theta_0}}\frac{1}{2}E(\varepsilon-\frakq)^2}_{=:\widehat \psi_e(\Bfrakc)} + \underbrace{\vphantom{\ln\frac{\theta}{\theta_0}} E\alpha_T(\frakq - \varepsilon)(\theta-\theta_0)}_{=:\widehat \psi_{e-\theta}(\Bfrakc,\theta)} + \underbrace{ C[(\theta-\theta_0) - \theta \ln\frac{\theta}{\theta_0} ]}_{=:\widehat \psi_{\theta}(\theta)} \,  ,
\label{0d_cons_1}
\end{equation}
that is not defined for nonpositive temperatures  and is concave in $\theta$.
In addition, we define the dissipation potential function associated with the two dashpot devices
\begin{equation}
\widehat \phi(\dot \Bfrakc)= \frac{1}{2}H_1 \dot \varepsilon^2+\frac{1}{2} H_2 \dot \frakq^2 \, .
\label{0d_cons_2} 
\end{equation}
\subsubsection[Section212]{Governing System of Equations}
\label{Section212}
The thermomechanical response of the material element within a time interval $\calT = (0,t^e)\subset \mathbb{R}_+$ is fully governed by the two constitutive functions \req{0d_cons_1} and \req{0d_cons_2} and characterized by four equations \req{0d_gov_1}--\req{0d_gov_4}. First, on the mechanical side, the equilibrium equation 
\begin{equation}
\partial_{\varepsilon} \widehat \psi + \partial_{\dot \varepsilon} \widehat \phi = \sigma_{ext}
\quad
\mbox{in }
\calT
\label{0d_gov_1}
\end{equation}
is a relationship between internal and external stresses. Second, the evolving internal strain state is governed by Biot's equation
\begin{equation}
\partial_{\frakq} \widehat \psi + \partial_{\dot \frakq} \widehat \phi = 0
\quad
\mbox{in }
\calT 
\label{0d_gov_2}
\end{equation}
and characterizes an internal dissipative mechanism,  see \citet{biot65}. Third, on the thermal side, the state equation
\begin{equation}
\eta = - \partial_\theta \widehat \psi
\quad
\mbox{in }
\calT
\label{0d_gov_3}
\end{equation}
determines the current entropy $\eta$ or in an inverse manner the current temperature $\theta$. Finally, the evolution of the entropy is governed by a nonnegative dissipation
\begin{equation}
\dot \eta = \frac{1}{\theta} \calD 
\quad
\mbox{in }
\calT
\WITH \calD = \partial_{\dot \Bfrakc} \widehat \phi \cdot \dot \Bfrakc \geq 0 \, .
\label{0d_gov_4}
\end{equation}
Equation \req{0d_gov_4}$_1$ together with the {\it constitutive definition} of the dissipation $\calD$ is a form of the first law of thermodynamics,
i.e., the balance of energy for an adiabatic process under
consideration. The inequality \req{0d_gov_4}$_2$ is the second law of thermodynamics
that is ensured {\it a priori} by a dissipation potential function $\widehat \phi(\cdot)$ that is (i) nonnegative, (ii) zero in the origin $\widehat \phi(\Bzero) = 0$ and (iii) convex,  see, e.g., \citet{fremond13}. In case of a dissipation potential function that is positively homogeneous of degree $\ell$, i.e., $\widehat \phi(\gamma \dot \Bfrakc) = \gamma  ^\ell \widehat \phi(\dot \Bfrakc)$ for all $\gamma > 0$, the dissipation  can be written alternatively as\footnote{In case of viscous dissipation we have $\ell = 2$, see \req{0d_cons_2}.} $\calD = \ell \widehat \phi(\dot \Bfrakc)$. 
Aim of the subsequent treatment is the construction of variational
principles that account for the governing equations
\req{0d_gov_1}--\req{0d_gov_4} of the thermomechanical device.  Starting point is the definition of a (maybe unusual) dissipation potential function that depends on the {\it entropy rate}.

\subsection[Section22]{ Canonical Variational Principle with Entropy Variable}
\label{Section22}
%
The main difficulty in the formulation of rate-type variational principles in thermo-mechanics is to account for the full update of the thermal state via \req{0d_gov_4}.  For the nonisothermal case, we start by modeling the internal energy $e$  whose corresponding function naturally depends on the entropy $\eta$.  Hence, our investigation of thermomechanics in generalized standard materials is based on assuming the constitutive functions
\begin{equation}
e(t) = \widehat e (\Bfrakc, \eta) 
\AND
v(t) = \widehat v (\dot \Bfrakc, \dot \eta; \Bfrakc, \eta) \, .
\label{0d_canonical_1}
\end{equation}
The internal energy function $\widehat e$ determines the current internal energy stored in the material element. The  newly introduced dissipation potential function $\widehat v$ is assumed to be a function of the rates $(\dot \Bfrakc,\dot \eta)$ and might additionally depend on the current thermomechanical state $(\Bfrakc, \eta )$. Alternatively, instead of the current entropy $\eta$ the dissipation potential function $\widehat v$ can  depend on the current temperature $\theta$ by making use of the constitutive relationship \req{0d_gov_3}. Comparing $\widehat v$ with the dissipation potential function $\widehat \phi$ defined in \req{0d_cons_2}, we observe an additional thermal slot with the dependence on $\dot \eta$. This assumed dependency is the { \it key ingredient} of the subsequent construction of variational principles in dissipative thermomechanics. Based on the two constitutive functions \req{0d_canonical_1}, we construct the rate-type potential
\begin{equation}
\fterm{
\widehat p (\dot \Bfrakc, \dot \eta; \Bfrakc, \eta)= \frac{d}{dt}\widehat e(\Bfrakc, \eta) + \widehat v(\dot \Bfrakc, \dot \eta; \Bfrakc, \eta)
}
\label{0d_canonical_2} 
\end{equation}
at given thermomechanical state $(\Bfrakc, \eta)$, that accounts for energetic as well as dissipative mechanisms. Next, we define an external load function 
\begin{equation*}
\widehat p_{ext}(\dot \varepsilon;t) = \sigma_{ext}(t) \dot \varepsilon 
\end{equation*}
depending on the given external stress $\sigma_{ext}$, that loads the material element, see Figure~\ref{rheo}. Then, the rate of the thermomechanical state at current time $t$ under adiabatic conditions is determined by the {\it canonical saddle point principle}  
\begin{equation}
\{ \dot \Bfrakc, \dot \eta \} =
\mbox{Arg} \{\;
\substackrel{\dot \Bfrakc }{\mbox{inf}}\,
\substackrel{\dot \eta}{\mbox{sup}} \;
[ \, 
\widehat p (\dot \Bfrakc, \dot \eta; \Bfrakc, \eta)-\widehat p_{ext}(\dot \varepsilon;t) \, ] 
\; \}
\, .
\label{0d_canonical_4}
\end{equation}
Note that the saddle point structure of this variational principle is governed by an {\it assumed concavity} of the dissipation potential function $\widehat v$ with respect to $\dot \eta$, see Section~\ref{Section232} and Footnote~\ref{footnote_dissipation}. Taking the first derivative of the potential \req{0d_canonical_2} yields as necessary conditions of the variational principle \req{0d_canonical_4} three Euler equations for the rate of the thermomechanical state of the material element at current time $t$, namely
\begin{equation*}
\parbox{14cm}{\centering 
$
\begin{array}{llll}
1.&{\mbox{ Evolving external state}} &
 \partial_{\dot \varepsilon} \widehat p \equiv \partial_{\varepsilon}\widehat e + \partial_{\dot \varepsilon} \widehat v &= \sigma_{ext} \,  ,\\[1ex]
2. &{\mbox{ Evolving internal state}} & 
\partial_{\dot \frakq} \widehat p \equiv \partial_{\frakq} \widehat e + \partial_{\dot \frakq} \widehat v &= 0   \, , \\[1ex]
3. &{\mbox{ Evolving thermal state}} &
\partial_{\dot \eta} \widehat p \equiv  \partial_{\eta} \widehat e + \partial_{\dot \eta} \widehat v &= 0 \, .
\end{array}
$
}
\end{equation*}
The first equation determines the rate $\dot \varepsilon$ of the external strain state and is a form of the stress equilibrium \req{0d_gov_1} without inertia terms. The second equation governs the rate $\dot \frakq$ of the internal variable and is identical to \req{0d_gov_2}. Finally, the third equation is a form of the balance of energy, i.e., the first law of thermodynamics \req{0d_gov_4} combined with an inverse representation of \req{0d_gov_3}, which determines the rate $\dot \eta$ of the entropy. However, this is a quite unusual representation. The main difficulty in this canonical setting is the formulation of the two constitutive functions $\widehat e$ and $\widehat v$ in terms of the entropy $\eta$ and entropy rate $\dot \eta$, respectively\footnote{For a single linear viscous dashpot with current temperature $\theta$ and entropy $\eta$, respectively, and rates $(\dot \varepsilon, \dot \eta)$, the internal energy function $\widehat e$ and the dissipation potential function $\widehat v$ take,  by evaluating the partial Legendre transformations \req{0d_dual_1} and \req{0d_dual_3}, the nonintuitive forms
\begin{equation*}
\widehat e (\eta) = C \theta_0 (e^{\eta/C} - 1) 
\AND
\widehat v (\dot \varepsilon, \dot \eta; \theta) = - \frac{\theta^2}{2H} \left(\frac{\dot \eta}{\dot \varepsilon}\right)^2 \, .
\end{equation*}
\label{footnote_dissipation}
}. Hence, a formulation in terms of the temperature is  more convenient.

\subsection[Section23]{Mixed Variational Principle with Temperature Variable}
\label{Section23}
%
 For practical applications, the above canonical variational principle is transformed into a mixed saddle point principle that incorporates the temperature as a primary variable. This transformation is achieved within two steps.

\subsubsection[Section231]{Variable dual to Entropy} \label{Section231}
First, we define the internal energy function occuring in the rate-type potential \req{0d_canonical_2} by a partial Legendre transformation
\begin{equation}
\widehat e (\Bfrakc, \eta) = \sup_{\theta} [\, \widehat \psi(\Bfrakc, \theta) + \theta \eta   \,] \, ,
\label{0d_dual_1}
\end{equation}
where $\widehat \psi$ is the free energy function depending on the thermal variable $\theta$.
The necessary optimality condition of this Legendre transformation defines the entropy
\begin{equation}
\eta = -\partial_{\theta} \widehat \psi(\Bfrakc, \theta)
\label{0d_dual_2}
\end{equation}
in terms of the thermal variable $\theta$.  The latter is the dual quantity to the entropy $\eta$, i.e., the {\it temperature} of the material element.
\subsubsection[Section232]{Variable dual to Entropy Rate} \label{Section232}
In a second step, we define the dissipation potential function $\widehat v$ occuring in the rate-type potential \req{0d_canonical_2} by a partial Legendre transformation
\begin{equation}
\widehat v (\dot \Bfrakc, \dot \eta; \Bfrakc, \eta) = \inf_{T} [\, \widetilde \phi(\dot \Bfrakc, T; \Bfrakc, \theta) - T \dot \eta   \,]
\label{0d_dual_3}
\end{equation}
in terms of the dissipation potential function $\widetilde \phi$, that depends on the thermal variable $T$ and, by making use of \req{0d_dual_2}, is defined at current state $( \Bfrakc, \theta )$. The necessary condition of this Legendre transformation defines the evolution of the entropy
\begin{equation}
\dot \eta = \partial_T \widetilde \phi(\dot \Bfrakc, T; \Bfrakc, \theta)
\label{0d_dual_4}
\end{equation}
in terms of the thermal variable $T$. We identify $T$ as the quantity dual to the entropy rate $\dot \eta$ and call it the {\it thermal driving force} that, up to this point, is rigorously distinguished from the temperature $\theta$ of the material element. For the adiabatic case under consideration, \req{0d_dual_4} {\it must} have the form of the balance of energy \req{0d_gov_4} and we identify
\begin{equation}
\partial_T \widetilde \phi  \overset{!}{=} \frac{1}{\theta} \partial_{\dot \Bfrakc} \widehat \phi \cdot \dot \Bfrakc \, .
\label{0d_dual_5}
\end{equation}
This relationship restricts the form of the functional dependence of the dissipation potential function $\widetilde \phi$ and is fulfilled if a multiplicative dependence on the thermal driving force $T$ is assumed
\begin{equation}
\widetilde{\phi}(\dot \Bfrakc, T; \Bfrakc, \theta) = \widehat \phi(\frac{T}{\theta} \dot \Bfrakc; \Bfrakc, \theta) \, .
\label{0d_dual_6}
\end{equation}
To arrive at \req{0d_dual_5}, we use the result that at equilibrium the thermal driving force $T$ can be identified with the temperature $\theta$, which is an outcome of the mixed variational principle \req{0d_dual_7} below. The scaling factor $T/\theta$ first arised  via a symmetry argument as integrating factor in the seminal work of \citet{yang+stainier+ortiz06}.

\subsubsection[Section233]{Mixed Variational Principle} \label{Section233} 
Starting from the rate-type potential \req{0d_canonical_2}, the two partial Legendre transformations \req{0d_dual_1} and \req{0d_dual_3} motivate the definition of a mixed rate-type potential\footnote{To keep notation short, we subsequently do not write explicitly the dependence of functions on given states.}  
\begin{equation}
\widehat \pi(\dot \Bfrakc, \dot \eta, \dot \theta, T) := \frac{d}{dt} \widehat \psi(\Bfrakc, \theta) + (\theta-T) \dot \eta + \eta \dot \theta + \widehat \phi(\frac{T}{\theta} \dot \Bfrakc)
\label{0d_dual_6a}
\end{equation}
at given thermomechanical state $( \Bfrakc, \eta, \theta )$ in terms of the two constitutive functions $\widehat \psi$ and $\widehat \phi$. Inserting the necessary condition \req{0d_dual_2} on the given thermomechanical state yields the reduced mixed rate-type potential
\begin{equation*}
\widehat \pi_{red}(\dot \Bfrakc, \dot \eta, T)=\partial_{\Bfrakc} \widehat \psi \cdot \dot \Bfrakc + (\theta-T) \dot \eta + \widehat \phi(\frac{T}{\theta} \dot \Bfrakc)
\end{equation*}
at given thermomechanical state $( \Bfrakc, \eta, \theta )$. Based on this definition, we obtain the mixed saddle point principle
\begin{equation}
\{ \dot \Bfrakc, \dot \eta, T \} =
\mbox{Arg} \{\;
\substackrel{\dot \Bfrakc }{\mbox{inf}}\,
\substackrel{\dot \eta}{\mbox{sup}} \,
\substackrel{T}{\mbox{inf}} \;
[ \, 
\widehat \pi_{red} (\dot \Bfrakc, \dot \eta,T)- \widehat p_{ext}(\dot \varepsilon;t) \, ] 
\; \} \,  ,
\label{0d_dual_7}
\end{equation}
that determines at current time $t$ the rates of the external strain, internal variable and entropy as well as the thermal driving force. The Euler equations of this saddle point principle are
\begin{equation}
\parbox{14cm}{\centering 
$
\begin{array}{llll}
1.&{\mbox{ Evolving external state}} &
 \partial_{\dot \varepsilon} \widehat \pi_{red} \equiv \partial_{\varepsilon}\widehat \psi + \partial_{\dot \varepsilon} \widehat \phi &= \sigma_{ext} \,  ,  \\[1mm]
2. &{\mbox{ Evolving internal state}} & 
\partial_{\dot \frakq} \widehat \pi_{red} \equiv \partial_{\frakq} \widehat \psi + \partial_{\dot \frakq} \widehat \phi &= 0 \,  , \\[1mm]
3. &{\mbox{ Thermal driving force}} &
\partial_{\dot \eta} \widehat \pi_{red} \equiv  T- \theta &= 0 \,  , \\[1mm]
4.&{\mbox{ Evolving thermal state}} &
 \partial_{T} \widehat \pi_{red} \equiv - \dot \eta + \partial_{T} \widehat \phi  &= 0 \, .
\end{array}
$
}
\label{0d_dual_8}
\end{equation}
Again, the first equation represents the stress equilibrium and the second one Biot's equation. The third relationship identifies the thermal driving force with the given temperature as reported in \citet{yang+stainier+ortiz06}. The last equation represents the first law of thermodynamics in the form of an evolution equation for the entropy. Clearly, within this mixed setting the entropy rate $\dot \eta$ plays the role of a {\it Lagrange multiplier}. With known rates $(\dot \Bfrakc, \dot \eta)$ of mechanical state and entropy, the temperature rate at current time $t$ can be computed by taking the time derivative of \req{0d_gov_3}.
A constrained minimization formulation dual to the variational principle \req{0d_canonical_4} is obtained by commuting the order of the infimum and the supremum in \req{0d_dual_7} (which we assume to be possible) such that
\begin{equation*}
\substackrel{\dot \Bfrakc }{\mbox{inf}}\,
\substackrel{\dot \eta}{\mbox{sup}} \,
\substackrel{T}{\mbox{inf}} \;
[ \, 
\widehat \pi_{red} (\dot \Bfrakc, \dot \eta,T)-\widehat p_{ext}(\dot \varepsilon;t)  \, ] 
=
\substackrel{\dot \Bfrakc }{\mbox{inf}}\,
\substackrel{T \in \req{0d_dual_8}_3}{\mbox{inf}} \;
[ \, 
\partial_{\Bfrakc}\widehat \psi \cdot \dot \Bfrakc + \widehat \phi(\frac{T}{\theta} \dot \Bfrakc)-\widehat p_{ext}(\dot \varepsilon;t) \, ] \, .
\end{equation*}
 For a general discussion of constructing mixed and dual variational principles we refer to the books \citet{maugin92}, \citet{nguyen00} and \citet{berdichevsky09}.

\subsection[Section24]{Incremental Variational Principles}
\label{Section24}
%
We consider a finite time interval $[t_n, t_{n+1}] \subset \overline \calT$ with step length $\tau:=t_{n+1}-t_n>0$ and assume all thermomechanical state quantities at time $t_n$ to be known. The goal is then to determine all  the approximate state quantities at time $t_{n+1}$ based on variational principles valid for the time increment under consideration. Subsequently, all variables without subscript are understood to be evaluated at time $t_{n+1}$.

\subsubsection[Section241]{Mixed Variational Principle} \label{Section241}
A mixed variational principle associated with the finite time interval $[t_n, t_{n+1}]$ is based on the incremental potential
\begin{equation}
\widehat \pi^{\tau}(\Bfrakc, \eta, \theta, T) = Algo \{ \ \int_{t_n}^{t_{n+1}} \widehat \pi(\dot \Bfrakc, \dot \eta, \dot \theta, T) \, dt   \ \}
\label{0d_inc_1}
\end{equation}
in terms of the continuous rate-type potential $\widehat \pi$ defined in \req{0d_dual_6a}. Here, $Algo$ stands for a  numerical integration algorithm in time within the interval $[t_n,t_{n+1}]$. The incremental potential $\widehat \pi^\tau$ is understood to be defined at given thermomechanical state $(\Bfrakc_n, \eta_n, \theta_n )$ at time $t_n$. Then, the finite-step sized incremental mixed variational principle 
\begin{equation}
\{ \Bfrakc, \eta, \theta, T \} =
\mbox{Arg} \{\;
\substackrel{\Bfrakc }{\mbox{inf}}\,
\substackrel{\eta, \theta}{\mbox{sup}} \,
\substackrel{T}{\mbox{inf}} \;
[ \, 
\widehat \pi^\tau (\Bfrakc, \eta, \theta, T)-\sigma_{ext}(t_{n+1}) \varepsilon \, ] 
\; \}
\label{0d_inc_2}
\end{equation}
determines the mechanical state, the entropy, the temperature and the thermal driving force at current time $t_{n+1}$. The $Algo$ operator is constructed such that the incremental variational principle \req{0d_inc_2} yields as Euler equations {\it consistent algorithmic counterparts} of the Euler equations \req{0d_dual_8} stemming from the rate-type variational principle \req{0d_dual_7}. In what follows two different forms of the $Algo$ operator are discussed.

\subsubsection[Section242]{Implicit Variational Update} \label{Section242}
As a first approach, we  consider a  numerical algorithm that performs an {\it exact incrementation} of the internal energy
\begin{equation}
\int_{t_n}^{t_{n+1}} \frac{d}{dt}[\, \widehat \psi(\Bfrakc, \theta) + \theta \eta \, ] \, dt =  \widehat \psi(\Bfrakc, \theta) + \theta \eta - \widehat \psi(\Bfrakc_n, \theta_n) - \theta_n \eta_n 
\label{0d_inc_3}
\end{equation}
and an incrementation of the dissipative term according to a backward Euler scheme
\begin{equation}
\int_{t_n}^{t_{n+1}}[ \, -T \dot \eta  +\widehat \phi(\frac{T}{\theta}\dot \Bfrakc)   \,] \, dt \approx -T (\eta - \eta_n) + \tau \widehat \phi(\frac{T}{\theta_n} \dot \Bfrakc^\tau) \, ,
\label{0d_inc_4}
\end{equation}
see \citet{yang+stainier+ortiz06}. Here, $\dot \Bfrakc^\tau = (\Bfrakc - \Bfrakc_n)/\tau$ denotes an algorithmic expression of the rate of the mechanical state. The incremental potential \req{0d_inc_1} takes the form
\begin{equation*}
\fterm{
\widehat \pi^\tau(\Bfrakc, \eta, \theta, T) = \widehat \psi(\Bfrakc, \theta) - \widehat \psi(\Bfrakc_n, \theta_n) + \theta \eta - \theta_n \eta_n  -T (\eta - \eta_n) + \tau \widehat \phi(\frac{T}{\theta_n} \dot \Bfrakc^\tau) \, .
}
\end{equation*}
Then, the incremental saddle point principle \req{0d_inc_2} gives the Euler equations
\begin{equation}
\parbox{14cm}{\centering 
$
\begin{array}{llll}
1.&{\mbox{ Update external state}} &
  \partial_{\varepsilon} \widehat \pi^\tau \equiv \partial_{\varepsilon}\widehat \psi + \tau \partial_{\varepsilon} \widehat \phi &= \sigma_{ext} \,  ,  \\[1mm]
2. &{\mbox{ Update internal state}} & 
 \partial_{ \frakq} \widehat \pi^\tau \equiv \partial_{\frakq}\widehat \psi + \tau \partial_{ \frakq} \widehat \phi &= 0 \,  , \\[1mm]
3. &{\mbox{ Thermal driving force}} &
   \partial_{ \eta} \widehat \pi^\tau \equiv T- \theta &= 0 \,  , \\[1mm]
4.&{\mbox{ Current thermal state}} &
 \partial_{ \theta} \widehat \pi^\tau \equiv \partial_{\theta}\widehat \psi + \eta   &= 0 \,  , \\[1mm]
5.&{\mbox{ Update entropy}} &
 \partial_{T} \widehat \pi^\tau \equiv -(\eta-\eta_n) + \tau \partial_T \widehat \phi   &= 0 \, .
\end{array}
$
}
\label{0d_inc_6}
\end{equation}
The fifth equation represents a time-discrete form of the energy equation and governs the update of the entropy $\eta$. On the thermal side, the key equations are \req{0d_inc_6}$_4$ and \req{0d_inc_6}$_5$ 
\begin{equation}
\eta = - \partial_{\theta} \widehat \psi (\Bfrakc, \theta) \, ,
\quad
\eta = \eta_n + \frac{\tau}{\theta} \, \ell  \widehat \phi (\frac{\theta}{\theta_n} \dot \Bfrakc^\tau) \, ,
\label{0d_inc_7}
\end{equation}
where $T$ has been eliminated by \req{0d_inc_6}$_3$ and a positively homogeneous dissipation potential function $\widehat \phi$ of degree $\ell$ is assumed. By inserting \req{0d_inc_7}$_1$ into \req{0d_inc_7}$_2$ an update equation for the temperature $\theta$ is obtained that includes contributions from thermoelastic heating and heating due to dissipation. The last term in \req{0d_inc_7}$_2$ contains the algorithmic counterpart of the dissipation \req{0d_gov_4}$_2$. However, the argument of the dissipation potential function $\widehat \phi$ in \req{0d_inc_7}$_2$ is {\it scaled} by a factor $\theta/\theta_n$ which is a consequence of the implicit variational update. The same scaling factor arises in the dissipative terms of the Euler equations \req{0d_inc_6}$_1$ and \req{0d_inc_6}$_2$. Again, we refer to \citet{yang+stainier+ortiz06}. Finally,  we observe that the incremental updates of the mechanical and thermal variables are {\it fully coupled}.

\subsubsection[Section243]{Semi-Explicit Variational Update} \label{Section243}
For the construction of a variational update alternative to the  above implicit approach, we  propose an {\it approximate incrementation} of the internal energy
\begin{equation}
\int_{t_n}^{t_{n+1}} \frac{d}{dt} [ \, \widehat \psi(\Bfrakc, \theta) + \theta \eta \, ] \, dt \approx \widehat \psi(\Bfrakc, \theta) - \widehat \psi(\Bfrakc_n, \theta_n) + (\theta - \theta_n) \eta_n + (\eta - \eta_n) \theta_n
\label{0d_inc_9}
\end{equation}
based on an explicit integration of the second thermal term in the integrand. Note that $\eta$ and $\theta$ are given quantities in the above rate-type formulation and as such, they are frozen in the (approximate) evaluation of the integral. By \req{0d_inc_9} and the incrementation \req{0d_inc_4} of the dissipative term, the incremental potential \req{0d_inc_1} takes the form 
\begin{equation}
\fterm{
\widehat \pi^\tau(\Bfrakc, \eta, \theta, T) = \widehat \psi(\Bfrakc, \theta) -  \widehat \psi(\Bfrakc_n, \theta_n)  + (\theta - \theta_n) \eta_n + (\theta_n - T)(\eta - \eta_n) + \tau \widehat \phi(\frac{T}{\theta_n} \dot \Bfrakc^\tau) \, .
}
\label{0d_inc_10}
\end{equation}
As pointed out in Section~\ref{Section244}, this potential contains an intrinsic operator split. The incremental saddle point principle \req{0d_inc_2} gives the Euler equations
\begin{equation}
\parbox{14cm}{\centering 
$
\begin{array}{llll}
1.&{\mbox{ Update external state}} &
  \partial_{\varepsilon} \widehat \pi^\tau \equiv \partial_{\varepsilon}\widehat \psi + \tau \partial_{\varepsilon} \widehat \phi &= \sigma_{ext} \,  ,  \\[1mm]
2. &{\mbox{ Update internal state}} & 
 \partial_{ \frakq} \widehat \pi^\tau \equiv \partial_{\frakq}\widehat \psi + \tau \partial_{\frakq} \widehat \phi &= 0 \,  , \\[1mm]
3. &{\mbox{ Thermal driving force}} &
   \partial_{ \eta} \widehat \pi^\tau \equiv T- \theta_n &= 0 \,  , \\[1mm]
4.&{\mbox{ Current thermal state}} &
 \partial_{ \theta} \widehat \pi^\tau \equiv \partial_{\theta}\widehat \psi + \eta_n   &= 0  \, , \\[1mm]
5.&{\mbox{ Update entropy}} &
 \partial_{T} \widehat \pi^\tau \equiv -(\eta-\eta_n) + \tau \partial_T \widehat \phi   &= 0 \, .
\end{array}
$
}
\label{0d_inc_11}
\end{equation}
Comparing this set of equations with \req{0d_inc_6}, we observe {\it essential modifications} in \req{0d_inc_11}$_3$ and \req{0d_inc_11}$_4$. These equations identify the current thermal driving force $T$ with the {\it given} temperature $\theta_n$ at time $t_n$ and define the current temperature $\theta$ in terms of the current mechanical state $\Bfrakc$ and the {\it given} entropy $\eta_n$ at time $t_n$. The key equations \req{0d_inc_11}$_4$ and \req{0d_inc_11}$_5$ on the thermal side can be recast into
\begin{equation}
 \partial_{\theta} \widehat \psi(\Bfrakc, \theta) = - \eta_n \, ,
\quad
\eta = \eta_n + \frac{\tau}{\theta_n} \,  \ell \widehat \phi(\dot \Bfrakc ^\tau) \, ,
\label{0d_inc_12}
\end{equation}
where $T$ has been eliminated by \req{0d_inc_11}$_3$ and again a positively homogeneous dissipation potential function $\widehat \phi$ of degree $\ell$ is assumed. One observes, that in sharp contrast to the governing equation \req{0d_inc_7}$_2$ stemming from the implicit variational approach, equation  \req{0d_inc_12}$_2$ does {\it not} contain the scaling factor. Hence, by the explicit variational update the { \it algorithmically   consistent} form $\ell \widehat \phi(\dot \Bfrakc ^\tau)$ of the  dissipation \req{0d_gov_4}$_2$ is obtained. Furthermore, also the dissipative terms in the stress equilibrium \req{0d_inc_11}$_1$ and Biot's equation \req{0d_inc_11}$_2$ do not have the scaling factor. Finally,  we observe that as a consequence of the temperature update \req{0d_inc_12}$_1$ for given entropy $\eta_n$, a {\it decoupled} incremental formulation is at hand.

\subsubsection[Section244]{Operator Split} \label{Section244}
The two update equations \req{0d_inc_12} characterize a
variational-based staggered algorithmic treatment of dissipative
thermomechanics based on the potential $\widehat \pi^\tau$ defined in \req{0d_inc_10}:  (i) update the mechanical state $\Bfrakc$ as well as
the temperature $\theta$ at frozen entropy $\eta_n$ and given (predicted) thermal
driving force $T=\theta_n$, and (ii) update the entropy $\eta$ as well
as the thermal driving force $T$ which turns out to be identical to the predicted one.  Hence, the  semi-explicit variational update in the time interval
$[t_n, t_{n+1}]$ can be seen as a composition of two {\it fractional
steps}
\begin{equation*}
Algo = Algo_{\eta,T} \circ Algo_{\Bfrakc, \theta} \,  ,
\end{equation*}
which are both variational. This  numerical algorithm falls under the category of {\it incrementally isentropic operator splits}, 
a specific form of which is considered in \citet{armero+simo92}.
Within the {\it isentropic predictor step} we first
solve the variational sub-problem
\begin{equation}
(Algo_{\Bfrakc, \theta}):
\quad 
\{ \Bfrakc^*, \theta^* \} =
\mbox{Arg} \{\;
\substackrel{\Bfrakc }{\mbox{inf}}\,
\substackrel{\theta}{\mbox{sup}} \;
[ \, 
\widehat \pi^\tau (\Bfrakc, \eta_n, \theta, \theta_n)-\sigma_{ext}(t_{n+1})\varepsilon \, ] 
\; \} \,  ,
\label{0d_algo_2}
\end{equation}
that gives the mechanical state and the temperature at time $t_{n+1}$. The resulting optimality conditions are the equations \req{0d_inc_11}$_{1-2}$  and \req{0d_inc_11}$_4$. Within a second step, the {\it entropy corrector}, the entropy as well as the thermal driving force at time $t_{n+1}$ are obtained via the variational sub-problem 
\begin{equation*}
(Algo_{\eta, T}):
\quad 
\{ \eta^*, T^* \} =
\mbox{Arg} \{\;
\substackrel{\eta }{\mbox{sup}}\,
\substackrel{T}{\mbox{inf}} \;
 \, 
\widehat \pi^\tau (\Bfrakc^*, \eta, \theta^*, T) 
\; \}
\end{equation*}
for given mechanical state $\Bfrakc^*$ and temperature $\theta^*$ stemming from the isentropic predictor \req{0d_algo_2}. The resulting optimality conditions are the equations \req{0d_inc_11}$_3$ and \req{0d_inc_11}$_5$.

\section{Thermomechanics of Gradient-Extended Dissipative Solids} \label{Section3}
%
 In this section, we shortly summarize the basic constitutive framework for the thermomechanics of gradient-extended dissipative solids in an abstract setting and adopt the notion of \citet{miehe11, miehe14}.
Let $\calB \subset \mathbb{R}^d$ with $d \in \{2,3 \}$ be the reference placement ($d$-manifold) of a material body $B$ into the Euclidean space with smooth boundary $\partial \calB$. We study the thermomechanical behavior of the body under mechanical as well as thermal loadings in a time interval $\calT = (0, t_e) \subset \mathbb{R}_+$. To this end, we focus on a multiscale viewpoint in the sense that we relate dissipative effects at $(\BX,t) \in \calB \times \calT$ to changes in the microstructure (besides the dissipative effect stemming from heat conduction). In the phenomenological context, we account for these microstructural mechanisms by {\it micro-motion fields} that generalize the classical concept of internal variables. The evolution of the related dissipative system is considered to be nonsmooth and therefore the governing equations are written as differential inclusions.
Without loss of generality, we assume within our treatment homogeneous solids only, i.e., material properties do not depend on the position $\BX \in \calB$. In what follows, we denote by $\dot{(\cdot)} = \frac{\partial}{\partial t} (\cdot)(\BX,t)$ the material time derivative and by $\nabla(\cdot) (\BX,t)$ the gradient on the reference manifold $\calB$. In addition, to keep notation short, we understand the operation $\Ba \cdot \Bb$ either as {\it duality product} between a vector and a one-form $\Ba \cdot \Bb = a^c \, b_c$ or as {\it inner product} between two vectors $\Ba \cdot \Bb = a^c \, b^d \, \delta_{cd}$ and two one-forms $\Ba \cdot \Bb = a_c \, b_d \, \delta^{cd}$, respectively. The Kronecker deltas $\delta_{ab}$ and $\delta^{ab}$ represent the coefficients of the metric and inverse metric tensor in a Cartesian coordinate system, respectively.  Finally, we give the usual definition of a double contraction between two second-order tensors, e.g., for two mixed-variant tensors $\BA : \BB = A^{c}{}_d \,  B_{c}{}^d$.

\subsection[Section31]{Basic Fields of a Solid with Microstructural Changes} \label{Section31}

\subsubsection[Section311]{Macro-Motion of a Body} \label{Section311}
Within the geometrically nonlinear theory, the macroscopic motion of the body is described by the {\it macro-motion field}
\begin{equation}
\Bvarphi : 
\left\{
\begin{array}{ll}
  \calB \times \calT \rightarrow \mathbb{R}^d \\
 (\BX, t)  \mapsto \Bx = \Bvarphi(\BX,t) \, ,
\end{array}
\right.
\label{kine_1}
\end{equation}
which maps at time $t \in \calT$ points $\BX \in \calB$ of the reference configuration $\calB$ to points $\Bx \in \Bvarphi_t(\calB)$ of the current configuration $\Bvarphi_t(\calB)$. Let $\BG = \delta_{AB} \BE^A \otimes \BE^B$ and $\Bg = \delta_{ab} \Be^a \otimes \Be^b$ denote the Euclidean metric tensors associated with the reference and current configuration, where the Kronecker symbols refer to Cartesian coordinate systems on both manifolds. The deformation gradient $\BF = D\Bvarphi(\BX,t)$ is defined as the tangent corresponding to the macro-motion \req{kine_1}. Next, the convected current metric $\BC = \BF^T \Bg \BF$ is introduced, that is often denoted as the right Cauchy-Green tensor. The macro-motion \req{kine_1} is constrained by a positive Jacobian $J=\sqrt{\det \BC}>0$ that rules out penetration of matter. With respect to mechanical loading of the solid, the boundary of the reference configuration is decomposed into nonoverlapping parts $\partial \calB_{\Bvarphi}$ and $\partial \calB_{\Bt}$ such that $\partial\calB=\partial \calB_{\Bvarphi} \cup \partial \calB_{\Bt}$. We prescribe the macro-motion (Dirichlet boundary condition)
\begin{equation}
\Bvarphi(\BX,t) = \bar \Bvarphi(\BX,t) \quad \mbox{on } \partial \calB_{\Bvarphi} \times \calT
\label{defo_dirichlet}
\end{equation}
and the {\it macro-tractions} (Neumann boundary condition) on $\partial \calB_{\Bt} \times \calT$ as specified in Section~\ref{Section34} below.
\inputfig[t]{floats/fields}{fields}%
%

\subsubsection[Section312]{Micro-Motion of a Body} \label{Section312}
Microstructural dissipative changes of the body are described by additional fields related to the concept of internal variables. These variables are assembled in the {\it micro-motion field} of the solid
\begin{equation*}
\Bfrakq : 
\left\{
\begin{array}{ll}
  \calB \times  \calT \rightarrow \mathbb{R}^\delta \\
 (\BX, t)  \mapsto \Bfrakq(\BX,t) \, .
\end{array}
\right.
\end{equation*}
The array $\Bfrakq$ with in total $\delta$ scalar valued entries may contain internal variables of any tensorial rank that describe in a homogenized sense the micro-motion of the material due to structural changes on lower scales.
Classical examples for members of $\Bfrakq$ are damage variables, plastic strains or phase fractions. In addition, we assume all tensorial elements of $\Bfrakq$ to be {\it Lagrangian objects}, i.e., they are not affected by rigid body macro-motions superimposed on the current configuration $\Bvarphi_t(\calB)$.  Like in \citet{miehe11,miehe14}, we distinguish between {\it long-range variables} $\Bfrakq^l$ (order parameters) that are governed by PDEs in form of additional balance equations and connected to given length scale parameters, and {\it short-range variables} $\Bfrakq^s$ that are determined by ODEs and represent the standard concept of {\it local internal variables}. 
We summarize both, long- as well as short-range variables in the array
\begin{equation*}
\Bfrakq = ( \Bfrakq^s, \Bfrakq^l ) \, .
\end{equation*}
For the treatment of each long-range micro-motion field, we decompose the boundary of the reference configuration into nonoverlapping parts $\partial \calB_{\Bfrakq}$ and $\partial \calB_{H}$ such that $\partial\calB=\partial \calB_{\Bfrakq} \cup \partial \calB_{H}$. We prescribe the  micro-motion (Dirichlet boundary condition)
\begin{equation}
 \Bfrakq^l(\BX,t) = \bar \Bfrakq^l(\BX) \quad \mbox{on } \partial \calB_{\Bfrakq} \times \calT
\label{micro_motion_1}
\end{equation}
and the {\it micro-tractions} (Neumann boundary condition) on $\partial \calB_{H} \times \calT$ as specified in Section~\ref{Section333} below.  Note that the given micro-motion \req{micro_motion_1} on the Dirichlet boundary is considered to be independent of time. This is because we assume these variables to be ``{\it passive}'' in the sense that the deformation process cannot be driven by externally applying them.

\subsubsection[Section313]{Thermal Driving Force} \label{Section313}
As a third primary field in the thermomechanics of dissipative materials, we introduce the (macroscopic) {\it thermal driving force field}
\begin{equation*}
T : 
\left\{
\begin{array}{ll}
  \calB \times  \calT \rightarrow \mathbb{R}_+ \\
 (\BX, t)  \mapsto T(\BX,t) \, .
\end{array}
\right.
\end{equation*}
As outlined in the motivating Section \ref{Section232}, the thermal driving force appears as the dual quantity to the entropy rate and is identified as the temperature. It is governed by the {\it  generalized Legendre transformation} \req{duality_entropyrate_1} introduced below. With respect to  thermal loading of a heat conducting  solid, we decompose the boundary of the reference configuration into nonoverlapping parts $\partial \calB_T$ and $\partial \calB_q$ such that $\partial \calB = \partial \calB_T \cup \partial \calB_q$. We prescribe the thermal driving force (Dirichlet boundary condition)
\begin{equation}
T(\BX,t) = \bar T(\BX,t) \quad \mbox{on } \partial \calB_T \times \calT
\label{thermal_dirichlet}
\end{equation}
and the {\it heat flux} (Neumann boundary condition) on $\partial \calB_q \times \calT$ as specified in Section~\ref{Section34} below. The primary fields, namely the macro-motion, the long-range micro-motion and, in case of heat conduction, the thermal driving force are depicted in Figure~\ref{fields}.

\subsection[Section33]{Constitutive Framework to Thermomechanics of Gradient-Extended Solids} \label{Section33}

\subsubsection{Free Energy Function}
The free energy storage in continua is governed by a free energy function.
We specify it by focusing on simple materials of grade one, i.e., we include as arguments the first gradients of the micro- and macro-motions\footnote{For the seek of a compact notation, we assume from now on the presence of long-range variables $\Bfrakq = \Bfrakq^l$ only. In case of short-range variables $\Bfrakq^s$, the corresponding gradients vanish and no micro-mechanical boundary conditions are prescribed.}
\begin{equation}
\psi(\BX,t) = \widehat \psi(\Bfrakc,\theta) 
\quad
\mbox{with}
\quad
\Bfrakc = ( \BC, \Bfrakq, \nabla \Bfrakq )
\label{free_energy_1}
\end{equation}
defining the mechanical constitutive state.  Together with the temperature, this state is invariant under any rigid body macro-motion superimposed on the current configuration $\Bvarphi_t(\calB)$, i.e., the function \req{free_energy_1}$_1$ is objective. We define  the energetic
contribution to the (contravariant) first Piola-Kirchhoff stress
tensor $\BP$ as 
\begin{equation*}
\BP^e = 2 \BF \, \partial_{\BC} \widehat \psi \, ,
\end{equation*}
and the well known Clausius-Planck inequality takes the form
\begin{equation}
\calD_{int}= \Bg \BP^d : \dot \BF - \partial_{\Bfrakq} \widehat \psi \cdot \dot \Bfrakq - \partial_{\nabla \Bfrakq} \widehat \psi \cdot \nabla \dot \Bfrakq \geq 0 \, .
\label{free_energy_2}
\end{equation}
Here, $ \BP^d =  \BP -  \BP^e$ denotes the dissipative part of the first Piola-Kirchhoff stress tensor. 

\subsubsection[Section333]{Dissipation Potential Functions} \label{Section333}
Dissipative mechanisms are described by dissipation potential functions. First, we take mechanical effects  into account via an objective intrinsic dissipation potential function
\begin{equation}
\phi_{int}(\BX,t)=\widehat \phi_{int}(\dot \Bfrakc; \Bfrakc,\theta)
\label{dissipation_potential_1}
\end{equation}
at given thermomechanical state $(\Bfrakc,\theta)$. For generality, we assume in this abstract setting $\widehat \phi_{int}$ to be nonsmooth
with respect to $\dot \Bfrakc$, e.g., positively homogeneous of degree one, $\widehat \phi_{int}(\gamma \dot \Bfrakc; \Bfrakc,\theta) = \gamma \widehat \phi_{int}(\dot \Bfrakc;\Bfrakc, \theta)$ for all $\gamma>0$, which characterizes in the adiabatic case a rate-independent evolution, see also Footnote~\ref{f_1homo}. With the intrinsic dissipation potential function at hand, we define the dissipative stress
\begin{equation*}
\BP^d = 2 \BF \, \partial_{\dot \BC} \widehat \phi_{int} \, , 
\label{dissipation_potential_2}
\end{equation*}
that arises in the reduced Clausius-Planck inequality \req{free_energy_2}.  
As suggested in \citet{miehe14}, we assume the evolution equations\footnote{\textbf{Variational Derivative.} Throughout this text, we use the notation of variational derivatives
\begin{equation*}
\delta_z \widehat y(z, \nabla z) = \partial_z \widehat y(z, \nabla z) - \Div [\, \partial_{\nabla z} \widehat y(z, \nabla z)  \, ] \, .
\end{equation*}} for the micro-motions together with the Neumann boundary conditions
\begin{equation}
\delta_{\Bfrakq} \widehat \psi + \delta_{\dot \Bfrakq} \widehat
\phi_{int} \ni \Bzero \mbox{ in } \calB \times \calT \quad \mbox{with}
\quad [\, \partial_{\nabla \Bfrakq} \widehat \psi + \partial_{\nabla
\dot \Bfrakq} \widehat \phi_{int} \,] \, \Bn_0 \ni \Bzero \mbox{ on
} \partial \calB_H \times \calT \, .
\label{dissipation_potential_3}
\end{equation}
In view of the nonsmoothness of the dissipation potential function \req{dissipation_potential_1}, we write the micro- and macroscopic balance equations \req{dissipation_potential_3} and \req{gov_equ_1} as differential inclusions, i.e., $\partial_{\dot \Bfrakc}(\cdot)$ is understood as subdifferential. Considering the global form of \req{free_energy_2}, doing integration by parts two times and taking into account the homogeneous rate-type Dirichlet boundary condition $\dot \Bfrakq = \Bzero$ on $\partial \calB_\Bfrakq$ according to \req{micro_motion_1}, yields the alternative representation
\begin{equation}
\calD_{int} = \partial_{\dot \Bfrakc}\widehat \phi_{int} \cdot \dot \Bfrakc \geq 0
\label{dissipation_potential_5}
\end{equation}
for the intrinsic dissipation, see again \citet{miehe14}. 
Note that the micro-force balance \req{dissipation_potential_3}$_1$ and the homogeneous Neumann boundary condition \req{dissipation_potential_3}$_2$ are outcomes of the variational principle \req{cont_vp_5} set up
below.\footnote{ In the external power functional \req{powerm_1}, we assume vanishing micro-tractions. The reason for that is in the same line as for the Dirichlet boundary condition \req{micro_motion_1}, i.e., we do not consider an external driving of the deformation process by the micro-motion and its force-like quantities, respectively.}

As a second contribution to entropy production, we take into account heat conduction via an objective  dissipation potential
\begin{equation*}
\phi_{con}(\BX,t) = \widehat \phi_{con}(\KIg;\Bfrakc,\theta)
\quad
\mbox{with}
\quad
\KIg = - \nabla \theta / \theta
\end{equation*}
at given thermomechanical state $( \Bfrakc, \theta )$. With such a function at hand, we define the material heat flux vector
\begin{equation*}
\Bq = \partial_{\KIg} \widehat \phi_{con} \,  ,
\label{dissipation_potential_8}
\end{equation*}
and the well known Fourier inequality takes the form
\begin{equation}
\calD_{con} = \partial_{\KIg} \widehat \phi_{con} \cdot \KIg \geq 0 \, .
\label{dissipation_potential_9}
\end{equation}
The inequalities \req{dissipation_potential_5} and \req{dissipation_potential_9} serve as fundamental physically-based constraints on the dissipation potential functions $\widehat \phi_{int}$ and $\widehat \phi_{con}$. These two conditions are satisfied {\it a priori} for dissipation potential functions that are (i) nonnegative $\widehat \phi(\cdot \, ;\Bfrakc,\theta) \geq 0$, (ii) zero in the origin $\widehat \phi(\Bzero;\Bfrakc,\theta) = 0$ and (iii) convex in $\dot \Bfrakc$ and $\KIg$, respectively.

\subsection[Section34]{Governing Equations for Thermomechanics of Gradient-\\Extended Solids} \label{Section34}
We now summarize the governing equations for the thermomechanics of gradient-extended dissipative solids undergoing large deformations. First, we have the balance of linear momentum 
\begin{equation}
\delta_{\Bvarphi} \widehat \psi + \delta_{\dot \Bvarphi} \widehat \phi_{int} \ni \Bg \bar \Bgamma \mbox{ in } \calB \times \calT
\quad
\mbox{with}
\quad
2 \BF  [\, \partial_{\BC} \widehat \psi + \partial_{\dot \BC} \widehat \phi_{int} \, ] \, \Bn_0 \ni  \bar \Bt \mbox{ on } \partial \calB_\Bt \times \calT \, ,
\label{gov_equ_1}
\end{equation}
where $\bar \Bgamma(\BX,t)$ and $\bar \Bt(\BX,t)$ are prescribed body force and nominal surface traction fields. The balance of angular momentum $\BP\BF^T = \BF \BP^T$ is {\it  a priori} satisfied
by the objectivity of the free energy function \req{free_energy_1}
and the dissipation potential function
\req{dissipation_potential_1}.  The evolving micro-motion is governed
by the micro-force balance equation
\begin{equation}
\delta_{\Bfrakq} \widehat \psi + \delta_{\dot \Bfrakq} \widehat
\phi_{int} \ni \Bzero \mbox{ in } \calB \times \calT \quad \mbox{with}
\quad [\, \partial_{\nabla \Bfrakq} \widehat \psi + \partial_{\nabla
\dot \Bfrakq} \widehat \phi_{int} \,] \, \Bn_0 \ni \Bzero \mbox{ on
} \partial \calB_H \times \calT 
\label{gov_equ_2}
\end{equation}
as already given in \req{dissipation_potential_3}. On the
thermal side, the entropy is defined by the local state equation
\begin{equation}
\eta = - \partial_{\theta} \widehat \psi(\Bfrakc,\theta) \quad
\mbox{in } \calB \times \calT \, .
\label{gov_equ_3}
\end{equation}
Finally, the evolution of the entropy is governed by the energy equation
\begin{equation}
-\dot \eta + \frac{1}{\theta} \, \partial_{\dot \Bfrakc} \widehat \phi_{int} \cdot \dot \Bfrakc - \delta_{\theta} \widehat \phi_{con} = - \frac{\bar r}{\theta} 
\mbox{ in }
\calB \times \calT
\quad
\mbox{with}
\quad
\frac{1}{\theta} \, \partial_{\KIg} \widehat \phi_{con} \cdot \Bn_0 = \frac{\bar q}{\theta}
\mbox{ on }
\partial \calB_q \times \calT \, ,
\label{gov_equ_4}
\end{equation}
where $\bar r(\BX,t)$  and $\bar q(\BX,t)$ are prescribed heat source and material surface heat flux fields. Note that if $\widehat \phi_{int}(\cdot \, ; \Bfrakc,\theta)$ is positively homogeneous of degree one\footnote{Since $\widehat \phi_{int}(\cdot \, ; \Bfrakc, \theta)$ is assumed to be nonnegative and convex, positive homogeneity of degree one implies that this function is either nonsmooth at $\dot \Bfrakc = \Bzero$ or zero everywhere.\label{f_1homo}}, $\partial_{\dot \Bfrakc} \widehat \phi_{int}$ is a set, but the expression  $\partial_{\dot \Bfrakc} \widehat \phi_{int} \cdot \dot \Bfrakc$, which represents the intrinsic dissipation, evaluates to a unique value since it coincides with $\widehat \phi_{int}(\dot \Bfrakc;\Bfrakc,\theta)$, see also Section~\ref{Section212}. The equations \req{gov_equ_1}--\req{gov_equ_4} generalize the ODEs
\req{0d_gov_1}--\req{0d_gov_4} for the rheological device to the large-strain
continuum setting and include intrinsic gradient-type dissipative
effects as well as heat conduction.

 \section[Section4]{Variational Principles for Thermomechanics of \\ Gradient-Extended Solids}
\label{Section4}

We now discuss the variational structure of thermomechanics of gradient-extended dissipative solids undergoing large deformations, whose Euler equations were summarized in Section~\ref{Section34}. Again, the starting point is the definition of {\it canonical} energy and dissipation potential functions. By means of (generalized) Legendre transformations, different {\it rate-type} and {\it incremental} mixed formulations are derived.

\subsection[Section41]{Canonical Energy and Dissipation Potential Functionals} \label{Section41}
We generalize the variational framework for the rheological model presented in Section~\ref{Section2} to the large-strain continuum setting of gradient-extended dissipative solids. To this end, based on the internal energy and dissipation potential functions
\begin{equation}
\fterm{
e(\BX,t) = \widehat e(\Bfrakc,\eta)
\quad
\mbox{and}
\quad
v(\BX,t) = \widehat v(\dot \Bfrakc,\dot \eta; \Bfrakc,\eta,\BX,t ) \, ,
}
\label{canonical_energy_1}
\end{equation}
we introduce the energy and dissipation potential functionals
\begin{equation}
E(\Bvarphi,\Bfrakq,\eta) = \int_{\calB} \widehat e(\Bfrakc, \eta) \, dV
\quad
\mbox{and}
\quad
V(\dot \Bvarphi, \dot \Bfrakq, \dot \eta;\Bvarphi,\Bfrakq,\eta,t) = \int_{\calB} \widehat v(\dot \Bfrakc, \dot \eta;\Bfrakc,\eta,\BX,t) \, dV \, .
\label{canonical_energy_2}
\end{equation}
$E$ represents the internal energy stored in the entire body $\calB$ due to coupled micro-macro deformations and thermal effects.  The introduced functional $V$ is related to intrinsic dissipative mechanisms and entropy production due to heat conduction. In analogy to \req{0d_canonical_1}, the entropy $\eta$ as well as the entropy rate $\dot \eta$ are used as canonical thermal variables. On the mechanical side, the constitutive functions \req{canonical_energy_1} are assumed to depend on the constitutive state array $\Bfrakc$ defined in \req{free_energy_1}$_2$ that makes those functions {\it a priori} objective. Note, that the dissipation potential function $\widehat v$ in general depends on the current state $(\Bfrakc, \eta)$ as well as explicitly on position and time $(\BX,t)\in \calB \times \calT$ stemming from a possible inhomogeneous and time-dependent thermal loading, see below. 

\subsection[Section42]{Energy and Dissipation Functionals in terms of Temperature} \label{Section42}
For practical modeling, we transform the above energy and dissipation potential functionals into functionals that depend additionally on the temperature. This we do in analogy to Section~\ref{Section23}.

\subsubsection[Section421]{Variable dual to Entropy} \label{Section421}
First, we define the internal energy functional \req{canonical_energy_2}$_1$  by the Legendre transformation
\begin{equation}
E(\Bvarphi,\Bfrakq, \eta) = \sup_{\theta} \, E^+(\Bvarphi, \Bfrakq, \eta, \theta)
\label{duality_entropy_1}
\end{equation}
in terms of the {\it mixed} internal energy functional
\begin{equation}
E^+(\Bvarphi,\Bfrakq,\eta,\theta) = \int_{\calB} \widehat e^+(\Bfrakc,  \eta, \theta) \, dV 
\WITH
\widehat e^+ = \widehat \psi(\Bfrakc, \theta) + \theta \eta  \, .
\label{duality_entropy_2}
\end{equation}
The necessary optimality condition of \req{duality_entropy_1} is the statement \req{gov_equ_3} which identifies the thermal state variable $\theta$ as the dual quantity to the entropy $\eta$, i.e., as the temperature.

\subsubsection[Section422]{Variable dual to Entropy Rate} \label{Section422}
In a second step, we define at time $t$ the dissipation potential functional by a {\it generalized Legendre transformation}\footnote{The generalized Legendre transformation is conceptually in line with \citet{miehe+hildebrand+boeger14}.}
\begin{equation}
V(\dot \Bvarphi,\dot \Bfrakq, \dot \eta;\Bvarphi,\Bfrakq,\eta,t) = \sup_{T } [\, V^{+}(\dot \Bvarphi,\dot \Bfrakq, \dot \eta, T;\Bvarphi,\Bfrakq,\theta) - P_{ext}^T(T; \theta, t) \, ]
\label{duality_entropyrate_1}
\end{equation}
in terms of the {\it mixed} dissipation potential functional
\begin{equation}
V^+(\dot \Bvarphi,\dot \Bfrakq,\dot \eta, T;\Bvarphi,\Bfrakq,\theta) = \int_{\calB} \widehat v^+(\dot \Bfrakc, \dot \eta, T, \nabla T; \Bfrakc, \theta) \, dV
\WITH
\widehat v^+ = \widetilde \phi(\dot \Bfrakc, T, \nabla T;\Bfrakc,\theta) - T \dot \eta 
\label{duality_entropyrate_2}
\end{equation}
governed by a dissipation potential function $\widetilde \phi$ and the thermal load functional
\begin{equation}
P_{ext}^T(T;\theta,t) =  \int_{\calB} \widehat b_T(T;\theta,\BX,t) \, dV + \int_{\partial \calB_q} \widehat s_T(T;\theta,\BX,t) \, dA \, .
\label{duality_entropyrate_3}
\end{equation}
Note that by use of the local state equation \req{gov_equ_3}, the mixed
dissipation potential functional $V^+$ depends on the current temperature.
The maximization in \req{duality_entropyrate_1} at time $t$ is performed under the constraint $T=\bar T$ on $\partial \calB_T$ and defines as necessary condition the evolution of the entropy along with a thermal boundary condition
\begin{equation}
\dot \eta = \delta_T \widetilde \phi - \partial_T \widehat b_T 
\mbox{ in }
\calB
\quad
\mbox{and}
\quad
\partial_{\nabla T} \widetilde \phi \cdot \Bn_0 = \partial_T \widehat s_T 
\mbox{ on }
\partial \calB_q \, .
\label{duality_entropyrate_4}
\end{equation}
These equations generalize the local statement \req{0d_dual_4} for the rheological device to a large-strain continuum setting including intrinsic gradient-type dissipative effects as well as heat conduction. As in Section \ref{Section232}, we call $T$ the {\it thermal driving force} that is dual to the entropy rate $\dot \eta$. The entropy evolution \req{duality_entropyrate_4}$_1$ must have the form \req{gov_equ_4}$_1$ and we identify
\begin{equation}
\delta_T \widetilde \phi \overset{!}{=}\frac{1}{\theta} \partial_{\dot \Bfrakc} \widehat \phi_{int} \cdot \dot \Bfrakc - \delta_{\theta} \widehat \phi_{con}
\quad
\mbox{and}
\quad
\partial_T \widehat b_T  \overset{!}{=} -\frac{\bar r}{\theta} \, .
\label{duality_entropyrate_5}
\end{equation}
The first of these conditions is fulfilled for $T=\theta$ in $\calB$ if the dissipation potential function is specified as
\begin{equation}
\widetilde \phi(\dot \Bfrakc, T, \nabla T; \Bfrakc,\theta) = \widehat \phi_{int}(\frac{T}{\theta} \dot \Bfrakc; \Bfrakc,\theta) - \widehat \phi_{con}(-\frac{1}{T} \nabla T;\Bfrakc,\theta) \, ,
\label{duality_entropyrate_6}
\end{equation}
which is in line with \citet{yang+stainier+ortiz06}, but derived in an alternative way. The second condition \req{duality_entropyrate_5}$_2$ is satisfied for a volumetric thermal loading function
\begin{equation}
\widehat b_T(T;\theta,\BX,t) = - \frac{T}{\theta} \, \bar r(\BX,t) \, .
\label{duality_entropyrate_7}
\end{equation}
It remains to find an expression for the thermal surface load function $\widehat s_T$. Note that $-\partial _{\nabla T} \widehat \phi_{con}=1/\theta \, \partial_{\KIg} \widehat \phi_{con}$ for $T=\theta$ in $\calB$ and we identify from \req{gov_equ_4}$_2$
\begin{equation*}
\partial_T \widehat s_T \overset{!}{=} \frac{\bar q}{\theta} \, ,
\end{equation*}
which is fulfilled for a thermal surface load function
\begin{equation}
\widehat s_T(T;\theta,\BX,t) = \frac{T}{\theta} \, \bar q(\BX,t) \, .
\label{duality_entropyrate_9}
\end{equation}

\subsubsection{Load Functionals}
Besides the mixed energy and dissipation potential functionals \req{duality_entropy_2} and \req{duality_entropyrate_2}, we have an external thermomechanical load functional
\begin{equation*}
P_{ext}(\dot \Bvarphi,T;\theta,t) = P_{ext}^\Bvarphi(\dot \Bvarphi;t) + P_{ext}^T(T;\theta,t)
\end{equation*}
with decoupled mechanical and thermal contributions. On the mechanical side, we define the load functional
\begin{equation}
 P_{ext}^\Bvarphi(\dot \Bvarphi;t) = \int_{\calB}  \bar \Bgamma(\BX,t) \cdot \dot \Bvarphi \, dV + \int_{\partial \calB_\Bt}  \bar \Bt(\BX,t) \cdot \dot \Bvarphi \, dA
\label{powerm_1}
\end{equation}
in terms of a given body force field $\bar \Bgamma$ and nominal surface traction field $\bar \Bt$. The thermal load functional \req{duality_entropyrate_3} attains with the identifications \req{duality_entropyrate_7} and \req{duality_entropyrate_9} the form
\begin{equation*}
P_{ext}^T(T;\theta,t) = - \int_{\calB} \frac{T}{\theta} \, \bar r(\BX,t) \, dV + \int_{\partial \calB_q} \frac{T}{\theta} \, \bar q (\BX,t) \, dA
\end{equation*}
in terms of a given heat source field $\bar r$ and material surface heat flux $\bar q$.

\subsection[Section43]{Fundamental Mixed Variational Principle for Thermomechanics} \label{Section43}

\subsubsection[Section431]{Rate-type Formulation} \label{Section431}
Based on the internal energy and dissipation potential functionals $E^+$ and $V^+$ and the thermomechanical load functional $P_{ext}$, we define at current time $t$ the rate-type potential\footnote{To keep notation short, we subsequently do not write explicitly the dependence of functions and corresponding functionals on given quantities.}
\begin{equation*}
\Pi^+(\dot \Bvarphi, \dot \Bfrakq, \dot \eta, \dot \theta, T) = \frac{d}{dt}E^+( \Bvarphi, \Bfrakq, \eta, \theta)+V^+(\dot \Bvarphi,\dot \Bfrakq,\dot \eta,T) - P_{ext}(\dot \Bvarphi,T)
\end{equation*}
with given state $(\Bvarphi,\Bfrakq,\eta,\theta)$. We write this potential with its internal and external contributions 
\begin{equation}
\Pi^+(\dot \Bvarphi,\dot \Bfrakq,\dot \eta, \dot \theta,T)=\int_{\calB} \widehat \pi^+(\dot \Bfrakc,\dot \eta,\dot \theta,T,\nabla T) \, dV - P_{ext}(\dot \Bvarphi,T)
\label{cont_vp_2}
\end{equation}
in terms of the internal potential density
\begin{equation*}
\fterm
{
\widehat \pi^+(\dot \Bfrakc,\dot \eta, \dot \theta,T,\nabla T)= \frac{d}{dt}\widehat e^+(\Bfrakc,\eta,\theta) +\widehat v^+(\dot \Bfrakc, \dot \eta, T, \nabla T)   \, .
}
\end{equation*}
Recalling the mixed functions \req{duality_entropy_2}$_2$ and \req{duality_entropyrate_2}$_2$ together with \req{duality_entropyrate_6} and inserting the necessary condition \req{gov_equ_3} on the given thermomechanical state, yields a reduced internal potential density of the form
\begin{equation}
\widehat \pi^+_{red}(\dot \Bfrakc,\dot \eta,T,\nabla T) = \partial_{\Bfrakc} \widehat \psi \cdot \dot \Bfrakc + (\theta-T) \dot \eta +  \widehat \phi_{int}(\frac{T}{\theta}\dot \Bfrakc) - \widehat \phi_{con}(-\frac{1}{T} \nabla T) \, .
\label{cont_vp_4}
\end{equation}
Then, the rates of the macro- and the micro-motion as well as the rate of the entropy and the thermal driving force at current time $t$ are governed by the  variational principle\footnote{ Note carefully that in case of an adiabatic process, the potential is optimized by minimizing with respect to $T$ and maximizing with respect to $\dot \eta$, see the motivating Section~\ref{Section233}.}
\begin{equation}
\{ \dot \Bvarphi, \dot \Bfrakq, \dot \eta, T \} =
\mbox{Arg} \{\;
\substackrel{\dot \Bvarphi, \dot \Bfrakq, \dot \eta }{\mbox{inf}}\,
\substackrel{T}{\mbox{sup}} \;
[ \, 
\int_{\calB} \widehat \pi^+_{red}(\dot \Bfrakc, \dot \eta, T,\nabla T) \, dV - P_{ext}(\dot \Bvarphi,T)\, ] 
\; \}
\, .
\label{cont_vp_5}
\end{equation}
Here, one has to account for the rate forms of the Dirichlet boundary conditions \req{defo_dirichlet} and \req{micro_motion_1} for the macro- and micro-motions, i.e.,
\begin{equation}
\dot \Bvarphi = \dot{\bar \Bvarphi} 
\quad
\mbox{on }
\partial \calB_\Bvarphi
\AND
\dot \Bfrakq = \Bzero 
\quad
\mbox{on }
\partial \calB_{\Bfrakq}
\label{cont_vp_5a}
\end{equation}
 as well as for the Dirichlet boundary condition \req{thermal_dirichlet} for the thermal driving force. The variational principle \req{cont_vp_5} is as stated in \citet{yang+stainier+ortiz06}, but extended by a long-range micro-motion field. By the first variation of the functional \req{cont_vp_2}, we have the necessary optimality conditions 
\begin{equation*}
\delta_{\dot \Bvarphi} \Pi^+ + \delta_{\dot \Bfrakq} \Pi^+ + \delta_{\dot \eta} \Pi^+ \geq 0 \, , 
\quad
\delta_{T} \Pi^+ \leq 0
\end{equation*}
for all admissible test functions $\delta \dot \eta$ and $(\delta \dot \Bvarphi, \delta \dot \Bfrakq, \delta T)$ fulfilling homogeneous forms of the Dirichlet boundary conditions. We get the Euler equations
\begin{equation}
\parbox{14cm}{\centering 
$
\begin{array}{llll}
1.&{\mbox{ Evolving macro-motion}} &
 \delta_{\dot \Bvarphi} \widehat \pi^+_{red} \equiv \delta_{\Bvarphi}\widehat \psi + \delta_{\dot \Bvarphi} \widehat \phi_{int} &\ni \Bg \bar \Bgamma \,  ,   \\[1ex]
2. &{\mbox{ Evolving micro-motion}} & 
\delta_{\dot \Bfrakq} \widehat \pi^+_{red} \equiv \delta_{\Bfrakq} \widehat \psi + \delta_{\dot \Bfrakq} \widehat \phi_{int} &\ni \Bzero \,   ,   \\[1ex]
3. &{\mbox{ Thermal driving force}} &
\partial_{\dot \eta} \widehat \pi_{red}^+ \equiv  T - \theta &= 0 \,  ,   \\[1ex]
4. &{\mbox{ Evolving thermal state}} & 
\delta_{T} \widehat \pi^+_{red} \equiv -\dot \eta + \partial_T \widehat \phi_{int} - \delta_T \widehat \phi_{con}&= - \bar r /\theta
\end{array}
$
}
\label{cont_vp_6}
\end{equation}
in $\calB$ along with the Neumann boundary conditions
\begin{equation}
2 \BF [\, \partial_{\BC} \widehat \psi + \tfrac{T}{\theta} \partial_{\frac{T}{\theta} \dot \BC} \widehat \phi_{int} \, ] \, \Bn_0 \ni \bar \Bt \, , 
\quad
[\, \partial_{\nabla \Bfrakq} \widehat \psi + \tfrac{T}{\theta}\partial_{\frac{T}{\theta}\nabla \dot \Bfrakq} \widehat \phi_{int} \, ] \, \Bn_0 \ni \Bzero \, ,
\quad
- \partial_{\nabla T} \widehat \phi_{con} \cdot \Bn_0 = \bar q/ \theta
\label{cont_vp_7}
\end{equation}
on $\partial \calB_\Bt$, $\partial \calB_{H}$ and $\partial \calB_q$, respectively. In contrast to \req{0d_dual_8}, the equations are now exclusively governed by {\it variational derivatives} of the reduced potential density $\widehat \pi^+_{red}$ defined in \req{cont_vp_4}. The central three field equations are the quasi-static mechanical equilibrium
\begin{equation*}
\delta_{\dot \Bvarphi} \widehat \pi^+_{red} \big \vert_{T=\theta} \equiv -\Div[\, 2\BF( \partial_{\BC}\widehat \psi + \partial_{\dot \BC} \widehat \phi_{int}(\dot \Bfrakc))] \ni \bar \Bgamma \, , 
\end{equation*}
that governs the rate $ \dot \Bvarphi$ of the macro-motion, the micro-force balance
\begin{equation}
\delta_{\dot \Bfrakq} \widehat \pi^+_{red}\big \vert_{T=\theta} \equiv [\, \partial_{\Bfrakq} \widehat \psi + \partial_{\dot \Bfrakq} \widehat \phi_{int}(\dot \Bfrakc)  \,] - \Div[\partial_{\nabla \Bfrakq} \widehat \psi + \partial_{\nabla \dot \Bfrakq} \widehat \phi_{int}(\dot \Bfrakc)] \ni \Bzero
\label{cont_vp_9}
\end{equation}
determining the rate $\dot \Bfrakq$ of the micro-motion and the energy equation
\begin{equation*}
\delta_{T} \widehat \pi^+_{red}\big \vert_{T=\theta} \equiv - \dot \eta + \frac{1}{\theta} \, \partial_{\dot \Bfrakc} \widehat \phi_{int}(\dot \Bfrakc) \cdot \dot \Bfrakc- \frac{1}{\theta} \, \Div[\partial_{\KIg} \widehat \phi_{con}(\KIg)]= - \frac{\bar r}{\theta}
\end{equation*}
for the evolution $\dot \eta$ of the entropy.

\subsubsection[Section432]{Incremental Formulation} \label{Section432}
Consider a finite time interval $[ t_n, t_{n+1}] \subset \overline \calT$ with step length $\tau = t_{n+1} - t_n >0$ and assume all thermomechanical field variables at time $t_n$ to be known. The goal is then to determine all the approximate  fields at time $t_{n+1}$ based on variational principles valid for the time increment under consideration. Subsequently all variables without subscript are understood to be evaluated at time $t_{n+1}$. We may formulate the incremental potential
\begin{equation*}
\Pi^{+ \tau}(\Bvarphi, \Bfrakq, \eta, \theta, T) = E^{+ \tau}(\Bvarphi, \Bfrakq, \eta, \theta) + V^{+ \tau}(\Bvarphi,\Bfrakq,\eta,T) - P_{ext}^\tau(\Bvarphi,T) \, ,
\end{equation*}
where $E^{+ \tau}$, $V^{+ \tau}$ and $P_{ext}^{ \tau}$ are incremental energy, dissipation and load functionals associated with the time interval $[t_n, t_{n+1} ]$. These functionals are defined at given state $(\Bvarphi_n,\Bfrakq_n,\eta_n,\theta_n)$ at time $t_n$. In analogy to the rate-type formulation \req{cont_vp_2}, we rewrite the incremental potential
\begin{equation}
\Pi^{+ \tau}(\Bvarphi, \Bfrakq, \eta, \theta,T) = \int_{\calB} \widehat \pi^{+ \tau}(\Bfrakc, \eta, \theta, T, \nabla T) \, dV - P_{ext}^{\tau}(\Bvarphi,T)
\label{inc_vp_2}
\end{equation}
in terms of an incremental internal potential density $\widehat \pi^{+ \tau}$ which is defined at given state $(\Bfrakc_n,\eta_n,\theta_n)$.  Such a function is obtained by a  numerical integration algorithm
\begin{equation}
\widehat \pi^{+ \tau}(\Bfrakc, \eta, \theta, T, \nabla T) = Algo\{ \ \int_{t_n}^{t_{n+1}} \widehat \pi^{+}(\dot \Bfrakc, \dot \eta, \dot \theta, T, \nabla T) \, dt \, \} \, ,
\label{inc_vp_3}
\end{equation}
that has to be constructed in such a way that the subsequent incremental variational principle gives consistent algorithmic counterparts of the Euler equations \req{cont_vp_6}. In what follows, we construct an implicit as well as a semi-explicit  numerical integration algorithm, compare Section~\ref{Section24}. As a typical example we consider an integration using the approximations of the rates of state quantities
\begin{equation}
\dot \Bfrakc^\tau = (\Bfrakc - \Bfrakc_n)/\tau \, , \quad
\dot \eta ^\tau = (\eta - \eta_n)/\tau \, , \quad
\dot \theta ^\tau = (\theta - \theta_n)/\tau
\label{inc_vp_4}
\end{equation}
and the incremental internal potential density
\begin{equation}
\fterm{
\widehat \pi^{+\tau} = 
\widehat \psi(\Bfrakc, \theta)  + \tau [ \, (\theta_k - T) \dot \eta^\tau +  \eta_n \dot \theta^\tau + \widehat \phi_{int}(\frac{T}{\theta_n} \dot \Bfrakc^\tau;\Bfrakc_n,\theta_n) - \widehat \phi_{con}(-\frac{1}{T} \nabla T;\Bfrakc_n,\theta_n)  \,]
}
\label{inc_vp_5}
\end{equation}
with $k=n+1$ for an implicit  numerical integration algorithm according to \req{0d_inc_3} and $k=n$ for a semi-explicit  numerical integration algorithm according to \req{0d_inc_9}. In \req{inc_vp_5} we dropped terms that are associated with previous time $t_n$. Then, defining the incremental load functional
\begin{equation}
P_{ext}^\tau(\Bvarphi,T;\Bvarphi_n,\theta_n,t_{n+1}) = P_{ext}^{\Bvarphi}(\Bvarphi - \Bvarphi_n;t_{n+1}) + \tau P_{ext}^T(T;\theta_n,t_{n+1})
\label{inc_vp_6}
\end{equation}
 the incremental variational principle
\begin{equation*}
\{ \Bvarphi, \Bfrakq, \eta, \theta, T \} =
\mbox{Arg} \{\;
\substackrel{\Bvarphi, \Bfrakq, \eta }{\mbox{inf}}\,
\substackrel{\theta,T}{\mbox{sup}} \;
\Pi^{+ \tau}(\Bvarphi,\Bfrakq,\eta,\theta,T)
\; \}
\end{equation*}
determines all thermomechanical fields at time $t_{n+1}$. Note, that the optimization has to be done considering the Dirichlet boundary conditions \req{defo_dirichlet}, \req{micro_motion_1} and \req{thermal_dirichlet} at time $t_{n+1}$. The corresponding Euler equations read
\begin{equation}
\parbox{14cm}{\centering 
$
\begin{array}{llll}
1.&{\mbox{ Update macro-motion}} &
 \delta_{\Bvarphi} \widehat \pi^{+ \tau} \equiv \delta_{\Bvarphi}\widehat \psi + \tau \delta_{\Bvarphi} \widehat \phi_{int} &\ni \Bg \bar \Bgamma  \,  , \\[1ex]
2. &{\mbox{ Update micro-motion}} & 
\delta_{\Bfrakq} \widehat \pi^{+ \tau} \equiv \delta_{\Bfrakq} \widehat \psi + \tau \delta_{\Bfrakq} \widehat \phi_{int} &\ni \Bzero  \,  ,  \\[1ex]
3. &{\mbox{ Thermal driving force}} &
\partial_{\eta} \widehat \pi^{+ \tau} \equiv  T - \theta_k &= 0  \,  ,  \\[1ex]
4. &{\mbox{ Current temperature}} & 
\partial_{\theta} \widehat \pi^{+ \tau} \equiv \partial_{\theta} \widehat \psi + \eta_k &= 0 \,  ,    \\[1ex]
5. &{\mbox{ Update entropy}} & 
\delta_T \widehat \pi^{+ \tau} \equiv -(\eta - \eta_n) + \tau \partial_T \widehat \phi_{int} - \tau \delta_T \widehat \phi_{con} &= -\tau \bar r / \theta_n
\end{array}
$
}
\label{inc_vp_8}
\end{equation}
in $\calB$ along with the Neumann boundary conditions
\begin{equation}
2 \BF [\, \partial_{\BC} \widehat \psi + \tfrac{T}{\theta_n}\partial_{\frac{T}{\theta_n}\dot \BC^\tau} \widehat \phi_{int} \, ] \, \Bn_0 \ni \bar \Bt \, , 
\quad
[\, \partial_{\nabla \Bfrakq} \widehat \psi + \tfrac{T}{\theta_n}\partial_{\frac{T}{\theta_n}\nabla \dot \Bfrakq^\tau} \widehat \phi_{int} \, ] \, \Bn_0 \ni \Bzero \, ,
\quad
-\partial_{\nabla T} \widehat \phi_{con} \cdot \Bn_0 = \frac{\bar q}{\theta_n}
\label{inc_vp_9}
\end{equation}
on $\partial \calB_\Bt$, $\partial \calB_{H}$ and $\partial \calB_q$, respectively.  As a {\it  fundamental difference} to the fully implicit  numerical algorithm, a semi-explicit update identifies the thermal driving force with the given temperature at time $t_n$. Hence, the scaling factor results in $T/\theta_n = 1$ in $\calB$ such that the  algorithmically consistent form of the intrinsic dissipation defined in \req{free_energy_2} and reformulated in \req{dissipation_potential_5} is obtained. Especially, the incremental energy equation reads
\begin{equation*}
\eta = \eta_n + \frac{\tau}{\theta_n} \partial_{\dot \Bfrakc^\tau} \widehat \phi_{int}(\dot \Bfrakc^\tau) \cdot \dot \Bfrakc^\tau + \frac{\tau}{\theta_n} (\bar r - \Div[\partial_{-\frac{1}{T} \nabla T} \widehat \phi_{con}]\big \vert_{T = \theta_n}) \, .
\end{equation*}
Additionally,  also the dissipative terms in the quasi-static equilibrium \req{inc_vp_8}$_1$, micro-force balance \req{inc_vp_8}$_2$ and the boundary conditions \req{inc_vp_9}$_{1-2}$ do not contain the scaling factor.  
However, there are two issues that arise when using the  proposed semi-explicit update for a {\it heat conduction process}: (i) on the thermal side it might be restricted to homogeneous Neumann boundary conditions \req{inc_vp_9}$_3$ on the whole boundary, i.e., $\bar q = 0$ on $\partial \calB$ and (ii) the Courant-Friedrichs-Lewy (CFL) condition gives, depending on the mesh size associated with the spatial discretization, an upper bound for the time step size $\tau$ in order to obtain a stable numerical solution.\footnote{ A derivation of the CFL condition for the suggested semi-explicit update scheme in case of a heat conduction process is beyond the scope of this paper. For a stability analysis of fractional step methods in thermomechanics we refer to  \citet{armero+simo92}.}
Note that the semi-explicit update can be seen as an {\it incrementally isentropic operator split} that consists of two {\it fractional steps}
\begin{equation*}
Algo = Algo_{\eta,T} \circ Algo_{\Bvarphi,\Bfrakq,\theta} \, .
\end{equation*}
First, in the {\it isentropic predictor step} we optimize the incremental potential \req{inc_vp_2} with respect to the macro- and micro-motions $\Bvarphi$ and $\Bfrakq$ as well as the temperature $\theta$, i.e.,
\begin{equation*}
(Algo_{\Bvarphi, \Bfrakq, \theta}):
\quad
\{ \Bvarphi^*, \Bfrakq^*, \theta^* \} =
\mbox{Arg} \{\;
\substackrel{\Bvarphi, \Bfrakq, \theta }{\mbox{stat}} \;
\Pi^{+ \tau}(\Bvarphi,\Bfrakq,\eta_n,\theta,\theta_n)
\; \} \, ,
\end{equation*}
where the entropy is frozen. 
Then, the entropy $\eta$ and the thermal driving force $T$ are updated via the {\it entropy corrector step} 
\begin{equation*}
(Algo_{\eta, T}):
\quad
\{ \eta^*, T^* \} =
\mbox{Arg} \{\;
\substackrel{\eta, T }{\mbox{stat}} \;
\Pi^{+ \tau}(\Bvarphi^*,\Bfrakq^*,\eta,\theta^*,T) 
\; \} \, .
\end{equation*}

\subsection[Section44]{Mixed Variational Principle with Mechanical Driving Forces} \label{Section44}
 We now put the focus on {\it mixed variational principles} for the thermomechanics of gradient-extended dissipative continua, where not only the microstructural variables, but also the corresponding local driving forces are taken into account and considered as additional variables. Following \citet{miehe11,miehe14}, we consider the equivalent representation of the intrinsic dissipation  \req{dissipation_potential_5}
\begin{equation*}
\calD_{int} = \partial_{\dot \Bfrakc} \widehat \phi_{int}(\dot \Bfrakc; \Bfrakc, \theta) \cdot \dot \Bfrakc \iff \calD_{int}= \Bfrakf \cdot \partial_{\Bfrakf} \widehat \phi_{int}^*(\Bfrakf;\Bfrakc,\theta)
\end{equation*}
by the dual intrinsic dissipation potential function $\widehat \phi_{int}^*$ depending on the mechanical driving forces
\begin{equation}
\Bfrakf = (\Bfrakm, \Bfrakd, \Bfrakg ) 
\quad
\mbox{conjugate to}
\quad
\Bfrakc = (\BC, \Bfrakq, \nabla \Bfrakq) \, .
\label{mixed_mech_2}
\end{equation}
The Legendre transformation
\begin{equation}
\widehat \phi_{int}(\dot \Bfrakc;\Bfrakc,\theta) = \sup_{\Bfrakf} [\,  \Bfrakf \cdot \dot \Bfrakc - \widehat \phi_{int}^*(\Bfrakf;\Bfrakc,\theta)  \,]
\label{mixed_mech_3}
\end{equation}
motivates the definition of an extended {\it mixed} dissipation potential functional
\begin{equation}
V^*(\dot \Bvarphi,\dot \Bfrakq, \dot \eta, T, \Bfrakf; \Bvarphi, \Bfrakq, \theta) = \int_{\calB} \widehat v^*(\dot \Bfrakc, \dot \eta, T, \nabla T, \Bfrakf; \Bfrakc,\theta) \, dV
\label{mixed_mech_4}
\end{equation}
in terms of the {\it mixed} dissipation potential function
\begin{equation}
\widehat v^* = \frac{T}{\theta} \Bfrakf \cdot \dot \Bfrakc - T\dot\eta - \widehat \phi_{int}^*(\Bfrakf;\Bfrakc,\theta) - \widehat \phi_{con}(-\frac{1}{T} \nabla T; \Bfrakc,\theta) \,  ,
\label{mixed_mech_5}
\end{equation}
which governs the subsequent extended mixed variational principle. The necessary condition of \req{mixed_mech_3} reads\footnote{\textbf{Perzyna-type dual dissipation potential function.} An important example of a {\it smooth} intrinsic dual dissipation potential function is
\begin{equation}
\widehat \phi_{int}^*(\Bfrakf;\Bfrakc, \theta) = \frac{1}{2  \eta_f} \langle \,  \widehat f(\Bfrakf; \Bfrakc, \theta) \,  \rangle_+^2
\label{perzyna_1}
\end{equation}
in terms of a function $\widehat f(\cdot \, ; \Bfrakc, \theta)$ that differs from a {\it gauge} just by an additive constant and serves as a threshold function in the (adiabatically) rate-independent setting considered in Section~\ref{Section45} below. $\eta_f > 0$ is a material parameter and $\langle \cdot \rangle_+ \colon \mathbb{R} \rightarrow \mathbb{R}_+$, $x \mapsto \frac{1}{2}(|x|+x)$ the ramp function. Note carefully that $\widehat \phi_{int}^*(\cdot \, ; \Bfrakc, \theta)$ defined in \req{perzyna_1} is not a homogeneous function and its dual $\widehat \phi_{int}(\cdot \, ; \Bfrakc, \theta)$ is still nonsmooth in general, i.e.,
\begin{equation}
\dot \Bfrakc = \partial_{\Bfrakf} \widehat \phi^*_{int}
\quad
\mbox{but}
\quad
\delta_{\Bvarphi} \widehat \psi + \delta_{\dot \Bvarphi} \widehat \phi_{int} \ni \Bg \bar \Bgamma \AND
\delta_{\Bfrakq} \widehat \psi + \delta_{\dot \Bfrakq} \widehat \phi_{int} \ni \Bzero \, .
\label{perzyna_2}
\end{equation}
We can write \req{perzyna_2}$_1$ in the particular form 
\begin{equation*}
\dot \Bfrakc = \lambda \, \partial_{\Bfrakf} \widehat f
\WITH
\lambda = \frac{1}{\eta_f} \langle \, \widehat f \, \rangle_+ \, ,
\end{equation*}
which regularizes the rate-independent structure \req{threshold_4}--\req{threshold_5} to be discussed later in Section~\ref{Section45}. This rate-dependent setting is in line with the formulation of visco-plasticity according to \citet{perzyna66}.
\label{footnote_perzyna}
}
\begin{equation*}
 \dot \Bfrakc \in \partial_{\Bfrakf} \widehat \phi_{int}^*(\Bfrakf;\Bfrakc,\theta) \, .
\end{equation*}

\subsubsection{Rate-type Formulation} Based on the internal energy and dissipation potential functionals $E^+$ in \req{duality_entropy_2} and $V^*$ in \req{mixed_mech_4}, we are in the position to formulate a mixed rate-type variational principle that accounts for the mechanical driving forces $\Bfrakf$. We define at current time $t$ the rate-type potential\footnote{To keep notation short, we subsequently do not write explicitly the dependence of functions and corresponding functionals on given states.}
\begin{equation*}
\Pi^*(\dot \Bvarphi, \dot \Bfrakq, \dot \eta, \dot \theta, T ,\Bfrakf) = \frac{d}{dt}E^+(\Bvarphi, \Bfrakq, \eta, \theta) + V^*(\dot \Bvarphi, \dot \Bfrakq, \dot \eta, T, \Bfrakf) - P_{ext}(\dot \Bvarphi,T)
\end{equation*}
with given state $(\Bvarphi, \Bfrakq, \eta, \theta)$. We write this potential with its internal and external contributions
\begin{equation}
\Pi^*(\dot \Bvarphi, \dot \Bfrakq, \dot \eta, \dot \theta, T ,\Bfrakf) = \int_{\calB} \widehat \pi^*(\dot \Bfrakc, \dot \eta, \dot \theta, T, \nabla T, \Bfrakf) \, dV - P_{ext}(\dot \Bvarphi,T)
\label{extended_cont_2}
\end{equation}
in terms of the extended internal potential density
\begin{equation*}
\fterm{
\widehat \pi^*(\dot \Bfrakc, \dot \eta, \dot \theta, T, \nabla T, \Bfrakf) = \frac{d}{dt}\widehat e^+(\Bfrakc,\eta,\theta)+ \widehat v^*(\dot \Bfrakc, \dot \eta, T, \nabla T, \Bfrakf)  \, .
}
\end{equation*}
Recalling the mixed functions \req{duality_entropy_2}$_2$ and \req{mixed_mech_5} and inserting the necessary condition \req{gov_equ_3} on the given thermomechanical state, yields a reduced internal potential density of the form
\begin{equation}
\widehat \pi_{red}^*(\dot \Bfrakc, \dot \eta, T, \nabla T, \Bfrakf) = \partial_{\Bfrakc} \widehat \psi \cdot \dot \Bfrakc+(\theta-T)\dot \eta + \frac{T}{\theta} \Bfrakf \cdot \dot \Bfrakc - \widehat \phi_{int}^*(\Bfrakf) - \widehat \phi_{con}(-\frac{1}{T} \nabla T) \, .
\label{extended_cont_4}
\end{equation}
Then, the rates of the macro- and micro-motion as well as the rate of the entropy and the thermal and mechanical driving forces at current time $t$ are governed by the variational principle
\begin{equation*}
\{ \dot \Bvarphi, \dot \Bfrakq, \dot \eta, T, \Bfrakf \} =
\mbox{Arg} \{\;
\substackrel{\dot \Bvarphi, \dot \Bfrakq, \dot \eta }{\mbox{inf}}\,
\substackrel{T,\Bfrakf}{\mbox{sup}} \;
[ \, 
\int_{\calB} \widehat \pi^*_{red}(\dot \Bfrakc, \dot \eta, T,\nabla T,\Bfrakf) \, dV - P_{ext}(\dot \Bvarphi,T)\, ] 
\; \}
\, .
\end{equation*}
Like in \req{cont_vp_5}, one has to account for the Dirichlet boundary conditions \req{cont_vp_5a} and \req{thermal_dirichlet}.  By the first variation of the functional \req{extended_cont_2}, we have the necessary optimality conditions
\begin{equation*}
\delta_{\dot \Bvarphi} \Pi^* + \delta_{\dot \Bfrakq} \Pi^* + \delta_{\dot \eta} \Pi^* \geq 0 \, , 
\quad
\delta_{T} \Pi^* + \delta_{\Bfrakf} \Pi^* \leq 0
\end{equation*}
for all admissible test functions $(\delta \dot \eta, \delta \Bfrakf)$ and $(\delta \dot \Bvarphi, \delta \dot \Bfrakq, \delta T)$ fulfilling homogeneous forms of the Dirichlet boundary conditions. We obtain the Euler equations
\begin{equation}
\parbox{14cm}{\centering 
$
\begin{array}{llll}
1.&{\mbox{ Evolving macro-motion}} &
 \delta_{\dot \Bvarphi} \widehat \pi^*_{red} \equiv \delta_{\Bvarphi}\widehat \psi + \delta_{\dot \Bvarphi} (\frac{T}{\theta} \Bfrakf \cdot \dot \Bfrakc) &= \Bg \bar \Bgamma  \,  ,   \\[1ex]
2. &{\mbox{ Evolving micro-motion}} & 
\delta_{\dot \Bfrakq} \widehat \pi^*_{red} \equiv \delta_{\Bfrakq} \widehat \psi + \delta_{\dot \Bfrakq} (\frac{T}{\theta} \Bfrakf \cdot \dot \Bfrakc) &= \Bzero  \,  ,   \\[1ex]
3. &{\mbox{ Thermal driving force}} &
\partial_{\dot \eta} \widehat \pi_{red}^* \equiv  T - \theta &= 0  \,  ,   \\[1ex]
4. &{\mbox{ Evolving thermal state}} & 
\delta_{T} \widehat \pi^*_{red} \equiv -\dot \eta + \partial_T (\frac{T}{\theta} \Bfrakf \cdot \dot \Bfrakc) - \delta_T \widehat \phi_{con}&= - \bar r /\theta \, ,  \\[1ex]
5. &{\mbox{ Mechanical driving forces}} & 
\partial_{\Bfrakf} \widehat \pi^*_{red} \equiv \frac{T}{\theta} \dot \Bfrakc - \partial_{\Bfrakf} \widehat \phi_{int}^* & \ni \Bzero
\end{array}
$
}
\label{extended_cont_6}
\end{equation}
in $\calB$ along with the Neumann boundary conditions
\begin{equation}
2 \BF [\, \partial_{\BC} \widehat \psi + \tfrac{T}{\theta}\Bfrakm \, ] \, \Bn_0 = \bar \Bt \, , 
\quad
[\, \partial_{\nabla \Bfrakq} \widehat \psi + \tfrac{T}{\theta} \Bfrakg \, ] \, \Bn_0 = \Bzero \, ,
\quad
-\partial_{\nabla T} \widehat \phi_{con} \cdot \Bn_0 = \bar q/\theta
\label{extended_cont_7}
\end{equation}
on $\partial \calB_\Bt$, $\partial \calB_{H}$ and $\partial \calB_q$, respectively. The central three field equations are the quasi-static equilibrium
\begin{equation}
 \delta_{\dot \Bvarphi} \widehat \pi^*_{red}\big \vert_{T=\theta} \equiv- \Div[\,2\BF( \partial_{\BC} \widehat \psi +  \Bfrakm) \,] =  \bar \Bgamma \, ,
\label{extended_cont_8}
\end{equation}
the micro-force balance 
\begin{equation}
\delta_{\dot \Bfrakq} \widehat \pi^*_{red}\big \vert_{T=\theta} \equiv [ \, \partial_{\Bfrakq} \widehat \psi + \Bfrakd  \,] - \Div[ \, \partial_{\nabla \Bfrakq} \widehat \psi + \Bfrakg  \,] = \Bzero
\label{extended_cont_9}
\end{equation}
and the energy equation
\begin{equation}
\delta_{T} \widehat \pi^*_{red} \big \vert_{T=\theta} \equiv -\dot \eta + \frac{1}{\theta} \, \Bfrakf \cdot \dot \Bfrakc - \frac{1}{\theta} \Div[ \, \partial_{\KIg}\widehat \phi_{con}(\KIg) \, ]= - \bar r/\theta \, .
\label{extended_cont_10}
\end{equation}
They are complemented by the inverse definitions \req{extended_cont_6}$_5$ of the mechanical driving forces which split into three evolution equations
\begin{equation}
\dot \BC \in \partial_{\Bfrakm} \widehat \phi_{int}^* \, , 
\quad
\dot \Bfrakq \in \partial_{\Bfrakd} \widehat \phi_{int}^* \, ,
\quad
\nabla \dot \Bfrakq \in \partial_{\Bfrakg} \widehat \phi_{int}^* \, ,
\label{extended_cont_10e}
\end{equation}
as already given in \citet{miehe14}, where the identity $T=\theta$ in $\calB$ from \req{extended_cont_6}$_3$ has been used. Note carefully that the driving forces $\Bfrakf$ are variables and the variational derivatives in \req{extended_cont_6}$_{1-2}$ with respect to $\dot \Bvarphi$ and $\dot \Bfrakq$ are understood in the ordinary sense. Hence, \req{extended_cont_6}$_{1-2}$ as well as the boundary conditions \req{extended_cont_7}$_{1-2}$ do not represent differential inclusions. 

\subsubsection[Section443]{Incremental Formulation} \label{Section443}
Within a time interval $[t_n,t_{n+1}] \subset \overline \calT$ a variational principle can be constructed by the same avenue as outlined in Section~\ref{Section432}. 
Using the algorithmic approximations \req{inc_vp_4} of rates of state quantities, we get the incremental internal potential density
\begin{equation}
\fterm{
\widehat \pi^{* \tau} =
\widehat \psi(\Bfrakc, \theta)  + \tau [ \, (\theta_k - T) \dot \eta^\tau +  \eta_n \dot \theta^\tau + \frac{T}{\theta_n} \, \Bfrakf \cdot \dot \Bfrakc^\tau - \widehat \phi^*_{int}(\Bfrakf;\Bfrakc_n,\theta_n) - \widehat \phi_{con}(-\frac{1}{T} \nabla T;\Bfrakc_n, \theta_n)  \,] \,  , 
}
\label{extended_inc_2}
\end{equation}
where again $k=n+1$ corresponds to a fully implicit and $k=n$ to a semi-explicit update. In addition, we dropped terms that are associated with previous time $t_n$. Then, with the use of the incremental load functional \req{inc_vp_6} the incremental variational principle 
\begin{equation}
\{ \Bvarphi, \Bfrakq, \eta, \theta, T, \Bfrakf \} =
\mbox{Arg} \{\;
\substackrel{\Bvarphi, \Bfrakq, \eta }{\mbox{inf}}\,
\substackrel{\theta,T,\Bfrakf}{\mbox{sup}} \;
[ \, 
\int_{\calB} \widehat \pi^{* \tau}(\Bfrakc, \eta, \theta, T,\nabla T,\Bfrakf) \, dV - P_{ext}^\tau(\Bvarphi,T)\, ] 
\; \}
\label{extended_inc_3}
\end{equation}
determines all thermomechanical fields at time $t_{n+1}$. The corresponding Euler equations in $\calB$ are time-discrete forms of \req{extended_cont_6} together with \req{inc_vp_8}$_4$ stemming from the variation with respect to the temperature $\theta$. As before, the  proposed semi-explicit integration of the rate of the internal energy yields for the scaling factor $T/\theta_n = 1$ in $\calB$ and one obtains the  algorithmically consistent form of the intrinsic dissipation, i.e., the incremental energy equation reads
\begin{equation*}
\eta = \eta_n + \frac{\tau}{\theta_n} \, \Bfrakf \cdot \dot \Bfrakc^\tau + \frac{\tau}{\theta_n}(\bar r - \Div[\partial_{- \frac{1}{T} \nabla T} \widehat \phi_{con}] \vert_{T = \theta_n}) \, .
\end{equation*}
In addition, the dissipative terms in the time-discrete forms of the quasi-static equilibrium \req{extended_cont_6}$_1$, the micro-force balance \req{extended_cont_6}$_2$ and the boundary conditions \req{extended_cont_7}$_{1-2}$ as well as the dissipative terms in the time discrete forms of the evolution equations \req{extended_cont_6}$_5$ do not contain the scaling factor. The isentropic operator split is modified by an additional optimization in the isentropic predictor step with respect to the mechanical driving forces $\Bfrakf$.

\subsection[Section45]{Mixed Variational Principle with Threshold Function} \label{Section45}
Intrinsic dissipation potential functions are often modelled by the {\it principle of maximum  dissipation}. For the classical local theories of plasticity, this principle can be traced back among others to the work of \citet{hill50}, see also \citet{moreau76}, \citet{simo88} and  \citet{lubliner08}. For a general discussion of this principle and its connection to evolution laws governed by dissipation potentials we refer to \citet{hackl+fischer07} and \citet{hackl+fischer+svoboda10,hackl+fischer+svoboda_ad11}. 
The intrinsic dissipation potential function  for thermomechanics is defined by the constrained maximum principle
\begin{equation}
\widehat \phi_{int}(\dot \Bfrakc; \Bfrakc,\theta) = \sup_{\Bfrakf \in \mathbb{E}(\Bfrakc,\theta)} \Bfrakf \cdot \dot \Bfrakc \, ,
\label{threshold_1}
\end{equation}
that includes at given  state $(\Bfrakc,\theta)$ a domain of admissible mechanical driving forces which we assume to be smooth\footnote{For example, in case of multisurface plasticity one considers a set of $m \geq 2$ intersecting  functions $\widehat f_\alpha$, $\alpha=1,...,m$  which generate a convex but maybe nonsmooth domain $\mathbb{E}$, see, e.g., \citet{simo+kennedy+govindjee88}.  }
\begin{equation}
\mathbb{E}(\Bfrakc,\theta) = \{ \, \Bfrakf \, \vert \, \widehat f(\Bfrakf; \Bfrakc, \theta) \leq 0  \, \} \, .
\label{threshold_2}
\end{equation}
Clearly, the function $\widehat \phi_{int}(\cdot \, ; \Bfrakc, \theta)$ defined in \req{threshold_1} is positively homogeneous of degree one. The set \req{threshold_2} is governed by a threshold function $\widehat f(\Bfrakf; \Bfrakc,\theta) = \widehat w(\Bfrakf; \Bfrakc, \theta) - \widehat c(\Bfrakc,\theta)$, where $\widehat c(\Bfrakc,\theta) > 0$ is a positive threshold constant that might depend on the given thermomechanical state and $\widehat w$ a level set function that is a {\it gauge}, i.e., (i) nonnegative $\widehat w(\cdot \, ; \Bfrakc, \theta) \geq 0$, (ii) zero in the origin $\widehat w(\Bzero;\Bfrakc,\theta) = 0 $, (iii) convex in $\Bfrakf$ and (iv) positively homogeneous of degree one in $\Bfrakf$. As a result of (iii), the constrained optimization problem \req{threshold_1} has a {\it unique} solution. By the use of the Lagrange multiplier method, we can put \req{threshold_1} into the form
\begin{equation}
\widehat \phi_{int}(\dot \Bfrakc; \Bfrakc, \theta) = \sup_{\Bfrakf} \inf_{\lambda \geq 0} \, [\, \Bfrakf \cdot \dot \Bfrakc - \lambda \widehat f(\Bfrakf; \Bfrakc,\theta)  \, ] \,  ,
\label{threshold_3}
\end{equation}
 whose necessary optimality condition defines the evolution of the constitutive state 
\begin{equation}
\dot \Bfrakc = \lambda \, \partial_{\Bfrakf} \widehat f \, ,
\label{threshold_4}
\end{equation}
where the Lagrange multiplier $\lambda$ satisfies classical loading-unloading conditions in Kuhn-Tucker form
\begin{equation}
\lambda \geq 0 \, , 
\quad
\widehat f(\Bfrakf;\Bfrakc,\theta) \leq 0 \, ,
\quad
\lambda \, \widehat f(\Bfrakf;\Bfrakc,\theta) = 0 \, .
\label{threshold_5}
\end{equation}
This is the typical structure of flow rules associated with rate-independent material behavior. Note that even though $\partial_\Bfrakf \widehat f (\Bzero; \Bfrakc, \theta)$ is a subdifferential (see Footnote~\ref{f_1homo}), we get for  $\Bfrakf = \Bzero$ the unique value $\dot \Bfrakc = \Bzero$.
The optimization problem \req{threshold_3} now motivates the definition of a modified {\it mixed} dissipation potential functional
\begin{equation}
V_{\lambda}^*(\dot \Bvarphi,\dot \Bfrakq,\dot \eta, T, \Bfrakf, \lambda; \Bvarphi, \Bfrakq, \theta) = \int_{\calB} \widehat v^{*}_{\lambda}(\dot \Bfrakc, \dot \eta, T, \nabla T, \Bfrakf, \lambda; \Bfrakc,\theta) \, dV
\label{threshold_6}
\end{equation}
in terms of the {\it mixed} dissipation potential function
\begin{equation}
\widehat v^{*}_{\lambda} = \frac{T}{\theta}  \Bfrakf \cdot \dot \Bfrakc - T \dot \eta - \lambda  \widehat f(\Bfrakf;\Bfrakc,\theta) - \widehat \phi_{con}(-\frac{1}{T} \nabla T; \Bfrakc,\theta) \, ,
\label{threshold_7}
\end{equation}
that governs the subsequent modified mixed variational principle.

It should be mentioned that unlike for local theories, the Lagrange multiplier $\lambda$ for $\widehat f = 0$ cannot be determined at current time by a local consistency condition  in {\it terms of rates of the external quantities deformation and temperature}. Hence, as, e.g., mentioned in \citet{deborst+muehlhaus92} a nonlocal version has to be elaborated, see Section~\ref{Section451}.

\subsubsection[Section451]{Rate-type Formulation} \label{Section451}
Based on the internal energy and dissipation potential functionals $E^+$ in \req{duality_entropy_2} and $V^*_{\lambda}$ in \req{threshold_6}, we are in the position to formulate a mixed rate-type variational principle that accounts for dissipative threshold mechanisms. We define at current time $t$ the rate-type potential
\begin{equation*}
\Pi^*_\lambda(\dot \Bvarphi,\dot \Bfrakq, \dot \eta, \dot \theta,T, \Bfrakf, \lambda) = \frac{d}{dt} E^+(\Bvarphi,\Bfrakq,\eta,\theta) + V^*_{\lambda}(\dot \Bvarphi, \dot \Bfrakq, \dot \eta, T, \Bfrakf, \lambda) - P_{ext}(\dot \Bvarphi, T)
\label{rate_vp_threshold_1}
\end{equation*}
with given state $(\Bvarphi,\Bfrakq, \eta, \theta)$. We write this potential with its internal and external contributions
\begin{equation}
\Pi^*_{\lambda}(\dot \Bvarphi, \dot \Bfrakq, \dot \eta, \dot \theta,T,\Bfrakf,\lambda) = \int_{\calB} \widehat \pi^*_{\lambda}(\dot \Bfrakc, \dot \eta, \dot \theta,T,\nabla T, \Bfrakf,\lambda) \, dV - P_{ext}(\dot \Bvarphi,T)
\label{rate_vp_threshold_2}
\end{equation}
in terms of the extended internal potential density
\begin{equation*}
\fterm{
\widehat \pi^*_{\lambda}(\dot \Bfrakc, \dot \eta, \dot \theta,T,\nabla T, \Bfrakf,\lambda) =  \frac{d}{dt}\widehat e^+(\Bfrakc,\eta,\theta)+ \widehat v^*_\lambda(\dot \Bfrakc, \dot \eta, T, \nabla T, \Bfrakf, \lambda) \, .
}
\label{rate_vp_threshold_3}
\end{equation*}
Recalling the mixed functions \req{duality_entropy_2}$_2$ and \req{threshold_7} and inserting the necessary condition \req{gov_equ_3} on the given thermomechanical state, yields a reduced internal potential density of the form
\begin{equation}
\widehat \pi^*_{\lambda,red}(\dot \Bfrakc, \dot \eta,T,\nabla T, \Bfrakf,\lambda) = \partial_{\Bfrakc} \widehat \psi \cdot \dot \Bfrakc+(\theta-T)\dot \eta + \frac{T}{\theta} \Bfrakf \cdot \dot \Bfrakc - \lambda \widehat f(\Bfrakf) - \widehat \phi_{con}(-\frac{1}{T} \nabla T) \, .
\label{rate_vp_threshold_4}
\end{equation}
Then, the rates of the macro- and micro-motion as well as the rate of the entropy, the thermal and mechanical driving forces and the Lagrange multiplier at current time $t$ are governed by the variational principle 
\begin{equation}
\{ \dot \Bvarphi, \dot \Bfrakq, \dot \eta, T, \Bfrakf, \lambda \} =
\mbox{Arg} \{\;
\substackrel{\dot \Bvarphi, \dot \Bfrakq, \dot \eta }{\mbox{inf}}\,
\substackrel{T,\Bfrakf}{\mbox{sup}} \,
\substackrel{\lambda \geq 0}{\mbox{inf}} \;
[ \, 
\int_{\calB} \widehat \pi^*_{\lambda,red}(\dot \Bfrakc, \dot \eta, T,\nabla T,\Bfrakf,\lambda) \, dV - P_{ext}(\dot \Bvarphi,T)\, ] 
\; \}
\, .
\label{rate_vp_threshold_5}
\end{equation}
Like in \req{cont_vp_5}, one has to account for the Dirichlet boundary conditions \req{cont_vp_5a} and \req{thermal_dirichlet}. By the first variation of the functional \req{rate_vp_threshold_2} we have the necessary optimality conditions
\begin{equation}
\delta_{\dot \Bvarphi} \Pi^*_\lambda + \delta_{\dot \Bfrakq} \Pi^*_\lambda + \delta_{\dot \eta} \Pi^*_\lambda \geq 0 \, , 
\quad
\delta_{T} \Pi^*_\lambda + \delta_{\Bfrakf} \Pi^*_\lambda  \leq 0 \, ,
\quad
\delta_{\lambda} \Pi^*_\lambda \geq 0
\label{rate_vp_threshold_6}
\end{equation}
for all admissible test functions $(\delta \dot \eta, \delta \Bfrakf, \delta \lambda)$ with $\lambda + \delta \lambda \geq 0$ in $\calB$ and $(\delta \dot \Bvarphi, \delta \dot \Bfrakq, \delta T)$ fulfilling homogeneous forms of Dirichlet boundary conditions. We obtain the Euler equations
\begin{equation}
\parbox{14cm}{\centering 
$
\begin{array}{llll}
1.&{\mbox{ Evolving macro-motion}} &
 \delta_{\dot \Bvarphi} \widehat \pi^*_{\lambda,red} \equiv \delta_{\Bvarphi}\widehat \psi + \delta_{\dot \Bvarphi} (\frac{T}{\theta} \Bfrakf \cdot \dot \Bfrakc) &= \Bg \bar \Bgamma  \,  ,    \\[1ex]
2. &{\mbox{ Evolving micro-motion}} & 
\delta_{\dot \Bfrakq} \widehat \pi^*_{\lambda,red} \equiv \delta_{\Bfrakq} \widehat \psi + \delta_{\dot \Bfrakq} (\frac{T}{\theta} \Bfrakf \cdot \dot \Bfrakc) &= \Bzero     \,  ,   \\[1ex]
3. &{\mbox{ Thermal driving force}} &
\partial_{\dot \eta} \widehat \pi_{\lambda,red}^* \equiv  T - \theta &= 0 \,  ,   \\[1ex]
4. &{\mbox{ Evolving thermal state}} & 
\delta_{T} \widehat \pi^*_{\lambda,red} \equiv -\dot \eta + \partial_T (\frac{T}{\theta} \Bfrakf \cdot \dot \Bfrakc) - \delta_T \widehat \phi_{con}&= - \bar r /\theta  \,  ,    \\[1ex]
5. &{\mbox{ Mechanical driving forces}} & 
\partial_{\Bfrakf} \widehat \pi^*_{\lambda,red} \equiv \frac{T}{\theta} \dot \Bfrakc - \lambda \, \partial_{\Bfrakf}\widehat f & = \Bzero  \,  ,   \\[1ex]
6. &{\mbox{ Loading Conditions}} & 
\partial_{\lambda} \widehat \pi^*_{\lambda,red} \equiv -\widehat f \geq 0 \, , \, \lambda \geq 0 \, , \, \lambda \widehat f = 0  &
\end{array}
$
}
\label{rate_vp_threshold_7}
\end{equation}
in $\calB$ along with the Neumann boundary conditions \req{extended_cont_7}. Note, that the condition $\widehat f \leq 0$ in \req{rate_vp_threshold_7}$_6$ follows from \req{rate_vp_threshold_6}$_3$ if we set $\lambda = 0$, which necessarily demands $\delta \lambda \geq 0$. On the other hand, by choosing $\lambda>0$ the test function $\delta \lambda$ can have any sign and we obtain from \req{rate_vp_threshold_6}$_3$ the equality $\widehat f = 0$, or in summary $\lambda \widehat f = 0$ as given in \req{rate_vp_threshold_7}$_6$. The central three field equations are identical to \req{extended_cont_8}--\req{extended_cont_10} and are complemented by the set of evolution equations \req{rate_vp_threshold_7}$_5$ which read
\begin{equation}
\dot \BC = \lambda \,  \partial_{\Bfrakm} \widehat f \, , \quad
\dot \Bfrakq = \lambda \, \partial_{\Bfrakd} \widehat f \, , \quad
\nabla \dot \Bfrakq = \lambda \, \partial_{\Bfrakg} \widehat f \, ,
\label{rate_vp_threshold_7e}
\end{equation}
 as already given in \citet{miehe14}, where the identity $T=\theta$ in $\calB$ from  \req{rate_vp_threshold_7}$_3$ has been used. Note that the only formal difference to \req{extended_cont_10e} is that the nonsmooth evolution is not written in form of  differential inclusions but is governed by the Karush-Kuhn-Tucker conditions \req{rate_vp_threshold_7}$_6$ containing the scalar variable $\lambda$. We also mention that the representation \req{rate_vp_threshold_7e} is possible because of the assumption of a smooth elastic domain $\mathbb{E}$, i.e., its boundary only consists of regular points.

In addition, we have for $\widehat f = 0$ the {\it nonlocal consistency conditions}
\begin{equation}
\lambda \geq 0 \, ,
\quad 
\frac{d}{dt}\widehat f(\Bfrakf; \Bfrakc, \theta) \leq 0 \, ,
\quad
\lambda \, \frac{d}{dt} \widehat f (\Bfrakf; \Bfrakc, \theta) = 0
\quad
\mbox{in } \calB \, ,
\label{cons_1}
\end{equation}
which at current time $t$ are supplemented by rate forms of  the equations \req{gov_equ_3} and \req{rate_vp_threshold_7}$_{1-2}$, the evolution equations \req{rate_vp_threshold_7}$_5$ and the rate forms of  the Neumann boundary conditions \req{extended_cont_7} related to the macro- and micro-motion, respectively. To see conditions \req{cons_1}, consider at current time $t$ a nonzero evolution of the mechanical constitutive state, i.e., $\lambda_t > 0$. Then, the first variation of the dissipation potential functional \req{threshold_6} with respect to the Lagrange mulitplier vanishes
\begin{equation}
\delta_\lambda V^*_\lambda \vert_t = \int_{\calB} - \delta \lambda \, \widehat f \vert_t \, dV = 0
\label{cons_2} 
\end{equation}
for all $\delta \lambda$. Next, at time $t+ \tau$, $\tau \geq 0$ we  consider a state with $\lambda_{t+\tau} \geq 0$ and the first variation of the dissipation potential functional \req{threshold_6} with respect to the Lagrange multiplier is nonnegative
\begin{equation}
\delta_\lambda V^*_\lambda \vert_{t + \tau} = \int_\calB -\delta \lambda \, \widehat f \vert_{t+ \tau} \, dV \geq 0
\label{cons_3}
\end{equation}
for all $\lambda_{t+\tau} + \delta \lambda \geq 0$. Subtracting \req{cons_2} from \req{cons_3} and dividing by $\tau$ yields
\begin{equation*}
\frac{1}{\tau} [ \,\delta_\lambda V^*_\lambda \vert_{t+\tau} - \delta_\lambda V_\lambda^* \vert_t  \,] = \int_\calB [\, - \delta \lambda \, ( \frac{d}{dt} \widehat f + \frac{\calO(\tau^2)}{\tau})  \, ] \, dV \geq 0
\end{equation*}
for all $\lambda_{t+\tau} + \delta \lambda \geq 0$. For $\tau \rightarrow 0$ we have $\lambda_{t+\tau} \rightarrow \lambda_t >0$, $\calO(\tau^2)/\tau \rightarrow 0$ and get $\frac{d}{dt}\widehat f = 0$ since $\delta \lambda$ can have any sign. For $\tau$ small enough, we assume $\lambda_{t+\tau} = 0$ and $\delta \lambda$ must be nonnegative yielding $\frac{d}{dt}\widehat f \leq 0$. When summarizing, we arrive at the nonlocal consistency conditions \req{cons_1}. From this condition the Lagrange multiplier {\it field}  can be determined in terms of the rates $(\dot \Bvarphi, \dot \theta)$ of the external fields deformation and temperature.

\subsubsection[Section452]{Incremental Formulation} \label{Section452}
Within a time interval $[t_n, t_{n+1}] \subset \overline \calT$ a variational principle can be constructed by the same avenue as outlined in Section~\ref{Section432}. 
Using the algorithmic approximations \req{inc_vp_4} of rates of state quantities, we get the incremental internal potential density
\begin{equation}
\fterm{
\widehat \pi^{* \tau}_\lambda =
\widehat \psi(\Bfrakc, \theta)  + \tau [ \, (\theta_k - T) \dot \eta^\tau +  \eta_n \dot \theta^\tau + \frac{T}{\theta_n} \, \Bfrakf \cdot \dot \Bfrakc^\tau - \lambda \widehat f(\Bfrakf; \Bfrakc_n, \theta_n) - \widehat \phi_{con}(-\frac{1}{T} \nabla T;\Bfrakc_n,\theta_n)  \,] \,  ,
}
\label{inc_vp_threshold_2}
\end{equation}
where like before $k=n+1$ corresponds to a fully implicit and $k=n$ to a semi-explicit update. Again, terms that are associated with previous time $t_n$ are dropped. Then, with the use of the incremental load functional \req{inc_vp_6} the incremental variational principle 
\begin{equation}
\{ \Bvarphi, \Bfrakq, \eta, \theta, T, \Bfrakf, \lambda  \} =
\mbox{Arg} \{\;
\substackrel{\Bvarphi, \Bfrakq, \eta }{\mbox{inf}}\,
\substackrel{\theta,T,\Bfrakf}{\mbox{sup}} \,
\substackrel{\lambda \geq 0}{\mbox{inf}} \;
[ \, 
\int_{\calB} \widehat \pi^{* \tau}_\lambda(\Bfrakc, \eta, \theta, T,\nabla T,\Bfrakf,\lambda) \, dV - P_{ext}^\tau(\Bvarphi,T)\, ] 
\; \}
\label{inc_vp_threshold_3}
\end{equation}
determines all fields at time $t_{n+1}$. The corresponding Euler equations in $\calB$ are time-discrete forms of \req{rate_vp_threshold_7} together with \req{inc_vp_8}$_4$ stemming from the variation with respect to the temperature $\theta$. Considering the semi-explicit integration, the incrementally isentropic operator split is modified by an additional optimization in the isentropic predictor step with respect to the Lagrange multiplier $\lambda \geq 0$. 
 \section[Section5]{Representative Model Problems}\label{Section5}

We now apply the general variational setting for thermomechanics of gradient-extended dissipative solids to complex multifield problems related to evolutions of {\it species content} in elastic solids, {\it damage} and {\it plasticity}. Especially, we highlight a new {\it minimization structure} for Cahn-Hilliard diffusion with respect to the species flux. The focus is put on incremental potentials which also inherently contain the proposed operator split. The latter is applied to a numerical example which shows the formation of an adiabatic shear band.

\subsection[Section51]{Cahn-Hilliard Diffusion coupled with Temperature Evolution} \label{Section51}
We consider as a first model problem the Cahn-Hilliard theory of diffusive phase separation in a rigid solid coupled with temperature evolution. In the following $c \colon \calB \times \calT \rightarrow [0,1]$ denotes a dimensionless concentration field  whose evolution is governed by the local form of conservation of species content
\begin{equation}
\dot c = - \Div \KIH \, ,
\label{ch_1}
\end{equation}
where $\KIH \colon \calB \times \calT \rightarrow \mathbb{R}^d$ is the species flux vector field. To give the concentration field the character of an order parameter $\Bfrakq^l = (c)$, we impose {\it homogeneous} boundary conditions $\dot c = 0$ on $\partial \calB_\Bfrakq$ and $\partial_{\nabla c} \widehat \psi \cdot \Bn_0 = 0$ on $\partial \calB_H$, see Figure~\ref{fields}. Hence, the dynamic process is only driven by an initial inhomogeneous distribution $c_0(\BX)$ of the concentration field in the domain $\calB$. In the following, we neglect the phenomenon of thermal diffusion (Soret effect) that is species flow caused by a temperature gradient, see \citet{groot+mazur13} and the recent contribution \citet{nateghi+keip21}. The free energy function decomposes into a local, nonlocal and purely thermal contribution
\begin{equation}
\widehat \psi(\Bfrakc, \theta) = \widehat \psi_l(c) + \widehat \psi_{\nabla}(\nabla c) + \widehat \psi_\theta(\theta) \, ,
\label{ch_3}
\end{equation}
where the last term is given in \req{0d_cons_1}. As considered in \citet{cahn+hilliard58} we choose
\begin{equation*}
\widehat \psi_l(c) = A[c \ln c + (1-c) \ln(1-c)] + Bc(1-c)
\AND
\widehat \psi_\nabla(\nabla c) = \frac{L}{2} | \nabla c |^2
\label{ch_4}
\end{equation*}
in terms of the threshold and mixing energy parameters $A$ and $B$ as well as the diffuse interface parameter $L$. Note the nonconvexity of $\widehat \psi_l$ for $B>2A$ which is related to phase  separation. The evolution of the concentration is driven by the chemical potential $\mu$ given by
\begin{equation}
\mu = \delta_c \widehat \psi = A \ln \frac{c}{1-c} + B(1-2c) - L \Delta c \, .
\label{ch_4a}
\end{equation}
It can be understood as a constitutive representation of a {\it micro-force balance} in the sense of \citet{gurtin94}.

\subsubsection[Section511]{Rate-type Minimization Principles in Isothermal Case} \label{Section511}
Point of departure is the definition of the energy and dissipation functionals 
\begin{equation}
E(c) = \int_\calB \widehat \psi(\Bfrakc) \, dV
\AND
D(\dot c;c) = \int_\calB \widehat \phi(\dot c; c) \, dV \, ,
\label{ch_5}
\end{equation}
where $\widehat \phi$ is the  smooth dissipation potential function accounting for diffusion mechanisms. A {\it minimization} of the corresponding rate-type potential with respect to $\dot c$ then yields as Euler equation the mass balance in the form 
\begin{equation*}
\delta_c \widehat \psi + \partial_{\dot c} \widehat \phi = 0 \quad \mbox{in } \calB \, .
\end{equation*}
Alternative to this setting, we now propose a new {\it minimization} formulation in terms of the species flux vector. In line with \citet{miehe+mauthe+teichtmeister15}, we reformulate the rate of the energy functional \req{ch_5}$_1$ at current concentration
\begin{equation}
\frac{d}{dt}E(c) = \calE(\KIH;c) = -\int_\calB \delta_c \widehat \psi \, \Div \KIH \, dV - \int_{\partial \calB} (\partial_{\nabla c} \widehat \psi \cdot \Bn_0) \, \Div \KIH \, dA \, ,
\label{ch_6}
\end{equation}
where we inserted the balance equation \req{ch_1}. The dissipation functional is defined as
\begin{equation}
X(\KIH; c) = \int_\calB \widehat \chi(\KIH;c) \, dV
\label{ch_7}
\end{equation}
in terms of the dissipation potential function $\widehat \chi$ which has the simple quadratic form
\begin{equation}
\widehat \chi(\KIH;c) = \frac{1}{2} \frac{1}{Mc(1-c)} \KIH \cdot \KIH \, .
\label{ch_8}
\end{equation}
Here, $M>0$ is a so-called  mobility parameter. Note that $\widehat \chi(\cdot \, ; c)$ is a positively homogeneous function of degree two and its image coincides with {\it half} the dissipation in a material element, see below. With the functionals \req{ch_6} and \req{ch_7} at hand, we define the potential
\begin{equation*}
\Pi(\KIH;c) = \calE(\KIH;c) + X(\KIH;c)
\label{ch_9}
\end{equation*}
with given concentration $c$. Its {\it minimization} with regard to {\it homogeneous} Dirichlet-type  boundary conditions $\KIH \cdot \Bn_0 = 0$ and $-\Div \KIH = 0$ on $\partial \calB_\Bfrakq$ determines  the current species flux field. We obtain the Euler equations
\begin{equation}
\parbox{14cm}{\centering 
$
\begin{array}{rlrll}
1.&{\mbox{Species flux }} &
 \nabla \delta_c \widehat \psi + \partial_\KIH \widehat \chi &= \Bzero
&\mbox{ in }\calB \,  , \\[1ex]
2. &{\mbox{Vanishing chemical potential}} & 
-\delta_c \widehat \psi  &= 0
&\mbox{ on } \partial \calB_H  \,  ,  \\[1ex]
3. &{\mbox{Vanishing micro-force}} &
\partial_{\nabla c} \widehat \psi \cdot \Bn_0 &= 0
&\mbox{ on }\partial\calB_H  \, .
\end{array}
$
}
\label{ch_10}
\end{equation}
They contain necessary {\it compatibility conditions} for the given concentration field. As a post-processing step, the current rate $\dot c$ of the concentration is determined via \req{ch_1}. Starting from \req{free_energy_2}, we now calculate the dissipation  whose unique source is diffusion
\begin{equation*}
\int_{\calB} \calD \, dV = \int_{\calB} \KIB \cdot \KIH \, dV 
\WITH
\KIB = - \nabla \delta_c \widehat \psi \, ,
\end{equation*}
where we performed integration by parts two times and inserted the balance \req{ch_1} as well as the homogeneous boundary conditions. Using \req{ch_10}$_1$, we can express the dual to the species flux vector by the dissipation potential $\widehat \chi$ and obtain with \req{ch_8} for the dissipation associated with a volume element $\calD = 2 \widehat \chi(\KIH;c) \geq 0$.

\subsubsection[Section512]{Incremental Formulation with Temperature Evolution} \label{Section512}
In addition to  \req{ch_8}, we consider the conductive dissipation potential function
\begin{equation*}
\widehat \phi_{con} (\KIg; \theta) = k \, \theta (\KIg \cdot \KIg)/2 \, ,
\label{ch_11}
\end{equation*}
which is convex in $\KIg = - \nabla \theta /\theta$ with $k>0$ being the thermal conductivity. For the incremental setting, we consider the implicit update of the species concentration
\begin{equation}
c = c_n - \tau \, \Div \KIH \, ,
\label{ch_12}
\end{equation}
and specify the incremental internal potential density \req{inc_vp_5}  as
\begin{equation*}
\begin{split}
\widehat \pi^{+ \tau} &=  \widehat \psi(c_n-\tau \, \Div \KIH, \nabla c_n - \tau \, \nabla [ \, \Div \KIH \, ], \theta) \\
 &+ \tau \, [ \, (\theta_k - T) \dot \eta^\tau + \eta_n \dot \theta^\tau+ \widehat \chi(\frac{T}{\theta_n} \KIH; c_n ,\theta_n) - \widehat \phi_{con}(- \frac{1}{T} \nabla T; \theta_n) \, ] \, .
\end{split}
\label{ch_13}
\end{equation*}
Here, $\widehat \chi$ is the dissipation potential function \req{ch_8} related to diffusion mechanisms, where an additional temperature dependence $M=\widehat M(\theta)$ of the mobility parameter is taken into account. We obtain the Euler equations
\begin{equation*}
\parbox{14cm}{\centering 
$
\begin{array}{lll}
\delta_{\KIH} \widehat \pi^{+\tau} \equiv \nabla \mu + [ \, \widehat M(\theta_n)c_n(1-c_n) \, ]^{-1} (\frac{T}{\theta_n})^2 \KIH &= \Bzero \,  , \\[1ex]
\partial_{\eta}\widehat \pi^{+\tau} \equiv \theta_k -T &= 0 \,  , \\[1ex]
\partial_{\theta} \widehat \pi^{+\tau} \equiv -C \ln \frac{\theta}{\theta_0} + \eta_k &= 0 \,  , \\[1ex]
\delta_{T} \widehat \pi^{+ \tau} \equiv - (\eta-\eta_n) + \frac{\tau}{\theta_n} [ \, \frac{T}{\theta_n} \frac{1}{\widehat M(\theta_n)c_n(1-c_n)} |\KIH|^2 \, ] +\frac{\tau}{T} \, \Div[ \, k \, \theta_n \frac{1}{T} \nabla T \,] &=- \tau \bar r / \theta_n 
\end{array}
$
}
\label{ch_14}
\end{equation*}
in $\calB$, where we recall the definition of the chemical potential \req{ch_4a} together with \req{ch_12}.

\subsection[Section52]{Thermomechanics of Gradient Damage} \label{Section52}
We consider as a second application the thermomechanics of a gradient damage model with an {\it elastic stage}. The scalar micro-motion field $d \colon \calB \times \calT \rightarrow [0,1]$, referred to as damage variable, measures at a macroscopic point $\BX \in \calB$ the ratio between an arbitrary oriented area of microcracks and a representative reference surface in which the mentioned crack surfaces are embedded, see, e.g., \citet{lemaitre96}. In this sense, a value $d=0$ characterizes an unbroken state, whereas $d=1$ represents a fully broken state. The irreversibility of micro-cracking is usually expressed by the inequality constraint $\dot d(\BX,t) \geq 0$ on the evolution of the damage variable. The mechanical constitutive state is  specified as
\begin{equation*}
\Bfrakc = (\BC, d, \nabla d) 
\label{gd_1}
\end{equation*}
and contains the right Cauchy-Green tensor $\BC$, the damage variable as well as its first gradient. In addition, we introduce the elastic right Cauchy-Green tensor $\BC^e = \BF^{eT}\Bg \BF^e$ that is based on the definition of an elastic, stress producing part $\BF^e=J_\theta ^{-1/3} \BF$ of the deformation gradient in terms of a volumetric thermal expansion $J_{\theta} = \exp[3 \alpha_T (\theta-\theta_0)]$, see \citet{lu+pister75}. One can then write $\BC^e$ by means of $\BC$ as $\BC^e = J_\theta^{-2/3} \BC$. A simple model of thermo-gradient-damage at large deformations may then be based on the objective free energy function
\begin{equation}
\widehat \psi(\Bfrakc,\theta)=\widehat g (d) \widehat \psi_e(\BC^e) + \widehat \psi_{\theta}(\theta)
\WITH
\widehat \psi_e = \frac{\mu}{2}(\tr \BC^e-3) + \frac{\mu}{\delta}[\,(\det \BC^e)^{-\frac{\delta}{2}} - 1 \,] \, ,
\label{gd_2}
\end{equation}
where $\widehat g(d)=(1-d)^2$ is a degradation function and $\widehat \psi_{\theta}$ the purely thermal contribution as given in \req{0d_cons_1}. Note that the gradient of damage does not arise in this constitutive function but will exclusively enter the dissipation potential function, see below. $\mu>0$ is the shear modulus and the second parameter $\delta >0$ models a weak compressibility. Note that only the full elastic energy storage is degraded. From \req{gd_2} we obtain the tensor functions
\begin{equation}
\begin{array}{llllll}
& \Bg \BP^e&=&  \partial_{\BF} \widehat \psi &=&  \widehat g(d) J_\theta ^{-\frac{2}{3}}  \Bg \BF [\,\mu \BG^{-1} - \mu (\det \BC^e)^{- \frac{\delta}{2}} \BC^{e-1} \, ] \,  ,  \\[1ex]
& \beta^e&=&  \partial_d \widehat \psi &=&  - 2(1-d) [ \, \frac{\mu}{2}(\tr \BC^e-3) + \frac{\mu}{\delta}((\det \BC^e)^{-\frac{\delta}{2}} - 1) \, ] \,  , \\[1ex]
& \tilde \eta &=& \partial_{\theta} \widehat \psi &=& - \widehat g(d) \alpha_T [ \, \mu \tr \BC^e - \mu (\det \BC^e)^{-\frac{\delta}{2}} \BC^{e}:\BC^{e-1}  \,] - C \ln \frac{\theta}{\theta_0} \, ,
\end{array}
\label{gd_3}
\end{equation}
that represent constitutive relationships for driving forces. 

\subsubsection[Section521]{Rate-type Formulation based on Indicator Function} \label{Section521}
We consider the intrinsic dissipation potential function
\begin{equation}
\widehat \phi_{int}(\dot d, \nabla \dot d; \nabla d, \theta) = \widehat \phi_{l}(\dot d; \theta)+\widehat \phi_{\nabla} (\nabla \dot d; \nabla d)
 \label{gd_4}
\end{equation}
as the sum of a nonsmooth local and a smooth nonlocal part
\begin{equation}
\widehat \phi_{l}(\dot d; \theta) = \widehat c(\theta) \dot d +  I_+(\dot d)
\AND
\widehat \phi_{\nabla} (\nabla \dot d; \nabla d) = \mu l^2 \nabla d \cdot \nabla \dot d \, .
\label{gd_4a}
\end{equation}
Here, irreversibility of damage is ensured by the indicator function $I_+(\dot d)$ of the set of positive real numbers defined as
\begin{equation}
I_+(\dot d) =
\left\{
\begin{array}{ll}
 0 \quad & \mbox{for }\dot d \geq 0 \,  ,\\
 + \infty \quad & \mbox{otherwise}
\end{array}
\right.
\AND
\partial I_+ (\dot d) =
\left\{
\begin{array}{ll}
 0 \quad & \mbox{for }\dot d > 0 \, ,\\
 \mathbb{R}_{-} \quad & \mbox{for } \dot d = 0 \, ,\\
 \emptyset \quad & \mbox{otherwise}
\end{array}
\right.
\label{gd_5}
\end{equation}
with $\partial$ denoting the subdifferential. The parameter $\widehat c(\theta)>0$ is a temperature-dependent force-like threshold value for the onset of damage with $\frac{d}{d \theta} \widehat c<0$ and $l$ a length scale parameter. Note that \req{gd_4} is a positively homogeneous function of degree one in $(\dot d, \nabla \dot d)$ and hence models for an adiabatic process a rate-independent evolution of damage. In addition, we have the
conductive dissipation potential function
\begin{equation}
\widehat \phi_{con}(\KIg; \BC,d, \theta) = \widehat k(d) \theta \, \BC^{-1}: (\KIg \otimes \KIg)/2 \, ,
\label{gd_6}
\end{equation}
where the thermal conductivity is a function of the damage variable and may take the simple form $\widehat k(d) = \widehat g(d)k_b$ in terms of the heat conduction coefficient $k_b > 0$ of the undamaged bulk. With the constitutive functions \req{gd_2}, \req{gd_4} and \req{gd_6} at hand, we specify the internal potential density \req{cont_vp_4} as
\begin{equation*}
\widehat \pi^+_{red} = \partial_{\BC} \widehat \psi : \dot \BC+\beta^e \dot d +(\theta-T) \dot \eta + \frac{T}{\theta}[\, \widehat c(\theta)\dot d + \mu l^2 \nabla d \cdot \nabla \dot d + I_+(\dot d) \, ]  - \widehat \phi_{con}(-\frac{1}{T} \nabla T) \, .
\label{gd_7}
\end{equation*}
Then, the variational principle \req{cont_vp_5} determines at current time $t$ the rates of deformation, damage and entropy as well as the thermal driving force and gives the Euler equations
\begin{equation}
\parbox{14cm}{\centering 
$
\begin{array}{lll}
\delta_{\dot \Bvarphi} \widehat \pi^+_{red} \equiv   - \Div \, \Bg \BP^e &= \Bg \bar \Bgamma \,  , \\[1ex]
 \delta_{\dot d} \widehat \pi^+_{red} \equiv \beta^e + \frac{T}{\theta}\widehat c(\theta)- \mu l^2 \Delta(\frac{T}{\theta}d) + \partial I_+(\dot d) &\ni 0 \,  , \\[1ex]
\partial_{\dot \eta}\widehat \pi^+_{red} \equiv \theta -T &= 0 \,  , \\[1ex]
\delta_{T} \widehat \pi^+_{red} \equiv -\dot \eta + \frac{1}{\theta} \widehat \phi_{int}(\dot d, \nabla \dot d) +\frac{1}{T} \, \Div[ \, \widehat k(d) \theta \, \BC^{-1} \frac{1}{T} \nabla T \,] &=-\bar r /\theta

\end{array}
$
}
\label{gd_8}
\end{equation}
in $\calB$. Observe, that the evolution of the entropy is driven by the rates of damage and gradient of damage. The nonsmooth evolution of the damage variable is governed by the differential inclusion \req{gd_8}$_2$. From there and the relation \req{gd_8}$_3$, we conclude for the determination of the rates of deformation and damage the nonlocal consistency conditions\footnote{These also follow from \req{cons_1} by using the threshold function defined in \req{gd_11} and the identifications
\begin{equation*}
\lambda = \dot d 
\AND 
\beta = - \beta^e + \Div (\partial_{\nabla \dot d} \widehat \phi_{\nabla})
\end{equation*}
for the Lagrange multiplier and the nonlocal driving force considered below.
 }
\begin{equation*}
\dot d \geq 0 \, ,
\quad
-\dot \beta^e - \tfrac{d}{d\theta}\widehat c(\theta) \, \dot \theta + \mu l^2 \Delta \dot d \leq 0 \, ,
\quad
\dot d \, [\, -\dot \beta^e - \tfrac{d}{d\theta}\widehat c(\theta) \, \dot \theta + \mu l^2 \Delta \dot d  \, ] = 0
\quad
\mbox{in }
\calB \, ,
\label{gd_8a}
\end{equation*}
where $\dot \beta^e = \partial^2_{d \BC} \widehat \psi : \dot \BC + \partial_{dd}^2 \widehat \psi \, \dot d + \partial_{d \theta}^2 \widehat \psi \, \dot \theta$ and $\dot \theta$ follows from taking the time derivative of the state equation \req{gov_equ_3}. In addition, the rate form
\begin{equation*}
\Div [ \, \partial_\BF \BP^e : \dot \BF + \partial_d \BP^e \, \dot d + \partial_\theta \BP^e \, \dot \theta  \,] = \dot{\bar \Bgamma} 
\label{gd_8b}
\end{equation*}
of mechanical equilibrium \req{gd_8}$_1$ has to be considered together with the  boundary conditions
\begin{equation}
[\,  \partial_\BF \BP^e : \dot \BF + \partial_d \BP^e \, \dot d + \partial_\theta \BP^e \, \dot \theta  \,] \, \Bn_0 = \dot{\bar \Bt} 
\AND
\nabla \dot d \cdot \Bn_0 = 0
\label{gd_8c}
\end{equation}
on  $\partial \calB_\Bt$ and $\partial \calB_H$, respectively. \req{gd_8c} are rate forms of the Neumann boundary conditions $\BP^e \, \Bn_0 = \bar \Bt$ on $\partial \calB_\Bt$ and $\nabla d \cdot \Bn_0 = 0$ on $\partial \calB_H$  which are outcomes of the variational principle.

From \req{dissipation_potential_5} the intrinsic dissipation is found to be
\begin{equation*}
\int_\calB \calD_{int} \, dV = -\int_\calB \beta^e \dot d \, dV \geq 0 \,  ,
\label{gd_9}
\end{equation*}
where we used integration by parts and the homogeneous boundary conditions. Hence, thermodynamic consistency is shown.

\subsubsection[Section522]{Incremental Formulation based on Threshold Function} \label{Section522}
We specify the array \req{mixed_mech_2} of dissipative driving forces  as
\begin{equation*}
\Bfrakf = (\Bzero, \beta, \Bzero) \, ,
\label{gd_10}
\end{equation*}
where $\beta$ is the quantity conjugate to $d$. The Legendre transform of the local part  \req{gd_4a}$_1$ of the intrinsic dissipation potential function reads
\begin{equation}
\widehat \phi_{l}^*(\beta;\theta) = \sup_{\dot d} [\, (\beta-\widehat c) \dot d - I_+(\dot d)  \,  ] = \sup_{\dot d \geq 0} [ \, (\beta - \widehat c) \dot d \,] = 
\left\{
\begin{array}{ll}
 0 &\mbox{ for } \beta - \widehat c \leq 0 \, ,\\
 + \infty &\mbox{ for } \beta - \widehat c > 0
\end{array}
\right.
\label{gd_10a}
\end{equation}
and enters the incremental internal potential density \req{extended_inc_2}. Note  that $\widehat \phi_l^*$ is the indicator function of the set \req{threshold_2} of admissible driving forces governed by the threshold function\footnote{An alternative intrinsic local dissipation potential function may be given by
\begin{equation*}
\widehat \phi_l(\dot d; d, \theta) = \widehat c(\theta) d \dot d + I_+(\dot d) \, ,
\end{equation*}
that in contrast to \req{gd_4a}$_1$ depends explicitly on the given damage state $d$. We obtain the threshold function
\begin{equation*}
\widehat f(\beta; d, \theta) = \beta - \widehat c(\theta) d \, ,
\end{equation*}
which, due to occurrence of the damage variable in the resistance term, corresponds to a model {\it without an elastic stage}, i.e., the damage starts to evolve from $d=0$ at the instant of loading. 
  }
\begin{equation}
\widehat f(\beta;\theta) = \beta - \widehat c(\theta) \, .
\label{gd_11}
\end{equation}
The latter defines the local intrinsic dissipation potential function by the constrained optimization problem
\begin{equation*}
\widehat \phi_{l}(\dot d; \theta) = \sup_\beta \inf_{\lambda \geq 0} [ \, \beta \dot d - \lambda \widehat f(\beta;\theta)  \, ] \, .
\label{gd_12}
\end{equation*}
Then, with the exact time integration of the nonlocal term \req{gd_4a}$_2$ of the intrinsic dissipation potential function
\begin{equation}
\widehat \phi_{\nabla}^\tau(\nabla d; \nabla d_n)=\int_{t_n}^{t_{n+1}}  \mu l^2 \nabla d \cdot \nabla \dot  d \, dt = \frac{1}{2}\mu l^2 (|\nabla d|^2 - |\nabla d_n|^2)  \, ,
\label{gd_13}
\end{equation}
the incremental internal potential density  \req{inc_vp_threshold_2} takes the form
\begin{equation*}
\widehat \pi^{*\tau}_\lambda = \widehat \psi(\Bfrakc,\theta) + \tau [\,(\theta_k - T) \dot \eta^\tau +  \eta_n \dot \theta^\tau + \frac{T}{\theta_n} \, ( \beta \dot d^\tau + \frac{1}{\tau} \widehat \phi_{\nabla}^\tau ) - \lambda \widehat f(\beta;\theta_n)  - \widehat \phi_{con}(\KIg;\BC_n,d_n,\theta_n) \,] \, .
\label{gd_14}
\end{equation*}
As a result, the incremental variational principle \req{inc_vp_threshold_3} gives the Euler equations 
\begin{equation}
\parbox{14cm}{\centering 
$
\begin{array}{lll}
\delta_{\Bvarphi} \widehat \pi^{+\tau}_{\lambda} \equiv - \Div \, \Bg \BP^e &= \Bg \Bar \Bgamma  \,  ,  \\[1ex]
 \delta_{d} \widehat \pi^{+\tau}_{\lambda} \equiv \beta^e + \frac{T}{\theta_n}\beta - \mu l^2 \Delta(\frac{T}{\theta_n}d) &= 0  \,  ,  \\[1ex]
\partial_{\eta}\widehat \pi^{+\tau}_{\lambda} \equiv \theta_k -T &= 0   \,  ,    \\[1ex]
\partial_{\theta} \widehat \pi^{+\tau}_{\lambda} \equiv \tilde \eta + \eta_k &= 0 \,  ,   \\[1ex]
\delta_{T} \widehat \pi^{+ \tau}_{\lambda} \equiv - (\eta-\eta_n) + \frac{1}{\theta_n} [ \, \beta (d-d_n) + \widehat \phi_{\nabla}^\tau \, ] +\frac{\tau}{T} \, \Div[ \, \theta_n \widehat k(d_n) \, \BC^{-1}_n \frac{1}{T} \nabla T \,] &=- \tau \bar r / \theta_n  \,  ,  \\[1ex]
\partial_\beta \widehat \pi^{+ \tau}_{\lambda} \equiv \frac{T}{\theta_n} (d-d_n) -\tau  \lambda &= 0  \,  ,   \\[1ex]
\partial_\lambda \widehat \pi^{+ \tau}_{\lambda} \equiv -\tau [\, \beta -\widehat c(\theta_n) \, ] \geq 0 \, , \,  \lambda \geq 0 \, , \, \tau \lambda [\, \beta - \widehat c(\theta_n) \,] = 0

\end{array}
$
}
\label{gd_15}
\end{equation}
in $\calB$, which represent time-discrete forms of the general equations \req{gov_equ_3} and \req{rate_vp_threshold_7}. We can reduce this set of equations by expressing for $k=n$ the dissipative driving force $\beta = - \beta^e + \mu l^2 \Delta d$ via \req{gd_15}$_2$ and the Lagrange multiplier $\tau \lambda = d - d_n$ via \req{gd_15}$_6$ yielding the explicit nonlocal form of the Karush-Kuhn-Tucker conditions
\begin{equation*}
d  \geq d_n \, ,
\quad 
\mu l^2 \Delta d - \beta^e - \widehat c(\theta_n) \leq 0 \, ,
\quad
(d-d_n) [\, \mu l^2 \Delta d - \beta^e - \widehat c(\theta_n) \,] = 0 \, ,
\label{gd_16}
\end{equation*}
where we recall the definition of the driving force \req{gd_3}$_2$.

As an alternative to this setting, which is fully rate-independent in the adiabatic case, we may consider a (regularized) viscous over-force formulation based on the {\it smooth} dual local intrinsic dissipation potential function, see Footnote~\ref{footnote_perzyna},
\begin{equation*}
\widehat \phi^{\eta *}_{l}(\beta; \theta) = \frac{1}{2 \eta_f} \langle \, \widehat f(\beta;\theta) \, \rangle^2_+  \, ,
\label{gd_16a}
\end{equation*}
that approaches \req{gd_10a} for the vanishing viscosity limit $\eta_f \rightarrow 0$. Then, with \req{gd_13} the incremental internal potential density \req{inc_vp_threshold_2} takes the form
\begin{equation*}
\begin{split}
\widehat \pi^{*\tau} = \widehat \psi + \tau [\,(\theta_k - T) \dot \eta^\tau +  \eta_n \dot \theta^\tau  + & \frac{T}{\theta_n} \, ( \beta \dot d^\tau + \frac{1}{\tau} \widehat \phi_{\nabla}^\tau ) \\ &- \frac{1}{2 \eta_f} \langle \widehat f(\beta;\theta_n) \rangle^2_+   - \widehat \phi_{con}(-\frac{1}{T}\nabla T;\BC_n,d_n, \theta_n) \,] \, ,
\end{split}
\end{equation*}
and the incremental variational principle \req{extended_inc_3} yields \req{gd_15}$_1$--\req{gd_15}$_5$ together with
\begin{equation*}
\partial_\beta \widehat \pi^{* \tau} \equiv \frac{T}{\theta_n}(d-d_n) - \frac{\tau}{\eta_f} \langle \beta - \widehat c(\theta_n) \rangle_+ = 0
\label{gd_19}
\end{equation*}
as Euler equations in $\calB$. 
We can reduce this set of equations by expressing the dissipative driving force $\beta$ from \req{gd_15}$_2$ yielding for $k=n$ the nonlocal (regularized) update equation
\begin{equation*}
d = d_n + \frac{\tau}{\eta_f} \langle  \mu l^2 \Delta d - \beta^e - \widehat c(\theta_n) \rangle_+
\end{equation*}
for the damage variable.

\subsection[Section53]{Thermomechanics of Additive Gradient Plasticity} \label{Section53}
As third model problem, we consider a thermomechanically coupled formulation of additive gradient plasticity. Besides the standard metrics $\BG$ and $\Bg$, we introduce on the reference configuration the (covariant) plastic metric tensor $\BG^p \in Sym_+(3)$ that evolves in time starting from the initial state $\BG^p(\BX,t_0) = \BG$. Following \citet{miehe+apel+lambrecht02} and as visualized in Figure~\ref{mappings}, a Lagrangian {\it elastic strain variable} may be based on an explicit dependence on the right Cauchy-Green tensor $\BC$, that is the current metric pulled back to the reference configuration, and the plastic metric $\BG^p$ in an {\it additive} format
\begin{equation*}
\Bvarepsilon^e = \Bvarepsilon - \Bvarepsilon^p \,  ,
\label{gp_1}
\end{equation*}
where the total and plastic Hencky strain tensors
\begin{equation*}
\Bvarepsilon = \frac{1}{2} \ln \BC
\AND
\Bvarepsilon^p = \frac{1}{2} \ln \BG^p
\label{gp_2}
\end{equation*}
are introduced. Hence, instead of $\BG^p$ we consider in what follows the logarithmic plastic strain $\Bvarepsilon^p$ as the local internal variable  whose evolution from $\Bvarepsilon^p(\BX,t_0) = \Bzero$ is governed by a standard flow rule. Note that within this framework the condition of plastic incompressibility $\det \BG^p = 1$ is equivalent to a standard statement of vanishing trace $\tr \Bvarepsilon^p = 0$. 
%
\inputfig[t]{floats/mappings}{mappings}%
%
We specify the mechanical constitutive state 
\begin{equation*}
\Bfrakc = (\Bvarepsilon, \Bvarepsilon^p, \alpha, \nabla \alpha) \, ,
\label{gp_3}
\end{equation*}
which contains a scalar hardening variable $\alpha$ as well as its first gradient. In the following, we focus on metal plasticity characterized by small elastic but large plastic deformations and consider the free energy function
\begin{equation}
\widehat \psi(\Bfrakc, \theta) = \widehat \psi_e(\Bvarepsilon, \Bvarepsilon^p) + \widehat \psi_p(\alpha, \theta) + \frac{1}{2} \mu l^2 | \nabla \alpha |^2 + \widehat \psi_{e-\theta}(\Bvarepsilon, \theta) + \widehat \psi_\theta(\theta) \, .
\label{gp_4}
\end{equation}
Here, $\widehat \psi_e$ is the purely elastic contribution that is assumed to have the quadratic form
\begin{equation}
\widehat \psi_e(\Bvarepsilon, \Bvarepsilon^p) = \frac{\kappa}{2} (\tr \Bvarepsilon)^2 + \mu | \, \mbox{Dev} \, \Bvarepsilon^e \,|^2 \, ,
\label{gp_5}
\end{equation}
where $\kappa>0$ and $\mu>0$ are the bulk and shear modulus, respectively. $\widehat \psi_p$ is an isotropic hardening function that also takes into account thermally induced softening. The gradient extension related to a length scale parameter $l$ is assumed to affect the scalar hardening variable only. The coupled thermoelastic response is modelled by the constitutive function
\begin{equation*}
\widehat \psi_{e-\theta}(\Bvarepsilon, \theta) = - \kappa \alpha_T (\tr \Bvarepsilon) (\theta-\theta_0)
\label{gp_6}
\end{equation*}
in line with \req{0d_cons_1}, where also the pure thermal contribution $\widehat \psi_\theta$ is specified. Note, that the function \req{gp_5}, known as Hencky energy, is {\it not polyconvex}\footnote{\citet{neff+ghiba+lankeit15} introduced an exponentiated Hencky energy which is polyconvex in the two-dimensional case.} with respect to the deformation gradient $\BF = D \Bvarphi$ in the sense of \citet{ball76}, but {\it rank-one convex} for a moderately high elastic deformation range, see \citet{bruhns+etal01}. Hence, it is applicable to the typical scenario of metal plasticity. With the free energy function \req{gp_4} at hand, we can derive the tensor functions
\begin{equation}
\begin{array}{llllll}
& \Bsigma^e&=&  \partial_{\Bvarepsilon} \widehat \psi &=&  \kappa [\, \tr \Bvarepsilon - \alpha_T(\theta-\theta_0) \,] \, \BG^{-1} + 2\mu \, \mbox{Dev} [\BG^{-1} \Bvarepsilon^e \BG^{-1}] \,  ,  \\[1ex]
& \Bbeta^e&=&  \partial_{\Bvarepsilon^p} \widehat \psi &=& - 2\mu \, \mbox{Dev} [\BG^{-1} \Bvarepsilon^e \BG^{-1}]  \,  ,  \\[1ex]
& \beta^e&=&  \delta_\alpha \widehat \psi &=&  \partial_{\alpha} \widehat \psi_p - \mu l^2 \Delta \alpha  \,  ,   \\[1ex]
& \tilde \eta &=& \partial_{\theta} \widehat \psi &=& \partial_\theta \widehat \psi_p - \kappa \alpha_T \tr \Bvarepsilon - C \ln \frac{\theta}{\theta_0} \, ,
\end{array}
\label{gp_7}
\end{equation}
that represent constitutive relationships for driving forces.  Note the occurrence of the Laplacian term $\mu l^2 \Delta \alpha$ in \req{gp_7}$_3$ that is in line with the approach of Aifantis, see, e.g., \citet{aifantis87}.

\subsubsection[Section531]{Incremental Formulation based on Threshold Function} \label{Section531}
We specify the array \req{mixed_mech_2} of dissipative driving forces  as
\begin{equation*}
\Bfrakf = (\Bzero, \Bs, \beta, \Bzero) \, ,
\label{gp_8}
\end{equation*}
where $(\Bs, \beta)$ are the quantities conjugate to $(\Bvarepsilon^p, \alpha)$. For von Mises-type gradient plasticity we define the yield function
\begin{equation}
\widehat f (\Bs, \beta;\theta) = | \, \Bs \,| - \sqrt{\frac{2}{3}} \, [ \, \widehat y(\theta) - \beta \, ] \, ,
\label{gp_9}
\end{equation}
that restricts the set of admissible driving forces according to \req{threshold_2}. $\widehat y(\theta)$ is a temperature dependent yield stress function with $\frac{d}{d\theta} \widehat y < 0$. The corresponding intrinsic dissipation potential function is defined by the constrained optimization problem
\begin{equation}
\widehat \phi_{int}(\dot \Bvarepsilon ^p, \dot \alpha; \theta) = \sup_{\Bs, \beta} \inf_{\lambda \geq 0} [\, \Bs : \dot \Bvarepsilon^p + \beta \dot \alpha - \lambda \widehat f(\Bs, \beta; \theta)  \,] \, .
\label{gp_10}
\end{equation}
The intrinsic dissipation follows as
\begin{equation}
\calD_{int} = - ( \Bbeta^e : \dot \Bvarepsilon^p + \beta^e \dot \alpha  )   =  \lambda (\Bs:\partial_\Bs \widehat f + \beta \partial_\beta \widehat f) =  \sqrt{\frac{2}{3}} \, \lambda \widehat y(\theta)  \geq 0 \, ,
\label{gp_10_a}
\end{equation}
where we used the evolution laws for the plastic strain and hardening variable stemming from \req{gp_10} and the identities $\Bs = -\Bbeta^e$ and $\beta = -\beta^e$ which are outcomes of the rate-type variational principle \req{rate_vp_threshold_5}. Note that the last equality in \req{gp_10_a} follows from $\widehat w (\Bs, \beta) = |\, \Bs \, |+ \sqrt{2/3}\, \beta$ being a gauge, i.e., $\widehat w$ is a positively homogeneous function of degree one with the resulting property $\Bs : \partial_\Bs \widehat w + \beta \, \partial_\beta \widehat w = \widehat w$. With \req{gp_10} the incremental internal potential density  \req{inc_vp_threshold_2} is  specified as
\begin{equation*}
\begin{split}
\widehat \pi^{* \tau}_{\lambda} = \widehat \psi(\Bfrakc,\theta) + \tau [\,(\theta_k - T) \dot \eta^\tau +  \eta_n \dot \theta^\tau + \frac{T}{\theta_n} \, & ( \Bs : \dot \Bvarepsilon^{p \tau} +\beta \dot \alpha^{\tau} ) \\ & - \lambda \widehat f(\Bs,\beta;\theta_n)  - \widehat \phi_{con}(-\frac{1}{T} \nabla T;\BC_n, \theta_n) \,] \, ,
\end{split}
\label{gp_11}
\end{equation*}
where  $\widehat \phi_{con}$ is the conductive dissipation potential function defined in \req{gd_6}, however, with a constant heat conduction coefficient. Then, the incremental variational principle \req{inc_vp_threshold_3} gives the Euler equations
\begin{equation}
\parbox{14cm}{\centering 
$
\begin{array}{lll}
\delta_{\Bvarphi} \widehat \pi^{+\tau}_{\lambda} \equiv - \Div [\, \Bsigma^e:\partial_\BF \Bvarepsilon \, ] &= \Bg \Bar \Bgamma  \, ,  \\[1ex]
\partial_{\Bvarepsilon^p} \widehat \pi^{+\tau}_{\lambda} \equiv  \Bbeta^e + \frac{T}{\theta_n} \Bs&= \Bzero \,  ,   \\[1ex]
 \delta_{\alpha} \widehat \pi^{+\tau}_{\lambda} \equiv \beta^e + \frac{T}{\theta_n}\beta  &= 0 \,  ,   \\[1ex]
\partial_{\eta}\widehat \pi^{+\tau}_{\lambda} \equiv \theta_k -T &= 0 \,  ,   \\[1ex]
\partial_{\theta} \widehat \pi^{+\tau}_{\lambda} \equiv \tilde \eta + \eta_k &= 0 \,  , \\[1ex]
\delta_{T} \widehat \pi^{+ \tau}_{\lambda} \equiv - (\eta-\eta_n) + \frac{\tau}{\theta_n} [ \, \Bs:\dot \Bvarepsilon^{p \tau} + \beta \dot \alpha^\tau  \, ] +\frac{\tau}{T} \, \Div[ \,\theta_n k \BC^{-1}_n \frac{1}{T} \nabla T \,] &=-\tau \bar r / \theta_n  \,  ,   \\[1ex]
\partial_\Bs \widehat \pi^{+ \tau}_{\lambda} \equiv \frac{T}{\theta_n} (\Bvarepsilon^p-\Bvarepsilon^p_n) -\tau  \lambda \, \BG \Bs \BG / | \,  \Bs \, |&= \Bzero \, , \\[1ex]
\partial_\beta \widehat \pi^{+ \tau}_{\lambda} \equiv \frac{T}{\theta_n} (\alpha-\alpha_n) -\tau  \lambda \sqrt{2/3}&= 0 \,  , \\[1ex]
\partial_\lambda \widehat \pi^{+ \tau}_{\lambda} \equiv -\tau \widehat f(\Bs, \beta; \theta_n) \geq 0 \, , \, \lambda \geq 0 \, , \, \tau \lambda \widehat f (\Bs, \beta; \theta_n) = 0

\end{array}
$
}
\label{gp_12}
\end{equation}
in $\calB$ which represent time-discrete forms of the general equations  \req{gov_equ_3} and \req{rate_vp_threshold_7}. 
We can reduce the set of equations \req{gp_12} by expressing for $k=n$ the dissipative driving forces $\Bs = - \Bbeta^e$ and $\beta = - \beta^e$ via \req{gp_12}$_2$ and \req{gp_12}$_3$, respectively, and the Lagrange multiplier $\tau \lambda = \sqrt{3/2} \, (\alpha - \alpha_n)$ via \req{gp_12}$_8$ yielding the nonlocal form of the Karush-Kuhn-Tucker conditions in strain space
\begin{equation*}
\alpha \geq \alpha_n \, , 
\quad
| \, \Bbeta^e  \, | - \sqrt{2/3} \, (\widehat y(\theta_n) + \beta^e) \leq 0 \, ,
\quad
(\alpha - \alpha_n) [ \, | \, \Bbeta^e  \, | - \sqrt{2/3} \, (\widehat y(\theta_n) + \beta^e) \, ] = 0 \, ,
\label{gp_13}
\end{equation*}
where we recall the definitions of the driving forces \req{gp_7}$_{2-3}$. 

As an alternative to this setting, which is fully rate-independent in the adiabatic case, we may consider a (regularized) viscous over-force formulation based on the {\it smooth} dual intrinsic dissipation potential function
\begin{equation*}
\widehat \phi^{\eta *}_{int}(\Bs, \beta; \theta) = \frac{1}{2\eta_f} \langle \widehat f (\Bs, \beta; \theta) \rangle_+^2
\label{gp_14}
\end{equation*}
in terms of the threshold function defined in \req{gp_9}. Then, the incremental  internal potential density \req{inc_vp_threshold_2} takes the form
\begin{equation}
\begin{split}
\widehat \pi^{*\tau} = \widehat \psi + \tau [\,(\theta_k - T) \dot \eta^\tau +  &\eta_n \dot \theta^\tau + \frac{T}{\theta_n} \,  ( \Bs : \dot \Bvarepsilon^{p \tau} +\beta \dot \alpha^{\tau} )  \\ 
&- \frac{1}{2 \eta_f} \langle \widehat f(\Bs, \beta;\theta_n) \rangle_+^2   - \widehat \phi_{con}(-\frac{1}{T} \nabla T;\BC_n, \theta_n) \,] \, ,
\end{split}
\label{gd_18}
\end{equation}
and the incremental variational principle \req{extended_inc_3} yields \req{gp_12}$_{1-6}$  together with
\begin{equation*}
\parbox{14cm}{\centering 
$
\begin{array}{lll}
\partial_{\Bs} \widehat \pi^{*\tau} \equiv \frac{T}{\theta_n} (\Bvarepsilon^p - \Bvarepsilon^p_n) - (\tau / \eta_f) \, \langle  \widehat f (\Bs, \beta; \theta_n) \rangle_+ \,  \BG \Bs \BG / |  \Bs | &= \Bzero \,  , \\[1ex]
\partial_{\beta} \widehat \pi^{*\tau} \equiv \frac{T}{\theta_n} (\alpha - \alpha_n) - (\tau / \eta_f)\, \langle \widehat f (\Bs, \beta; \theta_n) \rangle_+ \, \sqrt{2/3}  &= 0 
\end{array}
$
}
\label{gp_19}
\end{equation*}
as Euler equations in $\calB$. This set of equations can again be reduced by expressing the dissipative driving forces $\Bs$ and $\beta$ from \req{gp_12}$_2$ and \req{gp_12}$_3$ yielding for $k=n$ the nonlocal (regularized) update equations
\begin{equation*}
\parbox{14cm}{\centering 
$
\begin{array}{lll}
\Bvarepsilon^p = \Bvarepsilon^p_n  - (\tau / \eta_f) \, \langle \,  | \, \Bbeta^e \, | - \sqrt{2/3} \, (\widehat y(\theta_n) + \beta^e) \, \rangle_+ \,  \BG \Bbeta^e \BG / | \Bbeta^e  | \,  , \\[1ex]
\alpha = \alpha_n  + (\tau / \eta_f) \,  \langle \,  | \, \Bbeta^e \, | - \sqrt{2/3} \, (\widehat y(\theta_n) + \beta^e) \, \rangle_+ \,\sqrt{2/3}
\end{array}
$
}
\label{gp_20}
\end{equation*}
for the plastic Hencky strain as well as the hardening variable. A local finite strain thermoplasticity model that uses the same additive strain kinematics together with the plastic configurational entropy in the sense of \citet{simo+miehe92} is proposed by \citet{ulz09}. For a comparison of rate-independent and rate-dependent formulations in isothermal gradient-plasticity of Fleck-Willis-type we refer to \citet{nielsen+niordson19}.

\inputtab[t]{examples/shearband/tab01}{tab01}%
%
\inputfig[t]{examples/shearband/bvp}{bvp}%
%
%
\inputfig[t]{examples/shearband/contour}{contour}%
%

\subsubsection[Section523]{Numerical Example: Cross Shear Localization} \label{Section523} 
For softening visco-plasticity of von Mises-type, we analyze the
development of shear bands in a rectangular plate $\calB = (0,L)
\times (0,2L)$ with $L=50$~mm subject to tensile loading under the
condition of plane strain. The geometric setup is depicted in
Figure~\ref{bvp}. We use the viscous over-force formulation of the
mixed large deformation setting from Section~\ref{Section531} above
and specify the isotropic hardening function in \req{gp_4}  as
\begin{equation*}
\widehat \psi_p(\alpha,\theta) = \frac{1}{2} \widehat h(\theta) \alpha^2 \, .
\label{cshear_1}
\end{equation*}
Here, $\widehat h$ is a temperature dependent hardening function which together with the temperature dependent yield stress function is  specified as
\begin{equation}
\widehat h (\theta) = h_0[ \, 1-w_h(\theta-\theta_0) \, ] 
\AND
\widehat y(\theta) = y_0[ \,1-w_0(\theta-\theta_0)  \,] \, ,
\label{cshear_2}
\end{equation}
see \citet{simo+miehe92}. The used material parameters are summarized
in Table~\ref{tab01}. To trigger plasticity in the center, the initial yield stress is reduced by $3$\% in the element shaded in dark grey in Figure~\ref{bvp}. For simplicity, we neglect the effect of heat conduction which is a reasonable assumption in case of a fast formation of the shear band generated by a high loading rate.  We stretch the specimen with a constant displacement rate $\dot{\bar u} = 5$~mm/s within the time interval $\calT = (0,1)$~s that is divided into 1000  equal increments\footnote{Since the process of heat conduction is neglected and the viscosity is chosen very low in order to have a formulation that is close to the nonsmooth setting, the overall thermomechanical material behavior is de facto rate-independent. Hence, the specific loading rate applied on the specimen is practically irrelevant.}. We use the  proposed semi-explicit variational update with index $k=n$
in the incremental internal potential density \req{gd_18}.  Since heat conduction is not taken into account, the time step size $\tau$ is not restricted by the CFL condition. Due to the
variational structure, the resulting stiffness matrix within a typical
Newton-Raphson iteration step is symmetric. As (global) primary fields
we take the macroscopic deformation $\Bvarphi$, the scalar
hardening/softening variable $\alpha$ and its dual driving force
$\beta$. The temperature $\theta$ is calculated via the implicit local equation \req{gp_12}$_5$. Due to symmetry, only one quarter of the domain is
discretized by $15 \times 30$, $20 \times 40$ and $25 \times 50$
quadrilateral finite elements. We use a Q$_1$E$_4$-Q$_1$-Q$_1$ element
pairing which bases on a (local) enhancement of the macroscopic
displacement gradient according to
\citet{simo+armero92}. Figure~\ref{contour} shows contour plots of the
equivalent plastic strain $\alpha$ and the temperature $\theta$ at
final displacement $\bar u = 5$~mm for three different plastic length
scale parameters $l \in \{0.05,0.1,0.2 \}$~mm. Clearly, the specimen
experiences a rise in the temperature in the region of plastic
dissipation. When increasing the plastic length scale, the equivalent plastic strain $\alpha$ as well as the temperature $\theta$ spread over more elements. At the same time, one observes decreasing maximum values of $\alpha$ and
$\theta$, see also
\citet{aldakheel+miehe17} and the references cited therein,
i.e., \citet{voyiadjis+faghihi12}. As widely known, in case of a local
theory  ($l=0$~mm) the plastic deformation localizes within one element
width. This mesh dependency also manifests itself in the
load-displacement curve of the structure as shown in
Figure~\ref{force-disp}a).  In contrast, the regularization provided by
the gradient theory yields mesh independent results, see Figure~\ref{force-disp}b). There, one also
observes the additional softening effect in the nonisothermal case due
to locally decreasing yield stresses according to \req{cshear_2}$_2$. At this point, we want to note that \citet{wcislo+pamin17} incorporate a gradient-enhanced relative temperature field in their formulation in order to regularize adiabatic thermal softening behavior.

For the used mixed setting of gradient plasticity, alternative finite element formulations are presented in \citet{miehe+welschinger+aldakheel14}. Note that in this context the construction of inf-sup stable finite element pairings is a difficult task which in detail has recently been investigated by \citet{krischok+linder19}. Especially, our chosen element pairing results in a nonphysical oscillatory behavior of the driving force field $\beta$. However, this instability seems to have no visible influence on the macroscopic deformation field $\Bvarphi$, the scalar plastic strain field $\alpha$ and the temperature field $\theta$. 
\inputfig[t]{examples/shearband/force-disp}{force-disp}%
%

\section[Section7]{Conclusion} \label{Section7}
We presented a unified framework for the fully coupled thermomechanics
of gradient-extended dissipative solids which undergo large
deformations. The key of the formulation is the consideration of the
entropy and the entropy rate as canonical variables which enter besides the gradient-extended mechanical state and the rate of this state,
respectively, the internal energy and dissipation functions. Here, the
rate-type formulation of local thermoplasticity outlined in
\citet{yang+stainier+ortiz06} is recovered by the concept of dual
variables. Especially, the quantity conjugate to the entropy is
rigorously distinguished from the quantity conjugate to the entropy
rate. The coupled macro- and micro-balances as well as the energy
equation are outcomes of the stated saddle point principle. Within
this setting, the entropy is also driven by additional gradient-type
dissipative effects.  Emphasis was also put on extended variational
formulations which incorporate dissipative mechanical driving forces
and temperature dependent threshold mechanisms. In addition, we
discussed variationally consistent time incrementations and suggested
a  semi-explicit numerical algorithm that renders the  algorithmically consistent form of the
intrinsic dissipation. It was shown that this algorithm has the structure of an operator split. A special focus was put on applying the framework to complex multifield problems which are of interest in the thermomechanics of solids. Three model problems, i.e., Cahn-Hilliard
diffusion, gradient damage and (additive) gradient plasticity strongly
coupled with temperature evolution, showed the capability of the
proposed formulation. In a numerical example we studied the formation
of a cross shear band under adiabatic condition. 
\\
\\

\textbf{Acknowledgement.} Support for this research was provided by the German Research Foundation (DFG) under Grant MI 295/20--1.

\bibliographystyle{elsarticle-harv} 
\bibliography{bibliography}

\begin{thebibliography}{90}
\expandafter\ifx\csname natexlab\endcsname\relax\def\natexlab#1{#1}\fi
\providecommand{\url}[1]{\texttt{#1}}
\providecommand{\href}[2]{#2}
\providecommand{\path}[1]{#1}
\providecommand{\DOIprefix}{doi:}
\providecommand{\ArXivprefix}{arXiv:}
\providecommand{\URLprefix}{URL: }
\providecommand{\Pubmedprefix}{pmid:}
\providecommand{\doi}[1]{\href{http://dx.doi.org/#1}{\path{#1}}}
\providecommand{\Pubmed}[1]{\href{pmid:#1}{\path{#1}}}
\providecommand{\bibinfo}[2]{#2}
\ifx\xfnm\relax \def\xfnm[#1]{\unskip,\space#1}\fi
\bibitem[{Aifantis(1987)}]{aifantis87}
\bibinfo{author}{Aifantis, E.C.}, \bibinfo{year}{1987}.
\newblock \bibinfo{title}{The physics of plastic deformation}.
\newblock \bibinfo{journal}{International Journal of Plasticity}
  \bibinfo{volume}{3}, \bibinfo{pages}{211--247}.
\bibitem[{Aldakheel \& Miehe(2017)}]{aldakheel+miehe17}
\bibinfo{author}{Aldakheel, F.}, \bibinfo{author}{Miehe, C.},
  \bibinfo{year}{2017}.
\newblock \bibinfo{title}{Coupled thermomechanical response of gradient
  plasticity}.
\newblock \bibinfo{journal}{International Journal of Plasticity}
  \bibinfo{volume}{91}, \bibinfo{pages}{1--24}.
\bibitem[{Armero \& Simo(1992)}]{armero+simo92}
\bibinfo{author}{Armero, F.}, \bibinfo{author}{Simo, J.C.},
  \bibinfo{year}{1992}.
\newblock \bibinfo{title}{A new unconditionally stable fractional step method
  for non-linear coupled thermomechanical problems}.
\newblock \bibinfo{journal}{International Journal for Numerical Methods in
  Engineering} \bibinfo{volume}{35}, \bibinfo{pages}{737--766}.
\bibitem[{Ball(1976)}]{ball76}
\bibinfo{author}{Ball, J.M.}, \bibinfo{year}{1976}.
\newblock \bibinfo{title}{Convexity conditions and existence theorems in
  nonlinear elasticity}.
\newblock \bibinfo{journal}{Archive for Rational Mechanics and Analysis}
  \bibinfo{volume}{63}, \bibinfo{pages}{337--403}.
\bibitem[{Bartels et~al.(2015)Bartels, Bartel, Canadija \&
  Mosler}]{bartels+etal15}
\bibinfo{author}{Bartels, A.}, \bibinfo{author}{Bartel, T.},
  \bibinfo{author}{Canadija, M.}, \bibinfo{author}{Mosler, J.},
  \bibinfo{year}{2015}.
\newblock \bibinfo{title}{On the thermomechanical coupling in dissipative
  materials: {A} variational approach for generalized standard materials}.
\newblock \bibinfo{journal}{Journal of the Mechanics and Physics of Solids}
  \bibinfo{volume}{82}, \bibinfo{pages}{218--234}.
\bibitem[{Berdichevsky(2009)}]{berdichevsky09}
\bibinfo{author}{Berdichevsky, V.L.}, \bibinfo{year}{2009}.
\newblock \bibinfo{title}{Variational Principles of Continuum Mechanics: {I}.
  Fundamentals}.
\newblock \bibinfo{publisher}{Springer}.
\bibitem[{Biot(1965)}]{biot65}
\bibinfo{author}{Biot, M.A.}, \bibinfo{year}{1965}.
\newblock \bibinfo{title}{Mechanics of Incremental Deformations}.
\newblock \bibinfo{publisher}{John Wiley \& Sons, Inc., New York}.
\bibitem[{Bruhns et~al.(2001)Bruhns, Xiao \& Meyers}]{bruhns+etal01}
\bibinfo{author}{Bruhns, O.T.}, \bibinfo{author}{Xiao, H.},
  \bibinfo{author}{Meyers, A.}, \bibinfo{year}{2001}.
\newblock \bibinfo{title}{Constitutive inequalities for an isotropic elastic
  strain-energy function based on {H}encky's logarithmic strain tensor}.
\newblock \bibinfo{journal}{Proceedings of the Royal Society of London. Series
  A: Mathematical, Physical and Engineering Sciences} \bibinfo{volume}{457},
  \bibinfo{pages}{2207--2226}.
\bibitem[{Cahn \& Hilliard(1958)}]{cahn+hilliard58}
\bibinfo{author}{Cahn, J.W.}, \bibinfo{author}{Hilliard, J.E.},
  \bibinfo{year}{1958}.
\newblock \bibinfo{title}{Free energy of a nonuniform system. {I}.
  {I}nterfacial free energy}.
\newblock \bibinfo{journal}{The Journal of chemical physics}
  \bibinfo{volume}{28}, \bibinfo{pages}{258--267}.
\bibitem[{Canadija \& Mosler(2011)}]{canadija+mosler11}
\bibinfo{author}{Canadija, M.}, \bibinfo{author}{Mosler, J.},
  \bibinfo{year}{2011}.
\newblock \bibinfo{title}{On the thermomechanical coupling in finite strain
  plasticity theory with non-linear kinematic hardening by means of incremental
  energy minimization}.
\newblock \bibinfo{journal}{International Journal of Solids and Structures}
  \bibinfo{volume}{48}, \bibinfo{pages}{1120--1129}.
\bibitem[{Capriz(2013)}]{capriz13}
\bibinfo{author}{Capriz, G.}, \bibinfo{year}{2013}.
\newblock \bibinfo{title}{Continua with microstructure}.
  volume~\bibinfo{volume}{35}.
\newblock \bibinfo{publisher}{Springer Science \& Business Media}.
\bibitem[{Carstensen et~al.(2002)Carstensen, Hackl \&
  Mielke}]{carstensen+hackl+mielke02}
\bibinfo{author}{Carstensen, C.}, \bibinfo{author}{Hackl, K.},
  \bibinfo{author}{Mielke, A.}, \bibinfo{year}{2002}.
\newblock \bibinfo{title}{Nonconvex potentials and microstructures in
  finite-strain plasticity}.
\newblock \bibinfo{journal}{Proceedings of the the Royal Society of London,
  Series A} \bibinfo{volume}{458}, \bibinfo{pages}{299--317}.
\bibitem[{Cosserat \& Cosserat(1896)}]{cosserat+cosserat1896}
\bibinfo{author}{Cosserat, E.}, \bibinfo{author}{Cosserat, F.},
  \bibinfo{year}{1896}.
\newblock \bibinfo{title}{Sur la th{\'e}orie de l'{\'e}lasticit{\'e}. {P}remier
  m{\'e}moire}, in: \bibinfo{booktitle}{Annales de la Facult{\'e} des sciences
  de Toulouse: Math{\'e}matiques}, pp. \bibinfo{pages}{I1--I116}.
\bibitem[{Dahlberg \& Ortiz(2019)}]{dahlberg+ortiz19}
\bibinfo{author}{Dahlberg, C.F.O.}, \bibinfo{author}{Ortiz, M.},
  \bibinfo{year}{2019}.
\newblock \bibinfo{title}{Fractional strain-gradient plasticity}.
\newblock \bibinfo{journal}{European Journal of Mechanics-A/Solids}
  \bibinfo{volume}{75}, \bibinfo{pages}{348--354}.
\bibitem[{De~Borst \& M{\"u}hlhaus(1992)}]{deborst+muehlhaus92}
\bibinfo{author}{De~Borst, R.}, \bibinfo{author}{M{\"u}hlhaus, H.B.},
  \bibinfo{year}{1992}.
\newblock \bibinfo{title}{Gradient-dependent plasticity: formulation and
  algorithmic aspects}.
\newblock \bibinfo{journal}{International Journal for Numerical Methods in
  Engineering} \bibinfo{volume}{35}, \bibinfo{pages}{521--539}.
\bibitem[{De~Groot \& Mazur(2013)}]{groot+mazur13}
\bibinfo{author}{De~Groot, S.R.}, \bibinfo{author}{Mazur, P.},
  \bibinfo{year}{2013}.
\newblock \bibinfo{title}{Non-equilibrium thermodynamics}.
\newblock \bibinfo{publisher}{Courier Corporation}.
\bibitem[{Dimitrijevic \& Hackl(2008)}]{dimitrijevic+hackl08}
\bibinfo{author}{Dimitrijevic, B.J.}, \bibinfo{author}{Hackl, K.},
  \bibinfo{year}{2008}.
\newblock \bibinfo{title}{A method for gradient enhancement of continuum damage
  models}.
\newblock \bibinfo{journal}{Technische Mechanik. Scientific Journal for
  Fundamentals and Applications of Engineering Mechanics} \bibinfo{volume}{28},
  \bibinfo{pages}{43--52}.
\bibitem[{Eringen(2012)}]{eringen12}
\bibinfo{author}{Eringen, A.C.}, \bibinfo{year}{2012}.
\newblock \bibinfo{title}{Microcontinuum field theories: I. Foundations and
  solids}.
\newblock \bibinfo{publisher}{Springer Science \& Business Media}.
\bibitem[{Fleck \& Hutchinson(2001)}]{fleck+hutchinson01}
\bibinfo{author}{Fleck, N.A.}, \bibinfo{author}{Hutchinson, J.W.},
  \bibinfo{year}{2001}.
\newblock \bibinfo{title}{A reformulation of strain gradient plasticity}.
\newblock \bibinfo{journal}{Journal of the Mechanics and Physics of Solids}
  \bibinfo{volume}{49}, \bibinfo{pages}{2245--2271}.
\bibitem[{Fleck \& Willis(2009)}]{fleck+willis09}
\bibinfo{author}{Fleck, N.A.}, \bibinfo{author}{Willis, J.R.},
  \bibinfo{year}{2009}.
\newblock \bibinfo{title}{A mathematical basis for strain-gradient plasticity
  theory. {P}art {I}: Scalar plastic multiplier}.
\newblock \bibinfo{journal}{Journal of the Mechanics and Physics of Solids}
  \bibinfo{volume}{57}, \bibinfo{pages}{161--177}.
\bibitem[{Fohrmeister et~al.(2018)Fohrmeister, Bartels \&
  Mosler}]{fohrmeister+bartles+mosler2018}
\bibinfo{author}{Fohrmeister, V.}, \bibinfo{author}{Bartels, A.},
  \bibinfo{author}{Mosler, J.}, \bibinfo{year}{2018}.
\newblock \bibinfo{title}{Variational updates for thermomechanically coupled
  gradient-enhanced elastoplasticity. {I}mplementation based on hyper-dual
  numbers}.
\newblock \bibinfo{journal}{Computer Methods in Applied Mechanics and
  Engineering} \bibinfo{volume}{339}, \bibinfo{pages}{239--261}.
\bibitem[{Forest(2009)}]{forest09}
\bibinfo{author}{Forest, S.}, \bibinfo{year}{2009}.
\newblock \bibinfo{title}{Micromorphic approach for gradient elasticity,
  viscoplasticity, and damage}.
\newblock \bibinfo{journal}{Journal of Engineering Mechanics}
  \bibinfo{volume}{135}, \bibinfo{pages}{117--131}.
\bibitem[{Forest \& Sievert(2003)}]{forest+sievert03}
\bibinfo{author}{Forest, S.}, \bibinfo{author}{Sievert, R.},
  \bibinfo{year}{2003}.
\newblock \bibinfo{title}{Elastoviscoplastic constitutive frameworks for
  generalized continua}.
\newblock \bibinfo{journal}{Acta Mechanica} \bibinfo{volume}{160},
  \bibinfo{pages}{71--111}.
\bibitem[{Francfort \& Mielke(2006)}]{francfort+mielke06}
\bibinfo{author}{Francfort, G.}, \bibinfo{author}{Mielke, A.},
  \bibinfo{year}{2006}.
\newblock \bibinfo{title}{Existence results for a class of rate-independent
  material models with nonconvex elastic energies}.
\newblock \bibinfo{journal}{Journal f{\"u}r die reine und angewandte Mathematik
  (Crelles Journal)} \bibinfo{volume}{2006}, \bibinfo{pages}{55--91}.
\bibitem[{Fr{\'e}mond(2013)}]{fremond13}
\bibinfo{author}{Fr{\'e}mond, M.}, \bibinfo{year}{2013}.
\newblock \bibinfo{title}{{N}on-smooth thermomechanics}.
\newblock \bibinfo{publisher}{Springer Science \& Business Media}.
\bibitem[{Gurtin(1996)}]{gurtin94}
\bibinfo{author}{Gurtin, M.E.}, \bibinfo{year}{1996}.
\newblock \bibinfo{title}{Generalized {G}inzburg-{L}andau and {C}ahn-{H}illiard
  equations based on a microforce balance}.
\newblock \bibinfo{journal}{Physica D} \bibinfo{volume}{92},
  \bibinfo{pages}{178--192}.
\bibitem[{Gurtin(2003)}]{gurtin03}
\bibinfo{author}{Gurtin, M.E.}, \bibinfo{year}{2003}.
\newblock \bibinfo{title}{On a framework for small-deformation viscoplasticity:
  free energy, microforces, strain gradients}.
\newblock \bibinfo{journal}{International Journal of Plasticity}
  \bibinfo{volume}{19}, \bibinfo{pages}{47--90}.
\bibitem[{Gurtin \& Anand(2005)}]{gurtin+anand05}
\bibinfo{author}{Gurtin, M.E.}, \bibinfo{author}{Anand, L.},
  \bibinfo{year}{2005}.
\newblock \bibinfo{title}{A theory of strain-gradient plasticity for isotropic,
  plastically irrotational materials. {P}art {I}: Small deformations}.
\newblock \bibinfo{journal}{Journal of the Mechanics and Physics of Solids}
  \bibinfo{volume}{53}, \bibinfo{pages}{1624--1649}.
\bibitem[{Hackl(1997)}]{hackl97}
\bibinfo{author}{Hackl, K.}, \bibinfo{year}{1997}.
\newblock \bibinfo{title}{Generalized standard media and variational principles
  in classical finite strain elastoplasticity}.
\newblock \bibinfo{journal}{Journal of the Mechanics and Physics of Solids}
  \bibinfo{volume}{45}, \bibinfo{pages}{667--688}.
\bibitem[{Hackl \& Fischer(2007)}]{hackl+fischer07}
\bibinfo{author}{Hackl, K.}, \bibinfo{author}{Fischer, F.D.},
  \bibinfo{year}{2007}.
\newblock \bibinfo{title}{On the relation between the principle of maximum
  dissipation and inelastic evolution given by dissipation potentials}.
\newblock \bibinfo{journal}{Proceedings of the Royal Society A: Mathematical,
  Physical and Engineering Sciences} \bibinfo{volume}{464},
  \bibinfo{pages}{117--132}.
\bibitem[{Hackl et~al.(2011a)Hackl, Fischer \&
  Svoboda}]{hackl+fischer+svoboda_ad11}
\bibinfo{author}{Hackl, K.}, \bibinfo{author}{Fischer, F.D.},
  \bibinfo{author}{Svoboda, J.}, \bibinfo{year}{2011}a.
\newblock \bibinfo{title}{Addendum. {A} study on the principle of maximum
  dissipation for coupled and non-coupled non-isothermal processes in
  materials}.
\newblock \bibinfo{journal}{Proceedings of the Royal Society A: Mathematical,
  Physical and Engineering Sciences} \bibinfo{volume}{467},
  \bibinfo{pages}{2422--2426}.
\bibitem[{Hackl et~al.(2011b)Hackl, Fischer \&
  Svoboda}]{hackl+fischer+svoboda10}
\bibinfo{author}{Hackl, K.}, \bibinfo{author}{Fischer, F.D.},
  \bibinfo{author}{Svoboda, J.}, \bibinfo{year}{2011}b.
\newblock \bibinfo{title}{A study on the principle of maximum dissipation for
  coupled and non-coupled non-isothermal processes in materials}.
\newblock \bibinfo{journal}{Proceedings of the Royal Society A: Mathematical,
  Physical and Engineering Sciences} \bibinfo{volume}{467},
  \bibinfo{pages}{1186--1196}.
\bibitem[{Halphen \& Nguyen(1975)}]{halphen+nguyen75}
\bibinfo{author}{Halphen, B.}, \bibinfo{author}{Nguyen, Q.S.},
  \bibinfo{year}{1975}.
\newblock \bibinfo{title}{Sur les mat\'eraux standards g\'en\'eralis\'es}.
\newblock \bibinfo{journal}{Journal de M\'ecanique} \bibinfo{volume}{40},
  \bibinfo{pages}{39--63}.
\bibitem[{Hill(1950)}]{hill50}
\bibinfo{author}{Hill, R.}, \bibinfo{year}{1950}.
\newblock \bibinfo{title}{The Mathematical Theory of Plasticity}.
\newblock \bibinfo{publisher}{Oxford University Press}.
\bibitem[{Krischok \& Linder(2019)}]{krischok+linder19}
\bibinfo{author}{Krischok, A.}, \bibinfo{author}{Linder, C.},
  \bibinfo{year}{2019}.
\newblock \bibinfo{title}{A generalized inf--sup test for multi-field
  mixed-variational methods}.
\newblock \bibinfo{journal}{Computer Methods in Applied Mechanics and
  Engineering} \bibinfo{volume}{357}, \bibinfo{pages}{112497}.
\bibitem[{Lancioni et~al.(2015)Lancioni, Yal{\c{c}}inkaya \&
  Cocks}]{lancioni+etal15}
\bibinfo{author}{Lancioni, G.}, \bibinfo{author}{Yal{\c{c}}inkaya, T.},
  \bibinfo{author}{Cocks, A.}, \bibinfo{year}{2015}.
\newblock \bibinfo{title}{Energy-based non-local plasticity models for
  deformation patterning, localization and fracture}.
\newblock \bibinfo{journal}{Proceedings of the Royal Society A: Mathematical,
  Physical and Engineering Sciences} \bibinfo{volume}{471},
  \bibinfo{pages}{20150275}.
\bibitem[{Le et~al.(2018)Le, Marigo, Maurini \& Vidoli}]{le+etal18}
\bibinfo{author}{Le, D.T.}, \bibinfo{author}{Marigo, J.J.},
  \bibinfo{author}{Maurini, C.}, \bibinfo{author}{Vidoli, S.},
  \bibinfo{year}{2018}.
\newblock \bibinfo{title}{Strain-gradient vs damage-gradient regularizations of
  softening damage models}.
\newblock \bibinfo{journal}{Computer Methods in Applied Mechanics and
  Engineering} \bibinfo{volume}{340}, \bibinfo{pages}{424--450}.
\bibitem[{Leismann \& Mahnken(2015)}]{leismann+mahnken15}
\bibinfo{author}{Leismann, T.}, \bibinfo{author}{Mahnken, R.},
  \bibinfo{year}{2015}.
\newblock \bibinfo{title}{Comparison of hyperelastic micromorphic, micropolar
  and microstrain continua}.
\newblock \bibinfo{journal}{International Journal of Non-Linear Mechanics}
  \bibinfo{volume}{77}, \bibinfo{pages}{115--127}.
\bibitem[{Lemaitre(1996)}]{lemaitre96}
\bibinfo{author}{Lemaitre, J.}, \bibinfo{year}{1996}.
\newblock \bibinfo{title}{A Course on Damage Mechanics}.
\newblock \bibinfo{publisher}{Springer}.
\bibitem[{Lorentz \& Andrieux(1999)}]{lorentz+andrieux99}
\bibinfo{author}{Lorentz, E.}, \bibinfo{author}{Andrieux, S.},
  \bibinfo{year}{1999}.
\newblock \bibinfo{title}{A variational formulation for nonlocal damage
  models}.
\newblock \bibinfo{journal}{International Journal of Plasticity}
  \bibinfo{volume}{15}, \bibinfo{pages}{119--138}.
\bibitem[{Lu \& Pister(1975)}]{lu+pister75}
\bibinfo{author}{Lu, S.C.H.}, \bibinfo{author}{Pister, K.S.},
  \bibinfo{year}{1975}.
\newblock \bibinfo{title}{Decomposition of the deformation and representation
  of the free energy function for isotropic thermoelastic solids}.
\newblock \bibinfo{journal}{International Journal of Solids and Structures}
  \bibinfo{volume}{11}, \bibinfo{pages}{927--934}.
\bibitem[{Lubliner(2008)}]{lubliner08}
\bibinfo{author}{Lubliner, J.}, \bibinfo{year}{2008}.
\newblock \bibinfo{title}{Plasticity theory}.
\newblock \bibinfo{publisher}{Courier Corporation}.
\bibitem[{Mainik \& Mielke(2009)}]{mainik+mielke09}
\bibinfo{author}{Mainik, A.}, \bibinfo{author}{Mielke, A.},
  \bibinfo{year}{2009}.
\newblock \bibinfo{title}{Global existence for rate-independent gradient
  plasticity at finite strain}.
\newblock \bibinfo{journal}{Journal of Nonlinear Science} \bibinfo{volume}{19},
  \bibinfo{pages}{221--248}.
\bibitem[{Mariano(2002)}]{mariano02}
\bibinfo{author}{Mariano, P.M.}, \bibinfo{year}{2002}.
\newblock \bibinfo{title}{Multifield theories in mechanics of solids}, in:
  \bibinfo{booktitle}{Advances in Applied Mechanics}.
  \bibinfo{publisher}{Elsevier}. volume~\bibinfo{volume}{38}, pp.
  \bibinfo{pages}{1--93}.
\bibitem[{Maugin(1990)}]{maugin90}
\bibinfo{author}{Maugin, G.A.}, \bibinfo{year}{1990}.
\newblock \bibinfo{title}{Internal variables and dissipative structures}.
\newblock \bibinfo{journal}{Journal of Non-Equilibrium Thermodynamics}
  \bibinfo{volume}{15}, \bibinfo{pages}{173--192}.
\bibitem[{Maugin(1992)}]{maugin92}
\bibinfo{author}{Maugin, G.A.}, \bibinfo{year}{1992}.
\newblock \bibinfo{title}{The Thermomechanics of Plasticity and Fracture}.
\newblock \bibinfo{publisher}{Cambridge University Press}.
\bibitem[{Maugin \& Muschik(1994)}]{maugin+muschik94}
\bibinfo{author}{Maugin, G.A.}, \bibinfo{author}{Muschik, W.},
  \bibinfo{year}{1994}.
\newblock \bibinfo{title}{Thermodynamics with internal variables. {P}art {I}.
  general concepts}.
\newblock \bibinfo{journal}{Journal of Non-Equilibrium Thermodynamics}
  \bibinfo{volume}{19}, \bibinfo{pages}{217--249}.
\bibitem[{Miehe(2002)}]{miehe02}
\bibinfo{author}{Miehe, C.}, \bibinfo{year}{2002}.
\newblock \bibinfo{title}{Strain-driven homogenization of inelastic
  microstructures and composites based on an incremental variational
  formulation}.
\newblock \bibinfo{journal}{International Journal for Numerical Methods in
  Engineering} \bibinfo{volume}{55}, \bibinfo{pages}{1285--1322}.
\bibitem[{Miehe(2011)}]{miehe11}
\bibinfo{author}{Miehe, C.}, \bibinfo{year}{2011}.
\newblock \bibinfo{title}{A multi-field incremental variational framework for
  gradient-extended standard dissipative solids}.
\newblock \bibinfo{journal}{Journal of the Mechanics and Physics of Solids}
  \bibinfo{volume}{59}, \bibinfo{pages}{898--923}.
\bibitem[{Miehe(2014)}]{miehe14}
\bibinfo{author}{Miehe, C.}, \bibinfo{year}{2014}.
\newblock \bibinfo{title}{Variational gradient plasticity at finite strains.
  {P}art {I}: Mixed potentials for the evolution and update problems of
  gradient-extended dissipative solids}.
\newblock \bibinfo{journal}{Computer Methods in Applied Mechanics and
  Engineering} \bibinfo{volume}{268}, \bibinfo{pages}{677--703}.
\bibitem[{Miehe et~al.(2002)Miehe, Apel \& Lambrecht}]{miehe+apel+lambrecht02}
\bibinfo{author}{Miehe, C.}, \bibinfo{author}{Apel, N.},
  \bibinfo{author}{Lambrecht, M.}, \bibinfo{year}{2002}.
\newblock \bibinfo{title}{Anisotropic additive plasticity in the logarithmic
  strain space. {M}odular kinematic formulation and implementation based on
  incremental minimization principles for standard materials}.
\newblock \bibinfo{journal}{Computer Methods in Applied Mechanics and
  Engineering} \bibinfo{volume}{191}, \bibinfo{pages}{5383--5425}.
\bibitem[{Miehe et~al.(2014a)Miehe, Hildebrand \&
  B\"oger}]{miehe+hildebrand+boeger14}
\bibinfo{author}{Miehe, C.}, \bibinfo{author}{Hildebrand, F.E.},
  \bibinfo{author}{B\"oger, L.}, \bibinfo{year}{2014}a.
\newblock \bibinfo{title}{Mixed variational potentials and inherent symmetries
  of the {C}ahn-{H}illiard theory of diffusive phase separation}.
\newblock \bibinfo{journal}{Proceedings of the Royal Society A: Mathematical,
  Physical and Engineering Science} \bibinfo{volume}{470},
  \bibinfo{pages}{20130641 1--16}.
\bibitem[{Miehe et~al.(2015)Miehe, Mauthe \&
  Teichtmeister}]{miehe+mauthe+teichtmeister15}
\bibinfo{author}{Miehe, C.}, \bibinfo{author}{Mauthe, S.},
  \bibinfo{author}{Teichtmeister, S.}, \bibinfo{year}{2015}.
\newblock \bibinfo{title}{Minimization principles for the coupled problem of
  {D}arcy-{B}iot-type fluid transport in porous media linked to phase field
  modeling of fracture}.
\newblock \bibinfo{journal}{Journal of the Mechanics and Physics of Solids}
  \bibinfo{volume}{82}, \bibinfo{pages}{186--217}.
\bibitem[{Miehe et~al.(2014b)Miehe, Welschinger \&
  Aldakheel}]{miehe+welschinger+aldakheel14}
\bibinfo{author}{Miehe, C.}, \bibinfo{author}{Welschinger, F.},
  \bibinfo{author}{Aldakheel, F.}, \bibinfo{year}{2014}b.
\newblock \bibinfo{title}{Variational gradient plasticity at finite strains.
  {P}art {II}: Local-global updates and mixed finite elements for additive
  plasticity in the logarithmic strain space}.
\newblock \bibinfo{journal}{Computer Methods in Applied Mechanics and
  Engineering} \bibinfo{volume}{268}, \bibinfo{pages}{704--734}.
\bibitem[{Mielke(2016)}]{mielke16}
\bibinfo{author}{Mielke, A.}, \bibinfo{year}{2016}.
\newblock \bibinfo{title}{Free energy, free entropy, and a gradient structure
  for thermoplasticity}, in: \bibinfo{booktitle}{Innovative Numerical
  Approaches for Multi-Field and Multi-Scale Problems}.
  \bibinfo{publisher}{Springer}, pp. \bibinfo{pages}{135--160}.
\bibitem[{Mielke \& Roub{\'\i}{\v{c}}ek(2006)}]{mielke+roubicek06}
\bibinfo{author}{Mielke, A.}, \bibinfo{author}{Roub{\'\i}{\v{c}}ek, T.},
  \bibinfo{year}{2006}.
\newblock \bibinfo{title}{Rate-independent damage processes in nonlinear
  elasticity}.
\newblock \bibinfo{journal}{Mathematical Models and Methods in Applied
  Sciences} \bibinfo{volume}{16}, \bibinfo{pages}{177--209}.
\bibitem[{Moreau(1976)}]{moreau76}
\bibinfo{author}{Moreau, J.J.}, \bibinfo{year}{1976}.
\newblock \bibinfo{title}{Application of convex analysis to the treatment of
  elastoplastic systems}, in: \bibinfo{booktitle}{Applications of methods of
  functional analysis to problems in mechanics}. \bibinfo{publisher}{Springer},
  pp. \bibinfo{pages}{56--89}.
\bibitem[{Mosler \& Bruhns(2010)}]{mosler+bruhns10}
\bibinfo{author}{Mosler, J.}, \bibinfo{author}{Bruhns, O.T.},
  \bibinfo{year}{2010}.
\newblock \bibinfo{title}{On the implementation of rate-independent standard
  dissipative solids at finite strain--{V}ariational constitutive updates}.
\newblock \bibinfo{journal}{Computer Methods in Applied Mechanics and
  Engineering} \bibinfo{volume}{199}, \bibinfo{pages}{417--429}.
\bibitem[{M{\"u}hlhaus \& Aifantis(1991)}]{muhlhaus+aifantis91}
\bibinfo{author}{M{\"u}hlhaus, H.B.}, \bibinfo{author}{Aifantis, E.C.},
  \bibinfo{year}{1991}.
\newblock \bibinfo{title}{A variational principle for gradient plasticity}.
\newblock \bibinfo{journal}{International Journal of Solids and Structures}
  \bibinfo{volume}{28}, \bibinfo{pages}{845--857}.
\bibitem[{Nateghi \& Keip(2021)}]{nateghi+keip21}
\bibinfo{author}{Nateghi, A.}, \bibinfo{author}{Keip, M.A.},
  \bibinfo{year}{2021}.
\newblock \bibinfo{title}{A thermo-chemo-mechanically coupled model for cathode
  particles in lithium--ion batteries}.
\newblock \bibinfo{journal}{Acta Mechanica} \bibinfo{volume}{232},
  \bibinfo{pages}{3041--3065}.
\bibitem[{Neff et~al.(2015)Neff, Ghiba \& Lankeit}]{neff+ghiba+lankeit15}
\bibinfo{author}{Neff, P.}, \bibinfo{author}{Ghiba, I.D.},
  \bibinfo{author}{Lankeit, J.}, \bibinfo{year}{2015}.
\newblock \bibinfo{title}{The exponentiated {H}encky-logarithmic strain energy.
  {P}art {I}: Constitutive issues and rank-one convexity}.
\newblock \bibinfo{journal}{Journal of Elasticity} \bibinfo{volume}{121},
  \bibinfo{pages}{143--234}.
\bibitem[{Nguyen(2000)}]{nguyen00}
\bibinfo{author}{Nguyen, Q.S.}, \bibinfo{year}{2000}.
\newblock \bibinfo{title}{Stability and Nonlinear Solid Mechanics}.
\newblock \bibinfo{publisher}{John Wiley \& Sons, New York}.
\bibitem[{Nguyen(2011)}]{nguyen11}
\bibinfo{author}{Nguyen, Q.S.}, \bibinfo{year}{2011}.
\newblock \bibinfo{title}{Variational principles in the theory of gradient
  plasticity}.
\newblock \bibinfo{journal}{Comptes Rendus Mecanique} \bibinfo{volume}{339},
  \bibinfo{pages}{743--750}.
\bibitem[{Nielsen \& Niordson(2019)}]{nielsen+niordson19}
\bibinfo{author}{Nielsen, K.}, \bibinfo{author}{Niordson, C.},
  \bibinfo{year}{2019}.
\newblock \bibinfo{title}{A finite strain {F}{E}-implementation of the
  {F}leck-{W}illis gradient theory: {R}ate-independent versus visco-plastic
  formulation}.
\newblock \bibinfo{journal}{European Journal of Mechanics-A/Solids}
  \bibinfo{volume}{75}, \bibinfo{pages}{389--398}.
\bibitem[{Oliferuk \& Raniecki(2018)}]{oliferuk+raniecki18}
\bibinfo{author}{Oliferuk, W.}, \bibinfo{author}{Raniecki, B.},
  \bibinfo{year}{2018}.
\newblock \bibinfo{title}{Thermodynamic description of the plastic work
  partition into stored energy and heat during deformation of polycrystalline
  materials}.
\newblock \bibinfo{journal}{European Journal of Mechanics-A/Solids}
  \bibinfo{volume}{71}, \bibinfo{pages}{326--334}.
\bibitem[{Ortiz \& Repetto(1999)}]{ortiz+repetto99}
\bibinfo{author}{Ortiz, M.}, \bibinfo{author}{Repetto, E.A.},
  \bibinfo{year}{1999}.
\newblock \bibinfo{title}{Nonconvex energy minimization and dislocation
  structures in ductile single crystals}.
\newblock \bibinfo{journal}{Journal of the Mechanics and Physics of Solids}
  \bibinfo{volume}{47}, \bibinfo{pages}{397--462}.
\bibitem[{Ortiz \& Stainier(1999)}]{ortiz+stainier99}
\bibinfo{author}{Ortiz, M.}, \bibinfo{author}{Stainier, L.},
  \bibinfo{year}{1999}.
\newblock \bibinfo{title}{The variational formulation of viscoplastic
  constitutive updates}.
\newblock \bibinfo{journal}{Computer Methods in Applied Mechanics and
  Engineering} \bibinfo{volume}{171}, \bibinfo{pages}{419--444}.
\bibitem[{Peerlings et~al.(1996)Peerlings, de~Borst, Brekelmans \&
  De~Vree}]{peerlings+etal96}
\bibinfo{author}{Peerlings, R.H.J.}, \bibinfo{author}{de~Borst, R.},
  \bibinfo{author}{Brekelmans, W.A.M.}, \bibinfo{author}{De~Vree, J.H.P.},
  \bibinfo{year}{1996}.
\newblock \bibinfo{title}{Gradient enhanced damage for quasi-brittle
  materials}.
\newblock \bibinfo{journal}{International Journal for Numerical Methods in
  Engineering} \bibinfo{volume}{39}, \bibinfo{pages}{3391--3403}.
\bibitem[{Perzyna(1966)}]{perzyna66}
\bibinfo{author}{Perzyna, P.}, \bibinfo{year}{1966}.
\newblock \bibinfo{title}{Fundamental problems in viscoplasticity}.
\newblock \bibinfo{journal}{Advances in Applied Mechanics} \bibinfo{volume}{9},
  \bibinfo{pages}{243--377}.
\bibitem[{Petryk(2003)}]{petryk03}
\bibinfo{author}{Petryk, H.}, \bibinfo{year}{2003}.
\newblock \bibinfo{title}{Incremental energy minimization in dissipative
  solids}.
\newblock \bibinfo{journal}{Comptes Rendus Mecanique} \bibinfo{volume}{331},
  \bibinfo{pages}{469--474}.
\bibitem[{Pham et~al.(2011)Pham, Amor, Marigo \& Maurini}]{pham+etal11}
\bibinfo{author}{Pham, K.}, \bibinfo{author}{Amor, H.},
  \bibinfo{author}{Marigo, J.J.}, \bibinfo{author}{Maurini, C.},
  \bibinfo{year}{2011}.
\newblock \bibinfo{title}{Gradient damage models and their use to approximate
  brittle fracture}.
\newblock \bibinfo{journal}{International Journal of Damage Mechanics}
  \bibinfo{volume}{20}, \bibinfo{pages}{618--652}.
\bibitem[{Placidi \& Barchiesi(2018)}]{placidi+barchiesi18}
\bibinfo{author}{Placidi, L.}, \bibinfo{author}{Barchiesi, E.},
  \bibinfo{year}{2018}.
\newblock \bibinfo{title}{Energy approach to brittle fracture in
  strain-gradient modelling}.
\newblock \bibinfo{journal}{Proceedings of the Royal Society A: Mathematical,
  Physical and Engineering Sciences} \bibinfo{volume}{474},
  \bibinfo{pages}{20170878}.
\bibitem[{Placidi et~al.(2018)Placidi, Barchiesi \&
  Misra}]{placidi+barchiesi+misra18}
\bibinfo{author}{Placidi, L.}, \bibinfo{author}{Barchiesi, E.},
  \bibinfo{author}{Misra, A.}, \bibinfo{year}{2018}.
\newblock \bibinfo{title}{A strain gradient variational approach to damage: a
  comparison with damage gradient models and numerical results}.
\newblock \bibinfo{journal}{Mathematics and Mechanics of Complex Systems}
  \bibinfo{volume}{6}, \bibinfo{pages}{77--100}.
\bibitem[{Placidi et~al.(2019)Placidi, Misra \&
  Barchiesi}]{placidi+misra+barchiesi19}
\bibinfo{author}{Placidi, L.}, \bibinfo{author}{Misra, A.},
  \bibinfo{author}{Barchiesi, E.}, \bibinfo{year}{2019}.
\newblock \bibinfo{title}{Simulation results for damage with evolving
  microstructure and growing strain gradient moduli}.
\newblock \bibinfo{journal}{Continuum Mechanics and Thermodynamics}
  \bibinfo{volume}{31}, \bibinfo{pages}{1143--1163}.
\bibitem[{Rosakis et~al.(2000)Rosakis, Rosakis, Ravichandran \&
  Hodowany}]{rosakis+rosakis+ravichandran+hodowany00}
\bibinfo{author}{Rosakis, P.}, \bibinfo{author}{Rosakis, A.J.},
  \bibinfo{author}{Ravichandran, G.}, \bibinfo{author}{Hodowany, J.},
  \bibinfo{year}{2000}.
\newblock \bibinfo{title}{A thermodynamic internal variable model for the
  partition of plastic work into heat and stored energy in metals}.
\newblock \bibinfo{journal}{Journal of the Mechanics and Physics of Solids}
  \bibinfo{volume}{48}, \bibinfo{pages}{581--607}.
\bibitem[{Simo(1988)}]{simo88}
\bibinfo{author}{Simo, J.C.}, \bibinfo{year}{1988}.
\newblock \bibinfo{title}{A framework for finite strain elastoplasticity based
  on maximum plastic dissipation and the multiplicative decomposition: {P}art
  {I}. {C}ontinuum formulation}.
\newblock \bibinfo{journal}{Computer Methods in Applied Mechanics and
  Engineering} \bibinfo{volume}{66}, \bibinfo{pages}{199--219}.
\bibitem[{Simo \& Armero(1992)}]{simo+armero92}
\bibinfo{author}{Simo, J.C.}, \bibinfo{author}{Armero, F.},
  \bibinfo{year}{1992}.
\newblock \bibinfo{title}{Geometrically nonlinear enhanced strain mixed methods
  and the method of incompatible modes}.
\newblock \bibinfo{journal}{International Journal for Numerical Methods in
  Engineering} \bibinfo{volume}{33}, \bibinfo{pages}{1413--1449}.
\bibitem[{Simo \& Honein(1990)}]{simo+honein90}
\bibinfo{author}{Simo, J.C.}, \bibinfo{author}{Honein, T.},
  \bibinfo{year}{1990}.
\newblock \bibinfo{title}{Variational formulation, discrete conservation laws,
  and path--domain independent integrals for elasto--viscoplasticity}.
\newblock \bibinfo{journal}{Journal of Applied Mechanics} \bibinfo{volume}{57},
  \bibinfo{pages}{488--497}.
\bibitem[{Simo et~al.(1988)Simo, Kennedy \&
  Govindjee}]{simo+kennedy+govindjee88}
\bibinfo{author}{Simo, J.C.}, \bibinfo{author}{Kennedy, J.G.},
  \bibinfo{author}{Govindjee, S.}, \bibinfo{year}{1988}.
\newblock \bibinfo{title}{Non-smooth multisurface plasticity and
  viscoplasticity. {L}oading/unloading conditions and numerical algorithms}.
\newblock \bibinfo{journal}{International Journal for Numerical Methods in
  Engineering} \bibinfo{volume}{26}, \bibinfo{pages}{2161--2185}.
\bibitem[{Simo \& Miehe(1992)}]{simo+miehe92}
\bibinfo{author}{Simo, J.C.}, \bibinfo{author}{Miehe, C.},
  \bibinfo{year}{1992}.
\newblock \bibinfo{title}{Associative coupled thermoplasticity at finite
  strains: Formulation, numerical analysis and implementation}.
\newblock \bibinfo{journal}{Computer Methods in Applied Mechanics and
  Engineering} \bibinfo{volume}{98}, \bibinfo{pages}{41--104}.
\bibitem[{Stainier \& Ortiz(2010)}]{stainier+ortiz10}
\bibinfo{author}{Stainier, L.}, \bibinfo{author}{Ortiz, M.},
  \bibinfo{year}{2010}.
\newblock \bibinfo{title}{Study and validation of a variational theory of
  thermo-mechanical coupling in finite visco-plasticity}.
\newblock \bibinfo{journal}{Internation Journal of Solids and Structures}
  \bibinfo{volume}{47}, \bibinfo{pages}{705--715}.
\bibitem[{Su et~al.(2014)Su, Stainier \& Mercier}]{su+stainier+mercier14}
\bibinfo{author}{Su, S.}, \bibinfo{author}{Stainier, L.},
  \bibinfo{author}{Mercier, S.}, \bibinfo{year}{2014}.
\newblock \bibinfo{title}{Energy-based variational modeling of fully formed
  adiabatic shear bands}.
\newblock \bibinfo{journal}{European Journal of Mechanics-A/Solids}
  \bibinfo{volume}{47}, \bibinfo{pages}{1--13}.
\bibitem[{Svendsen(2004)}]{svendsen04}
\bibinfo{author}{Svendsen, B.}, \bibinfo{year}{2004}.
\newblock \bibinfo{title}{On thermodynamic- and variational-based formulations
  of models for inelastic continua with internal lengthscales}.
\newblock \bibinfo{journal}{Computer Methods in Applied Mechanics and
  Engineering} \bibinfo{volume}{193}, \bibinfo{pages}{5429--5452}.
\bibitem[{Taylor \& Quinney(1934)}]{taylor+quinney34}
\bibinfo{author}{Taylor, G.I.}, \bibinfo{author}{Quinney, H.},
  \bibinfo{year}{1934}.
\newblock \bibinfo{title}{The latent energy remaining in a metal after cold
  working}.
\newblock \bibinfo{journal}{Proceedings of the Royal Society A: Mathematical,
  Physical and Engineering Sciences} \bibinfo{volume}{143},
  \bibinfo{pages}{307--326}.
\bibitem[{Timofeev et~al.(2021)Timofeev, Barchiesi, Misra \&
  Placidi}]{timofeev+etal2021}
\bibinfo{author}{Timofeev, D.}, \bibinfo{author}{Barchiesi, E.},
  \bibinfo{author}{Misra, A.}, \bibinfo{author}{Placidi, L.},
  \bibinfo{year}{2021}.
\newblock \bibinfo{title}{Hemivariational continuum approach for granular
  solids with damage-induced anisotropy evolution}.
\newblock \bibinfo{journal}{Mathematics and Mechanics of Solids}
  \bibinfo{volume}{26}, \bibinfo{pages}{738--770}.
\bibitem[{Ulz(2009)}]{ulz09}
\bibinfo{author}{Ulz, M.H.}, \bibinfo{year}{2009}.
\newblock \bibinfo{title}{A {G}reen--{N}aghdi approach to finite anisotropic
  rate-independent and rate-dependent thermo-plasticity in logarithmic
  {L}agrangean strain--entropy space}.
\newblock \bibinfo{journal}{Computer Methods in Applied Mechanics and
  Engineering} \bibinfo{volume}{198}, \bibinfo{pages}{3262--3277}.
\bibitem[{Voyiadjis \& Faghihi(2012)}]{voyiadjis+faghihi12}
\bibinfo{author}{Voyiadjis, G.Z.}, \bibinfo{author}{Faghihi, D.},
  \bibinfo{year}{2012}.
\newblock \bibinfo{title}{Thermo-mechanical strain gradient plasticity with
  energetic and dissipative length scales}.
\newblock \bibinfo{journal}{International Journal of Plasticity}
  \bibinfo{volume}{30}, \bibinfo{pages}{218--247}.
\bibitem[{Wcis{\l}o \& Pamin(2017)}]{wcislo+pamin17}
\bibinfo{author}{Wcis{\l}o, B.}, \bibinfo{author}{Pamin, J.},
  \bibinfo{year}{2017}.
\newblock \bibinfo{title}{Local and non-local thermomechanical modeling of
  elastic-plastic materials undergoing large strains}.
\newblock \bibinfo{journal}{International Journal for Numerical Methods in
  Engineering} \bibinfo{volume}{109}, \bibinfo{pages}{102--124}.
\bibitem[{Yang et~al.(2006)Yang, Stainier \& Ortiz}]{yang+stainier+ortiz06}
\bibinfo{author}{Yang, Q.}, \bibinfo{author}{Stainier, L.},
  \bibinfo{author}{Ortiz, M.}, \bibinfo{year}{2006}.
\newblock \bibinfo{title}{A variational formulation of the coupled
  thermo-mechanical boundary-value problem for general dissipative solids}.
\newblock \bibinfo{journal}{Journal of the Mechanics and Physics of Solids}
  \bibinfo{volume}{54}, \bibinfo{pages}{401--424}.
\bibitem[{Ziegler(1963)}]{ziegler63}
\bibinfo{author}{Ziegler, H.}, \bibinfo{year}{1963}.
\newblock \bibinfo{title}{Some extremum principles in irreversible
  thermodynamics with application to continuum mechanics}, in:
  \bibinfo{editor}{Sneddon, I.N.}, \bibinfo{editor}{Hill, R.} (Eds.),
  \bibinfo{booktitle}{Progress in Solid Mechanics, Vol. IV}.
  \bibinfo{publisher}{North--Holland Publishing Company, Amsterdam}.

\end{thebibliography}

\end{document}